\documentclass[a4paper,10pt]{article}

\usepackage[lmargin=3.5cm,rmargin=3.5cm,tmargin=4cm,bmargin=4cm]{geometry}
\usepackage[T1]{fontenc}
\usepackage[bookmarksdepth=3]{hyperref}

\usepackage[tbtags]{amsmath}
\usepackage{amsfonts,amssymb}

\usepackage[standard,amsmath,thmmarks,hyperref]{ntheorem}
\usepackage{abstract}

\usepackage{xspace}
\usepackage{etex}
\usepackage[english]{babel}
\usepackage{csquotes}
\usepackage[nocompress]{cite}

\usepackage[titletoc]{appendix}

\usepackage{graphicx}
\usepackage[labelformat=simple]{subcaption}

\usepackage{tikz}
\usetikzlibrary{shapes,calc,decorations.pathreplacing,decorations.markings,arrows}
\usepackage{tikz-cd}

\usepackage{enumitem}
\renewlist{enumerate}{enumerate}{1}
\setlist[enumerate]{label=(\arabic*),}

\newcommand{\DBF}[1]{\expandafter\newcommand\csname #1\endcsname{{\mathbf{ #1 }}}}

\DBF{C}
\DBF{R}
\DBF{Q}
\DBF{Z}
\DBF{N}
\DBF{T}
\DBF{1}

\newcommand{\DCF}[1]{\expandafter\newcommand\csname cal#1\endcsname{{\mathcal{ #1 }}}}

\DCF{A}
\DCF{B}
\DCF{D}
\DCF{E}
\DCF{F}
\DCF{G}
\DCF{L}
\DCF{M}
\DCF{N}
\DCF{P}
\DCF{S}
\DCF{X}
\DCF{Y}
\DCF{Z}

\newcommand{\DFF}[1]{\expandafter\newcommand\csname frak#1\endcsname{{\mathfrak{ #1 }}}}

\DFF{L}

\newcommand {\DMO}[1]{\expandafter\DeclareMathOperator\csname #1\endcsname {#1}}

\DMO{supp}
\DMO{id}
\DeclareMathOperator{\img}{im}
\DeclareMathOperator{\vspan}{span}
\DMO{ind}
\DMO{sgn}
\DeclareMathOperator{\homology}{H}
\DMO{HC}
\DMO{HTW}
\DeclareMathOperator{\clop}{cl}
\DMO{Vol}
\DMO{Inv}
\DMO{BInv}
\DMO{SF}

\newcommand{\sphere}{\mathbf{S}}
\newcommand{\pt}{\mathrm{pt}}
\newcommand{\setmin}{\smallsetminus}

\newcommand{\iso}{\cong}

\newcommand{\set}[2]{ \left\{ #1 \; : \; #2 \right\} }
\newcommand\rst[2]{ \left. {#1} \right|_{#2} }

\newcommand{\cl}[1]{ \ensuremath{\overline{#1}} }
\newcommand{\clin}[2]{ {\clop}_{#2}{#1} }
\newcommand{\bdy}[1]{\partial {#1} }
\newcommand{\inter}[1]{\mathrm{int}\,{#1}}

\newcommand{\real}{\mathfrak{Re}\,}
\newcommand{\imag}{\mathfrak{Im}\,}

\renewcommand{\d}[1]{ \ensuremath{ \operatorname{d}\!{#1} } }
\newcommand{\dn}[2]{ \ensuremath{ \operatorname{d}^{#1}\!{#2} } }
\newcommand{\dsub}[2]{ \ensuremath{ \operatorname{d}_{#1}\!{#2} } }

\newcommand{\ddi}[2]{\ensuremath{ \dfrac{\partial #1 }{\partial #2} }}
\newcommand{\ddx}[2]{\ensuremath{ \dfrac{\d #1 }{\d #2} }}

\newcommand{\ddxe}[1]{\ensuremath{ \dfrac{\d{ \hspace{1ex} }}{\d{#1} } }}

\newcommand{\modn}[1]{\quad (\text{mod } {#1})}

\renewcommand{\epsilon}{\varepsilon}
\renewcommand{\phi}{\varphi}

\renewcommand{\leq}{\leqslant}
\renewcommand{\geq}{\geqslant}
\renewcommand{\dots}{\ldots}
\renewcommand{\qed}{\hfill \ensuremath{\Box}}

\makeatletter
\newtheoremstyle{mynonumberplain}%
  {\item[\theorem@headerfont\hskip\labelsep ##1\theorem@separator]}%
  {\item[\theorem@headerfont\hskip \labelsep ##1\
    ##3\theorem@separator]}
\makeatother

\theoremstyle{break}
\theoremprework{} 
\theorembodyfont{\slshape}
\theorempostwork{}
\theoremseparator{:}
\newcounter{dummy}
\numberwithin{dummy}{section}
\newtheorem{mydef}[dummy]{Definition}
\newtheorem{mythm}[dummy]{Theorem}
\renewtheorem*{mythm*}{Theorem}
\newtheorem{mylemma}[dummy]{Lemma}
\renewtheorem*{mylemma*}{Lemma}

\theoremnumbering{Alph}
\newtheorem{mythmintro}{Theorem}
\theoremnumbering{arabic}
\theoremsymbol{\ensuremath{\qed}}
\newtheorem{myremark}[dummy]{Remark}

\theoremstyle{nonumberplain}
\theorembodyfont{\upshape}
\theoremsymbol{\ensuremath{\diamondsuit}}
\theoremseparator{:}
\theoremprework{}

\theoremstyle{mynonumberplain}
\theorembodyfont{\upshape}
\theoremsymbol{\rule{1.2ex}{1.2ex}}
\theoremseparator{:}
\newtheorem{myproof}{Proof}

\theoremstyle{nonumberplain}
\theoremheaderfont{\normalfont\bfseries}
\theoremseparator{:}
\theoremsymbol{}
\theorembodyfont{\slshape}

\theoremsymbol{\textbf{***}}

\begin{document}


\title{A Floer homology approach to travelling waves\\ in reaction-diffusion equations on cylinders
\footnotetext{This research was supported by NWO VICI grant 639.033.109 and  NWO TOP grant 613.001.351.}
}

\author{Bente Bakker\thanks{Corresponding author, \href{mailto:bente.h.bakker@gmail.com}{\texttt{bente.h.bakker@gmail.com}}} \and
Jan Bouwe van den Berg \and
Rob Vandervorst \\[1.4 ex]
{\normalsize \textit{Department of Mathematics, VU University Amsterdam} } \\ {\normalsize \textit{De Boelelaan 1081a, 1081 HV Amsterdam, The Netherlands}} }
\date{\normalsize \today}
\maketitle

\begin{abstract}
\normalsize
Travelling waves form a prominent feature in the dynamics of scalar reaction-diffusion equations on unbounded cylinders.
The travelling waves can be identified with the bounded solutions of the elliptic PDE
\begin{equation}
  \label{eq:TWE_abstract}
  \begin{cases}
  \partial_t^2 u - c \partial_t u + \Delta u + f(x,u) = 0 \qquad & t \in \mathbf{R},\; x \in \Omega, \\
  B(u) = 0 & t \in \R,\; x \in \partial\Omega,
  \end{cases}
\end{equation}
where $c \neq 0$ is the wavespeed, $\Omega \subset \mathbf{R}^d$ is a bounded domain, $\Delta$ is the Laplacian on $\Omega$,
and $B$ denotes Dirichlet, Neumann, or periodic boundary data.
We develop a new homological invariant for the dynamics of the bounded solutions of \eqref{eq:TWE_abstract}.
Restrictions on the nonlinearity $f$ are kept to a minimum, for instance, any nonlinearity exhibiting polynomial growth in $u$ can be considered.
In particular, the set of bounded solutions of the travelling wave PDE may not be uniformly bounded.
Despite this, the homology is invariant under lower order (but not necessarily small) perturbations of the nonlinearity $f$,
thus making the homology amenable for computation.
Using the new invariant we derive lower bounds on the number of bounded solutions of \eqref{eq:TWE_abstract},
thus obtaining existence and multiplicity results for travelling wave solutions of reaction-diffusion equations on unbounded cylinders.
\end{abstract}

\vfill

\paragraph{\footnotesize 2010 Mathematics Subject Classification.}
{\footnotesize
Primary: 35J20, 57R58; \hspace{1 pt} Secondary: 35K57, 37B30, 37B35.
}


\vspace{1 ex}

\thispagestyle{empty}
\newpage

\setcounter{tocdepth}{1}
\tableofcontents

\binoppenalty=\maxdimen
\relpenalty=\maxdimen
\sloppy

\section{Introduction}

A prominent feature of reaction-diffusion equations is the formation of spatial and temporal patterns.
The formation of spatial patterns is often observed to be in the form of a travelling wave
invading one state (e.g.\  a homogeneous distribution) and leaving behind another (more complicated) state (e.g.\  a spatial pattern).
In this paper we develop a topological invariant based on a Floer homology construction.
This gives insight in the structure of the solutions of the reaction-diffusion equations. Furthermore,
it demonstrates emphatically that 
Floer homology has applications to a broad class of evolutionary PDEs,
far beyond the realm of symplectic topology where it is traditionally employed.

Historically Floer homology is defined for the Hamilton action functional in order to develop a Morse type theory
for contractible \mbox{period-1} orbits.
In particular, this approach has led to the resolution of the Arnol'd conjecture in many settings, see~\cite{MR0690288, floer1989symplectic, MR1688434, MR1642105} and the references therein.
In the classical context Floer homology gives an algebraic invariant which is related to a weighted count of
critical points of the Hamilton action and is isomorphic to the singular homology of the symplectic manifold.
Floer homology has been developed for numerous nonclassical settings, including~\cite{angenent1999superquadratic, Fabert2015floer, hohloch2009hypercontact, Isobe2017morse}.
The basic idea in the construction is that solutions of a differential equation can be organized using  gradient dynamical systems.

The main message of this paper is that such an approach works for the much larger class of \emph{gradient-like} dynamical systems, including strongly indefinite ones, and
may be regarded as an extension of the Conley index for elliptic partial differential equations.
%
Indeed, in the theory of pattern formation the differential equations in question often display  canonical
gradient-like behaviour.

In this introduction we will start off with an overview of the main results and explain the advantages of the Floer homology approach,
followed by a summary of the Floer construction.
We conclude the introduction with an example of a classical travelling wave problem using the Conley index and point out the analogues
with the Floer homology approach.


\subsection{Main results}

We consider a scalar reaction-diffusion equation on an unbounded cylindrical 
domain $\R \times \Omega$%
\begin{equation*}\label{eq:reaction_diffusion}
  \partial_s \phi = \Delta_{\bar x} \phi + f(\bar x,\phi), \qquad \text{for} \quad s \in \R,\; \bar x \in \R\times\Omega,
  \tag{\textrm{RDE}}
\end{equation*}
together with Dirichlet, Neumann, or periodic boundary conditions for $\bar x \in \R \times \bdy \Omega$.
Here $u$ is a scalar function, and $\Omega \subset \R^d$ is a bounded domain with smooth boundary.
The operator $\Delta_{\bar x}$ denotes the Laplacian on $\R \times \Omega$, that is,
\[
\Delta_{\bar x} = \partial_{x_0}^2 + \partial_{x_1}^2 + \cdots + \partial_{x_d}^2, 
\qquad \text{for} \quad \bar x = (x_0,x) = (x_0,x_1,\dots,x_d) \in \R \times \Omega.
\]
We will also be using the Laplacian on $\Omega$, which we will denote by $\Delta$, that is,
\[
\Delta = \partial_{x_1}^2 + \cdots + \partial_{x_d}^2, \qquad \text{for} \quad x = (x_1,\dots,x_d) \in \Omega.
\]

Suppose the nonlinearity $f$ is homogeneous in the (unbounded) $x_0$ variable.
Then a natural class of solutions (often observed experimentally) 
of \eqref{eq:reaction_diffusion} to consider are of the form $\phi(s,\bar x) = u(x_0 + c s , x_1, \dots, x_d)$,
for some $c \neq 0$ (without loss of generality we will assume $c > 0$).
Then $u(t,x)$ (where $t \in \R$, $x \in \R^d$) satisfies the elliptic PDE
\begin{equation}
  \label{eq:TWE_unperturbed}
    \partial_t^2 u - c\partial_t u + \Delta u + f(x,u) = 0, \qquad \text{for} \quad t \in \R,\; x \in \Omega,
\end{equation}
together with Dirichlet, Neumann, or periodic boundary conditions at $x \in \bdy \Omega$.

If $u$ is a solution of \eqref{eq:TWE_unperturbed}, then $\phi(s,\bar x) = u(x_0 + c s , x_1 , \dots , x_d)$
is called a \emph{travelling wave} (but we will also refer to $u$ as such) 
if it converges (locally uniformly in $x$) as $s \to \pm\infty$ to stationary solutions of \eqref{eq:reaction_diffusion}.
To make this more precise we first need to define $\alpha$- and $\omega$-limit sets.
Let $\alpha(u)$ denote the set of all accumulation points in the $C^1_{\text{loc}}$ topology of the shifts $u(\cdot + \tau , \cdot)$ as $\tau \to - \infty$.
Similarly, let $\omega(u)$ denote all those accumulation points of shifts $u(\cdot + \tau , \cdot)$ as $\tau \to \infty$.
Then $u$ is a travelling wave when $\alpha(u) \cap \omega(u) = \emptyset$,
and $\alpha(u)$ and $\omega(u)$ consist solely of stationary solutions of \eqref{eq:TWE_unperturbed},
i.e.\ each $z \in \alpha(u) \cup \omega(u)$ satisfies $\Delta z + f(x,z) = 0$ and the same boundary conditions as were chosen in \eqref{eq:TWE_unperturbed}.
Later on we will see that, under suitable conditions on the nonlinearity $f$,
any bounded solution of \eqref{eq:TWE_unperturbed} is either a stationary solution or a travelling wave.

In section \ref{sec:perturbations} we will formulate precise conditions on the nonlinearity $f$ for which our theory works.
Special instances of such nonlinearities are of the form
\begin{equation}
  \label{eq:f_odd}
  f_{\text{odd},\pm}(x,u) = \pm \alpha(x) |u|^{p-1} u + h(x,u)
\end{equation}
or
\begin{equation}
  \label{eq:f_even}
  f_{\text{even},\pm}(x,u) = \pm \alpha(x) |u|^p + h(x,u).
\end{equation}
Here $\alpha \in C_b^4(\Omega)$ is such that $\inf_{x\in\Omega} \alpha(x) > 0$, 
and the lower order term $h \in C^4(\Omega \times \R)$ is such that
\[
\limsup_{|u|\to\infty} \sup_{x\in\Omega} \frac{ |h(x,u)| }{ |u|^p } = 0.
\]
It should be stressed that, although the names $f_{\text{odd},\pm}$ and $f_{\text{even},\pm}$ are suggestive, we do \emph{not} assume any symmetry of the lower order term $h$.
For the power $p$ we restrict attention to the ``subcritical range'': $1 < p < \infty$ if $\dim \Omega = 1$,
and $1 < p \leq 3$ if $\dim \Omega = 2$.
Extension of the theory which also deals with higher dimensional domains and bigger $p$ are subjects for future research,
see also remark \ref{remark:Lp_extension} in section \ref{sec:perturbations}.

Stationary solutions of \eqref{eq:TWE_unperturbed}, i.e.\  solutions which are independent of $t$,
solve a semilinear elliptic problem on a bounded domain $\Omega$.
We say that a stationary solution $z$ is hyperbolic if the only solution of
\[
\Delta v + f_u(x,z) v = 0, \qquad \text{for} \quad  x \in \Omega,
\]
together with the same boundary condition considered in \eqref{eq:TWE_unperturbed} is $v \equiv 0$.

In this paper we develop an algebraic/topological invariant which takes into account
solutions of \eqref{eq:TWE_unperturbed} which are stationary, as well as certain solutions which connect stationary solutions (i.e. certain travelling waves).
In terms of applications, the main result from this paper is the following theorem.
\begin{mythmintro}[theorem \ref{thm:forced_existence_of_waves} from section \ref{sec:application}]
\label{thm:intro_existence_thm}
Consider any wave speed $c \neq 0$, and let $k \geq 1$.
Then the following holds:
\begin{itemize}

\item
If $f = f_{\text{odd},-}$ and \eqref{eq:TWE_unperturbed} has at least $2k$ distinct hyperbolic stationary solutions,
then \eqref{eq:reaction_diffusion} has at least $k$ distinct travelling wave solutions of wave speed $c$.
More precisely, to each given hyperbolic stationary solution $z$ (but with the possible exception of at most one of them),
there corresponds at least one travelling wave $u$
such that $\alpha(u) = \{z\}$ or $\omega(u) = \{z\}$ (but it is possible that $\omega(u)$ resp.\ $\alpha(u)$ consist of non-hyperbolic stationary solutions).

\item
If either $f = f_{\text{odd},+}$, or $f = f_{\text{even},-}$, or $f = f_{\text{even},+}$,
and \eqref{eq:TWE_unperturbed} has at least $2k-1$ distinct hyperbolic stationary solutions,
then \eqref{eq:reaction_diffusion} has at least $k$ distinct travelling wave solutions of wave speed $c$.
More precisely, to each given hyperbolic stationary solution $z$,
there corresponds at least one travelling wave $u$
such that $\alpha(u) = \{z\}$ or $\omega(u) = \{z\}$ (but it is possible that $\omega(u)$ resp.\ $\alpha(u)$ consist of non-hyperbolic stationary solutions).

\end{itemize}
Furthermore, in each of these cases there exists at least one more stationary solution (which might be non-hyperbolic).
\end{mythmintro}
Here we consider two travelling wave solutions $\phi_1$, $\phi_2$ of \eqref{eq:reaction_diffusion} to be distinct from each other
if $\phi_1(s_1,\cdot) \neq \phi_2(s_2,\cdot)$ for all $s_1, s_2 \in \R$,
i.e.\  they are not simply time translates of one another.
We note that the problem of establishing hyperbolic stationary solutions is  amenable to computer-assisted proof techniques, making the assumptions in theorem~\ref{thm:intro_existence_thm} verifiable in practice (see e.g.~\cite{ArioliKoch,BGLV}). 

When $f = f_{\text{odd},-}$, equation \eqref{eq:TWE_unperturbed} is dissipative, meaning that the set of all bounded solutions is compact.
In that case, similar results have previously been obtained using different methods.
See for example \cite{gardner1986existence, mielke1994essential, fiedler1998large, gkeba1999conley}.
However, the methods used there break down when the set of bounded solutions is not compact.

We point out that for suitable conditions on $f$ (e.g.\  demanding that $f = f_{\text{odd},+}$ has the symmetry $f(x,-u) = - f(x,u)$)
equation \eqref{eq:TWE_unperturbed} has infinitely many stationary solutions.
This is well known, see for example \cite{bahri1992solutions, yang1998nodal, ramos2009bahri, yu2014solutions, hajlaoui2015morse},
although this can also be deduced from theorem \ref{thm:intro_existence_thm} directly.
Now theorem \ref{thm:intro_existence_thm} implies that if all the stationary solutions are hyperbolic 
(a condition which can be ensured by adding a small perturbation), 
for such anti-symmetric nonlinearities $f = f_{\text{odd},+}$ there exist infinitely many travelling wave solutions of \eqref{eq:reaction_diffusion} 
with any given nonzero wave speed.

When $\Omega$ is zero-dimensional, a classical approach to proving the existence of connecting orbits in \eqref{eq:TWE_unperturbed}
is by using Conley index theory, see also subsection \ref{subsec:Conley_example}.
However, since \eqref{eq:TWE_unperturbed} is a strongly indefinite problem, 
arguments based on Conley index theory cannot be applied directly.
Indeed, any index pair for a stationary solution is homotopy equivalent to a pointed infinite dimensional sphere,
hence the Conley index of any rest point is trivial.
In \cite{fiedler1998large} this problem was circumvented by, roughly speaking,
assigning an index to isolated invariant sets via the limit of Conley indices of finite dimensional approximations of \eqref{eq:TWE_unperturbed}.
In order for this limit to make sense, one needs global compactness results on the
set of all bounded solutions of \eqref{eq:TWE_unperturbed},
i.e.\  this method is only applicable for dissipative nonlinearities $f$.

Previous work by various authors has shown that Floer homology is capable of dealing with a larger class of problems
than the analogous Conley index approach:
compare e.g.\  \cite{conley1983birkhoff, conley1984morse} with \cite{floer1989symplectic, salamon1999lectures},
and \cite{izydorek2002conley} with \cite{angenent1999superquadratic}.
Inspired by this, we construct a Floer-type homology theory for \eqref{eq:TWE_unperturbed}.
This construction only requires a local compactness result on the space of bounded solutions of \eqref{eq:TWE_unperturbed},
hence our results also apply to nondissipative nonlinearities (e.g.\  $f = f_{\text{odd},+}$).

\subsection{Comparison to classical Floer theory}

We want to point out here that this work is not a straightforward application of the standard Floer theory for Hamiltonian systems.
Equation \eqref{eq:TWE_unperturbed} takes over the role of the (perturbed) Cauchy-Riemann equation in the standard Floer theory.
Although both \eqref{eq:TWE_unperturbed} and the Cauchy-Riemann equation are elliptic, there are some important differences.

Equation \eqref{eq:TWE_unperturbed} is not a (formal) gradient flow (but it is gradient-like), 
and it is not immediately be clear whether the construction of the homology
still works for gradient-like equations.
For example, the index theory for this problem becomes more involved.
The obtained index can be related to the classical Morse index of a related parabolic equation, thus making the index amendable for computations.
The existence of a Lyapunov function for \eqref{eq:TWE_unperturbed} guarantees that, just as in the classical situation, 
the moduli spaces of connecting orbits can be compactified by adding broken orbits.

In the case of standard Floer theory, transversality can be obtained by perturbing the Hamiltonian and the almost complex structure.
Thus the perturbed equation is again a PDE.
However, to us it seems that a natural way to achieve generic transversality is by adding a nonlocal term to equation \eqref{eq:TWE_unperturbed}.
Of course, the downside of this approach is that the perturbed equation is then no longer a PDE.

To illustrate why this is a sensible choice, it is best to rewrite the equation \eqref{eq:TWE_unperturbed} as a dynamical system
\begin{equation}
  \label{eq:dynsys_formulation}
  \partial_t U = A(U).
\end{equation}
Here $U = (u , \partial_t u)$, and $A(U)$ is a differential operator acting on $U$ plus a nonlinear term $f(x,u)$.
But from a different viewpoint, $A$ is a (densely defined) vector field on a function space $X$ 
(consisting of functions depending on $x \in \Omega$).
To obtain generic transversality, one should allow for perturbations of \eqref{eq:dynsys_formulation}
which are localised in the ``phase space'' $X$.
In general, there seems to be no reason to assume such a perturbation can be chosen to be a differential operator.
The fact that this is possible in classical Floer theory is because solutions of perturbed Cauchy-Riemann equations 
share most properties with holomorphic functions, see \cite{floer1995transversality}.

The choice to perturb the equation \eqref{eq:TWE_unperturbed} out of the class of PDEs does introduce
a number of new technical obstacles.
The biggest hurdle turns out to be the unique continuation theory developed in section \ref{sec:unique_continuation}.
To deal with the nonlocal perturbation we had to develop a new variety of Carleman estimates.

Finally, we note that with classical Floer theory one is interested in the generators of the homology.
The boundary operator, which counts connecting orbits for the gradient-flow, is merely introduced in order to define the homology.
In contrast to this, we are interested in the connecting orbits of a gradient-like equation.
Hence the boundary operator encodes the information we are actually interested in.
The latter is comparable to the connection matrix in Conley index theory.

\subsection{Future work}
We have chosen to present the theory only for $\Omega$ of dimension $1$ and $2$.
This allows us to work solely over Hilbert spaces.
The advantage of this becomes apparent especially in sections \ref{sec:compactness} and \ref{sec:unique_continuation}.
By replacing the various Sobolev spaces by their $L^p$ counterparts, 
one should be able to obtain similar results for higher dimensional $\Omega$.

We restrict to spaces $\Omega$ which are either tori, or smooth domains in $\R^d$.
It appears to be straightforward to generalize the current results to work for more general spaces $\Omega$.
A natural requirement would be that $\Omega$ is an orientable Riemannian manifold, either closed, or with cylindrical ends.
It would then be natural to allow for mixed boundary conditions at the cylindrical ends.

Another aim for future work is to extend the invariant to higher order equations.
One of the main technical hurdles will then be the extension of the 
unique continuation theory from section \ref{sec:unique_continuation}.

In this article we only consider travelling waves in a scalar reaction-diffusion equation.
However, the same construction will work, essentially without modifications, for systems of reaction-diffusion equations,
provided that the reaction term is of the form $f = \nabla F$.

The invariant that we develop only incorporates index $1$ connecting orbits, for any wave speed $c \neq 0$.
However, a similar invariant can be developed to detect index $0$ orbits, for a specific wave speed $c = c_*$.
For this, one should incorporate a slow drift on the wave speed $c$, which connects a slow wave speed system with a fast wave speed system.
Using similar techniques as presented here the existence of nontrivial orbits for such systems can be shown.
By combining this with a priori estimates on the solutions, we are convinced that one can prove the existence of index $0$ orbits for a specific wave speed $c = c_*$,
see~\cite{MischaikMrozek, Smoller}.


In \cite{Fabert2015non} (and, using different methods, in \cite{abbondandolo2015non}) symplectic nonsqueezing results are developed in the setting of Hamiltonian PDEs.
These rigidity properties are derived using the analysis of $J$-holomorphic curves.
The current paper demonstrates that Floer homology can be constructed using travelling waves instead of $J$-holomorphic curves.
This suggests the PDEs dealt with in this paper possess additional rigidity properties which are worth exploring in more detail.

\subsection{A classical example}

\begin{figure}[t]
  \centering
  \def\svgwidth{0.4\textwidth}
  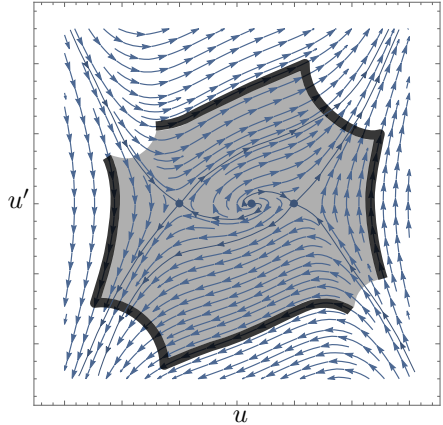
  \caption{Sketch of an index pair $(N,L)$ for the set of all bounded solutions of \eqref{eq:ODE_TWE},
    with $f(u) = (u - a)(1 - u^2)$, where $0 < a < 1$ and $c > 0$.
  The lightly shaded area indicates a possible choice of an isolating neighbourhood $N$,
  while the darkly shaded area indicates a possible choice for an exit set $L$.
  The number of heteroclinic connections depend on choices of $a$ and $c$.
}
  \label{fig:Conley_computation}
\end{figure}

\label{subsec:Conley_example}
To illustrate how a topological invariant can be used to deduce the existence of solutions of \eqref{eq:TWE_unperturbed},
we now briefly recall how this problem can be tackled using standard tools when $\Omega$ is zero-dimensional.
In that setting \eqref{eq:TWE_unperturbed} reduces to the ODE
\begin{equation}
  \label{eq:ODE_TWE}
  u'' - c u' + f(u) = 0.
\end{equation}
A topological approach to existence of solutions of this ODE dates back to work by Conley and Gardner,
see \cite{conley1980application, gardner1984existence}.

We first recall the definition of the Conley index.
Given a flow $(\phi_t)_t$ on a metric space $X$, a pair $(N,L)$ is called an index pair if, roughly speaking,
$L \subset N \subset X$ are compact subsets, such that $N$ and $N \setmin L$ are isolating neighbourhoods of the flow with $\Inv(N) = \Inv(N \setmin L)$,
and all orbits which leave $N$ must do so through $L$ without re-entering $N \setmin L$.
The (homological) Conley index of $(N,L)$ is then defined as the relative (singular) homology of the pair $(N,L)$.
In this example we will use $\Z_2$ coefficients for the homology.
It can be shown that any isolated invariant set $S$ for the flow $(\phi_t)_t$ admits an index pair $(N,L)$, 
i.e., an index pair for which $\Inv(N) = S$.
Moreover, if $(N_1,L_1)$ and $(N_2,L_2)$ are two index pairs for the same isolated invariant set $S$,
then the relative (singular) homologies of those pairs are isomorphic via a natural isomorphism.
Thus one can define the (homological) Conley index $\HC_*(S,\phi)$, up to natural isomorphism, as the relative homology of an index pair $(N,L)$ for $S$.
That is to say, $\HC_*(S,\phi)$ should be interpreted as an equivalence class of relative homologies.
This notion of defining $\HC_*(S,\phi)$ up to natural isomorphisms can be formalized by defining $\HC_*(S,\phi)$ as the inverse limit over all index pairs $(N,L)$ with $\Inv(N) = S$.

A crude way to apply Conley index theory to \eqref{eq:ODE_TWE} makes use of the direct sum property of the index.
This property states that, if an isolated invariant set $S$ for the flow $(\phi_t)_t$ can be written as the disjoint union
of $S_1$, and $S_2$, then
\[
\HC_*(S,\phi) \iso \HC_*(S_1,\phi) \oplus \HC_*(S_2,\phi).
\]
Consider \eqref{eq:ODE_TWE} as a dynamical system on $\R^2$,
and let $S \subset \R^2$ consist of all bounded orbits of this dynamical system.
Now note that for $c = 0$ the dynamical system is Hamiltonian,
and for $c > 0$ the system displays gradient-like behaviour,
with the original Hamiltonian function now strictly decreasing along nonstationary orbits.
This gradient-like behaviour implies that $S$ consists of stationary solutions and heteroclinic orbits.

Using the invariance property of the Conley index it can be seen that $\HC_*(S,\phi)$ is isomorphic to $\HC_*(\{0\},\psi)$,
where $(\psi_t)_t$ is the flow of \eqref{eq:ODE_TWE} with $f(u) = -u^3 - u$.
Since $0$ is a saddle point for $(\psi_t)_t$, it follows that $\HC_*(\{0\},\psi)$ is isomorphic to the reduced singular homology of a 1-sphere.
Hence
\[
\HC_n(S,\phi) \iso \HC_n(\{0\},\psi) \iso \homology_n(\sphere^1,\pt;\Z_2) \iso
\begin{cases}
  \Z_2 & \text{if } n = 1, \\
  0 & \text{otherwise}.
\end{cases}
\]
This is further illustrated in figure \ref{fig:Conley_computation}.

Now suppose the system does not possess any connecting orbits.
The gradient-like behaviour then implies that $S$ consists solely of rest points of the flow,
hence the direct sum property implies
\[
\HC_*(S,\phi) \iso \HC_*(\{(-1,0)\},\phi) \oplus \HC_*(\{(a,0)\},\phi) \oplus \HC_*(\{(1,0)\},\phi).
\]
But local phase plane analysis shows that all of the rest points have nontrivial Conley index.
Hence the rank of the direct sum is at least $3$, while $\HC_*(S,\phi)$ is of rank $1$.
This contradiction allows us to conclude that $S$ contains at least one connecting orbit.

\subsection{Outline of the paper}
As was already indicated, the proof of theorem \ref{thm:intro_existence_thm} is of the same spirit as the example sketched in the previous section.
Note that one of the essential ingredients in this approach is the gradient-like behaviour,
i.e., that the set of bounded solutions consists of stationary solutions and connecting orbits.
The other essential ingredient is the existence of an algebraic object associated to isolating neighbourhoods, such that:
\begin{enumerate}
\item It is amenable to computation, which follows from the invariance of the algebraic object under (not necessarily small) perturbations of the nonlinearity $f$.
In other words, it is a topological invariant.
\item It encodes dynamical information. In particular, it should satisfy a direct sum property.
\end{enumerate}
Since the equation we consider is infinite dimensional the Conley index is not applicable.
In this paper we develop a new topological invariant, which we call the travelling wave homology.

Let us now sketch the construction of the travelling wave homology as well as give an outline of the paper.
We begin by assigning an index $\mu_f(Z)$ to hyperbolic stationary solutions $Z$ of \eqref{eq:TWE_unperturbed}.
This index can be thought of as a generalization of the classical Morse index.
The existence of this (relative) index in our strongly indefinite setting relies on a version of the Fredholm alternative for \eqref{eq:TWE_unperturbed},
for which the hyperbolicity of the stationary solutions is needed,
and which is discussed in section \ref{sec:fredholm}.

The construction of the invariant then relies on a careful analysis of the spaces $\calM(Z_-,Z_+)$ of connecting orbits 
between fixed stationary solutions $Z_-$, $Z_+$.
One important observation is that these spaces are compact modulo ``broken trajectories''.
In section \ref{sec:moduli} we give precise definitions as well as a proof of this property.
Essential ingredients are the local compactness results from section \ref{sec:compactness},
as well as the existence of a Lyapunov function $\calE$.

Using the rapid decay of connecting orbits towards stationary solutions (discussed in section \ref{sec:exp_decay})
the spaces $\calM(Z_-,Z_+)$ can be described as the zero set of a differential operator defined between certain affine Hilbert spaces.
Thus, roughly speaking, if the image of this differential operator intersects the zero section transversely,
the implicit function theorem (making use of the Fredholm theory from section \ref{sec:fredholm})
gives us a manifold structure on $\calM(Z_-,Z_+)$.
In fact, this manifold is finite dimensional, with dimension equal to $\mu(Z_-) - \mu(Z_+)$.

As it turns out, the natural way to ensure transversality holds generically is by perturbing~\eqref{eq:TWE_unperturbed}
using a small \emph{nonlocal} term.
The perturbed equation takes the form
\begin{equation}
  \label{eq:TWE_perturbed}
  \partial_t^2 u(t,x) - c \partial_t u(t,x) + \Delta u(t,x) + f(x,u(t,x)) + g(x,(u(t,\cdot),\partial_t u(t,\cdot))) = 0.
\end{equation}
We stress that the perturbation $g$ depends on $u(t,\cdot)$ and $\partial_t u(t,\cdot)$ as functions on $\Omega$.
A typical example of such a perturbation is of the form displayed in \eqref{eq:g_representative}.
One particular part where this nonlocal term prevents us from applying known results is the unique continuation theory developed in 
section \ref{sec:unique_continuation}.
There we prove that if two solutions $(u,\partial_t u)$ and $(v,\partial_t v)$ of~\eqref{eq:TWE_perturbed} coincide at a certain time $t = t_0$,
they must in fact coincide for all $t \in \R$, i.e.\  $u \equiv v$.
In the absence of the term $g$ this would follow from classical Carleman estimates 
\cite{carleman1939probleme, aronszajn1956unique, aronszajn1962unique, kenig1986carleman, tataru2004unique, colombini2006strong}.
To deal with the nonlocal perturbation $g$, in section \ref{sec:unique_continuation} we have derived a version of the Carleman estimates
where the function is not required to be localised except for the $t$-direction.
This is possible, at the cost of the Carleman estimates no longer holding uniformly over all localised functions (as in the classical case),
but the way in which the estimates depend on the chosen function works well together with localising a solution of~\eqref{eq:TWE_unperturbed}
in the $t$-direction using cutoff functions.
This allows us to prove the desired uniqueness result.
In the end, we are able to show that for generic choices of $g$ the transversality condition is satisfied (see section \ref{sec:genericity}).

Of particular interest is then the case when $\calM(Z_-,Z_+)$ is of dimension $2$ and noncompact.
A careful analysis shows that this space can be compactified by adding, for each noncompact connected component of $\calM(Z_-,Z_+)$,
precisely two broken trajectories.
In the definition of the homology we shall also make use of isolating neighbourhoods
(the precise definition of which will be given in section~\ref{sec:homology}), which will play a similar role as in Conley theory.
Then, if one lets $C_n$ denote the group which is $\Z_2$-generated by index $n$ stationary solutions contained in a fixed isolating neighbourhood $N$,
and define homomorphisms $\partial_n : C_n \to C_{n-1}$ by counting (modulo~$2$) connecting orbits
which are contained in $N$,
it follows that $\partial_n \circ \partial_{n+1} \equiv 0$.
This way we arrive at the following theorem.
\begin{mythmintro}[theorem \ref{thm:homology_definition} from section \ref{sec:homology}]
One has $\partial_{n} \circ \partial_{n+1} = 0$, and consequently,
\[
\HTW_n(N,f,g,c) := \homology_n(C_*,\partial_*) = \frac{ \ker{ \partial_n } }{ \img{ \partial_{n+1} } }
\]
is well-defined.
\end{mythmintro}

The resulting homology $\HTW_*(N,f,g,c)$ is independent (up to natural isomorphisms) of the particular choice of $g$ (for small $g$).
Thus we obtain an invariant $\HTW_*(N,f,c)$ for~\eqref{eq:TWE_unperturbed}.
Furthermore, if $(f_\lambda,c_\lambda)$ is a homotopy between nonlinearities $(f_0,c_0)$ and $(f_1,c_1)$,
and $N$ satisfies an appropriate stability property with respect to this homotopy (see section \ref{sec:homology} for precise details),
then $\HTW_*(N,f_0,c_0) \iso \HTW_*(N,f_1,c_1)$.
In particular, when $N$ is the entire phase space and the nonlinearity $f$ is of the form given in \eqref{eq:f_odd} or \eqref{eq:f_even},
then the homology is invariant under homotopy on the lower order term $h$ and the coefficient $\alpha$,
as long as $\inf_{x\in\Omega} \alpha_\lambda(x) > 0$ uniformly in the homotopy parameter $\lambda$.
This is what allows us to determine explicitely the homology of the global dynamics in all four cases, see theorem~\ref{thm:computation_homology}.

The invariant $\HTW$ satisfies a direct sum property similar to that of the Conley index.
If $N = A \cup B$, where $A$ and $B$ are disjoint isolating neighbourhoods for the dynamics, then
\[
\HTW_*(N,f,c) \iso \HTW_*(A,f,c) \oplus \HTW_*(B,f,c).
\]
This allows us to prove theorem \ref{thm:intro_existence_thm} in a way analogous to the 
simplified example involving Conley index theory in section \ref{subsec:Conley_example}.

\section{The extended problem}
\label{sec:perturbations}

In this section we set up the extended problem into which our original problem can be embedded.
In later sections we will see that for generic choices from this class of extended problems we can set up the desired transversality theory.
It appears that this is not possible without considering the extended problem.

\subsection{Perturbations}
Let either $\Omega \subset \R^d$ be a bounded domain with smooth boundary, or $\Omega = \T^d = \R^d / \Z^d$.
We will restrict ourselves to $d \in \{1,2\}$, see remark \ref{remark:Lp_extension}.
Since we want to consider \eqref{eq:TWE_unperturbed} as a dynamical system we are going to 
rewrite it as a system of equations which involve first order derivatives of $t$ only.
We choose to incorporate the boundary conditions in the phase space.

We denote by $H^k(\Omega)$ the $L^2$ Sobolev space, defined as the closure of $C^\infty(\Omega)$ in the norm
\[
\| u \|_{H^k(\Omega)} = \left( \sum_{|\alpha|\leq k} \| \partial^\alpha u \|_{L^2(\Omega)}^2 \right)^{1/2}.
\]
When $\Omega \subset \R^d$ is a domain, let $B : C^\infty(\cl \Omega) \to C^\infty(\bdy\Omega)$ be given by either 
$B(u) = \rst{u}{\bdy\Omega}$ (Dirichlet), or $B(u) = \rst{\partial_\nu u}{\bdy\Omega}$ (Neumann), 
where $\nu$ denotes the outward pointing unit normal on $\bdy\Omega$.
When $\Omega$ is a torus (corresponding to periodic boundary conditions), we set $B = 0$.
For any $k \in \N_0$ we now define
\[
H^k_B(\Omega) := \clin{ \set{ u \in C^\infty(\cl\Omega) }{ B(u) = 0 } }{H^k(\Omega)}.
\]
Whenever it is convenient, the operator $B$ shall also be identified with its extension to Sobolev spaces $H^k(\Omega)$.

We can now introduce the phase spaces
\[
X^k := H^{k+1}_B(\Omega) \times H^k_B(\Omega), \qquad k \in \N_0.
\]
Together with the norms
\[
\| (u,v) \|_{X^k}^2 := \big( \| u \|_{H^{k+1}(\Omega)}^2 + \| v \|_{H^k(\Omega)}^2 \big)^{1/2}, \qquad (u,v) \in X^k
\]
these become separable Hilbert spaces.

Given a possibly unbounded open subset $J \subset \R$ define
\[
 W^{k,2}_{\text{loc}}(J ; X^0 , \dots, X^k) := W^{k,2}_{\text{loc}}(J , X^0) \cap W^{k-1,2}_{\text{loc}}(J , X^1) \cap \cdots 
\cap W^{1,2}_{\text{loc}}(J,X^{k-1}) \cap L^2_{\text{loc}}(J , X^k).
\]
We endow these spaces with the compact-open topology.
Convergence in this topology, which in fact makes $W^{k,2}_{\text{loc}}(J;X^0,\dots,X^k)$ into a Fr\'echet space,
is characterised as follows:
a sequence $(U_n)_n \subset W^{k,2}_{\text{loc}}(J;X^0,\dots,X^k)$ converges towards $U_\infty$ 
if and only if for any bounded open subset $J' \subset J$ it holds that
\begin{equation}
  \label{eq:weak_convergence_W^m}
  \max_{ 0 \leq \ell \leq k } \; \max_{ 0 \leq j \leq k - \ell } \; \int_{J'} \| \partial_t^j U_n(t) - \partial_t^j U_\infty(t) \|_{X^\ell}^2 \d t \to 0 
\qquad \text{as} \quad n \to \infty.
\end{equation}
The spaces $W^{m,2}(J;X^0,\dots,X^m)$ are defined in a similar fashion, where now $J'$ in \eqref{eq:weak_convergence_W^m}
is replaced by $J$.
The spaces $C^m(J;X^0,\dots,X^m)$ are defined by
\[
C^m(J;X^0,\dots,X^m) := C^m( J,X^0) \cap C^{m-1}( J,X^1) \cap \cdots \cap C^1( J,X^{m-1}) \cap C^0(J,X^m),
\]
where the topology is defined by uniform convergence of functions and their derivatives on compact subsets of $J$.
The spaces of bounded differentiable functions $C^m_b(J;X^0,\dots,X^m)$ are defined in an analogous manner,
where now the convergence is uniform over $J$ itself.
\begin{mylemma}
  One has the continuous embeddings 
  \begin{align*}
    W^{m,2}(J;X^0,\dots,X^m) &\hookrightarrow C_b^{m-1}(J;X^0,\dots,X^{m-1}), \\
    W_{\text{loc}}^{m,2}(J;X^0,\dots,X^m) &\hookrightarrow C^{m-1}(J;X^0,\dots,X^{m-1}).
  \end{align*}
Furthermore, the embeddings
\[
\begin{array}{l l l}
  W^{m,2}(J;X^0,\dots,X^m) &\hookrightarrow W^{m-1,2}(J;X^0,\dots,X^{m-1}), \qquad &\text{with $J$ bounded}, \\
 W_{\text{loc}}^{m,2}(J;X^0,\dots,X^m) &\hookrightarrow W_{\text{loc}}^{m-1,2}(J;X^0,\dots,X^{m-1}), \qquad &\text{for any $J$}
\end{array}
\]
are compact.
\end{mylemma}
\begin{myproof}
We give a sketch here, for more details we refer to \cite{amann1995linear}.
The first two statements are a consequence Morrey's inequality.
This relies on the integral representation $U(t) = U(t_0) + \int_{t_0}^t \partial_s U(s) \d s$, 
which is well-defined since the spaces $X^i$ are separable.
The compact embeddings follow from the Rellich–Kondrachov theorem for vector-valued Sobolev spaces. 
Here one uses that $X^i$ are seperable Banach spaces and the embeddings $X^i \hookrightarrow X^{i-1}$ are compact.
\end{myproof}

Let $A_{f,g,c} : X^1 \to X^0$ be the nonlinear operator defined by
\[
A_{f,g,c}( u , v )(x) :=
\begin{pmatrix}
  -v(x) \\
  \Delta u(x) + f(x,u(x)) - c v(x) + g(x,(u,v))
\end{pmatrix}
\]
for $(u,v) \in X^1$.
We stress here that the term $g(x,(u,v))$ depends on the functions $u$ and $v$, and not on the point $(u(x),v(x))$.
We assume $f : \Omega \times \R \to \R$ is of class $C^m$ with $m \geq 1$, 
and $g : \Omega \times X^0 \to \R$ is $C^m$ (in the Fr\'echet sense), 
and $c > 0$.
At each point in the development of the theory we will point out exactly how big $m$ needs to be,
but we already want to point out that all theorems hold for $m \geq 4$.
The nonlocal term $g$ will typically be a very small term.
Additional restrictions on $f$ and $g$ will be formulated in the next section.
For brevity we shall write $A$ instead of $A_{f,g,c}$ whenever this does not give rise to ambiguity.

We will study the behavior of the dynamical system 
\begin{equation*}
  \label{eq:TWE}
  \tag{\textrm{TWE}}
  \partial_t U + A_{f,g,c}(U) = 0, \qquad U \in W^{1,2}_{\text{loc}}(J ; X^0 , X^1).
\end{equation*}
Note that $U = (u,v)$ is a solution of \eqref{eq:TWE} if and only if $v = \partial_t u$, and
\[
\partial_t^2 u - c \partial_t u + \Delta u + f(x,u) + g(x,(u,\partial_t u)) = 0 \qquad \text{on} \quad J \times \Omega,
\]
and for each $t \in J$ the boundary condition $B(u(t,\cdot)) = 0$ is satisfied.
Unless mentioned otherwise, we assume that $J = \R$.
Note that as a consequence of the nonlocal perturbation $g$ \eqref{eq:TWE} is in general not a PDE.

\subsubsection{Conditions on \texorpdfstring{$f$}{f}, \texorpdfstring{$g$}{g}, and \texorpdfstring{$c$}{c}}
Henceforth we shall assume that $(f,g) \in C^m(\Omega\times\R) \times C^m(\Omega\times X^0)$ and $c > 0$
satisfy the following hypotheses.

\newcounter{f_counter}
\newcounter{f_altcounter}
\newcounter{g_counter}
\begin{enumerate}[label=\textbf{(f\arabic*)}, ref=(f\arabic*), series=flist]
\item\label{H:growth_bound_f}
There exists $C_f \geq 0$ such that $f$ satisfies the growth bounds
\[
\sup_{x\in\Omega} |f(x,u)| \leq C_f ( 1 + |u|^p ),
\]
where $1 \leq p < \infty$ if $\dim \Omega = 1$, and $1 \leq p \leq 3$ if $\dim \Omega = 2$.


\item\label{H:contact_condition_odd}
There exist some $-1 < \theta < 1$ and $C_f' \geq 0$ such that $f$ satisfies
\[
|F(x,u)| \leq C_f' + \frac{\theta}{2} f(x,u) u.
\]
Here $F(x,u) = \int_0^u f(x,s) \d s$.

\setcounter{f_altcounter}{\value{enumi}}
\end{enumerate}
As an alternative to \ref{H:contact_condition_odd} we can also consider (see also remark \ref{remark:practical_restrictions_f})
\begin{enumerate}[start=\value{f_altcounter}, label=\textbf{(f\arabic*$'$)}, ref=(f\arabic*$'$), series=faltlist]

\item \label{H:contact_condition_even}
There exist some $-1 < \theta < 1$ and $C_f' \geq 0$ such that $f$ satisfies
\[
|F(x,u)| \leq C_f' + \frac{\theta}{2} f(x,u) |u|.
\]
\end{enumerate}
When dealing with Neumann or periodic boundary conditions, we need an additional restriction on the nonlinearity,
given by \ref{H:superlinear_growth_f}.
\begin{enumerate}[resume, label=\textbf{(f\arabic*)}, ref=(f\arabic*), series=flist]

\item\label{H:superlinear_growth_f}
When considering Neumann or periodic boundary data, assume $f$ satisfies the superlinear growth condition
\[
\liminf_{|u|\to\infty} \inf_{x\in\Omega} \bigg| \frac{f(x,u)}{u} \bigg| > 0.
\]

\setcounter{f_counter}{\value{enumi}}
\end{enumerate}
Besides these restrictions on $f$ we need to put a few restrictions on $g$:
\begin{enumerate}[label=\textbf{(g\arabic*)}, ref=(g\arabic*), series=glist]

\item\label{H:growth_bound_g}
There exists a constant $C_{0,g}$ such that
\[
\sup_{x\in\Omega,\; U \in X^0} | g(x,U) | \leq C_{0,g},
\]
and for each $k \in \{1,\dots,m\}$ there exists a constant $C_{k,g}$ such that
\[
\sup_{x\in\Omega,\; U \in X^0} \| \dn k g(x,U) \|_{\calL( ( \R\times X^0 )^k ,\R)} \leq C_{k,g}.
\]

\item\label{H:cone_condition}
The perturbation $g$ satisfies the Lipschitz condition
\[
\sup_{x\in\Omega,\; u \in H^1_B(\Omega)} | g(x,u,v) | \leq \frac{ c }{ 2 \sqrt{ \Vol(\Omega) } } \| v \|_{L^2(\Omega)}.
\]
Note the dependence of the Lipschitz constant on the wave speed $c$.

\item\label{H:vanishing_condition}
The perturbation $g$ satisfies
\[
\d g(x,u,0) = 0 \qquad \text{for all} \quad x \in \Omega,\; u \in H_B^1(\Omega).
\]

\setcounter{g_counter}{\value{enumi}}
\end{enumerate}
\newcommand{\Fgroup}{(f1)--(f\arabic{f_counter})\xspace}
\newcommand{\Ggroup}{(g1)--(g\arabic{g_counter})\xspace}
\newcommand{\Hgroup}{(f1)--(f\arabic{f_counter}) and (g1)--(g\arabic{g_counter})\xspace}

\begin{myremark}
\label{remark:Lp_extension}
We want to stress here that both the dimensional restriction $d = \dim \Omega \leq 2$ 
and the growth restriction \ref{H:growth_bound_f} are merely technical.
The dimensional restriction ensures that the Sobolev embedding $H^2(J \times \Omega) \hookrightarrow C^0(J\times\Omega)$
holds, where $J \subset \R$ is a domain.
This fact is used in order to obtain the compactness results in section \ref{sec:compactness}.
The choice of $p$ in hypothesis \ref{H:growth_bound_f} ensures that the Sobolev embedding 
$H^1(J\times\Omega) \hookrightarrow L^{2p}(J\times\Omega)$ holds.
Consequently, the growth bound on $f$ ensures that the map $A_{f,g,c} : X^1 \to X^0$ is indeed well-defined,
bounded, and continuous (see \cite{appell1990nonlinear}).
Both these conditions can be relaxed by replacing all the spaces $H^k = W^{k,2}$ by $W^{k,q}$ spaces, 
for appropriately chosen $q$.
Since we mainly want to convey the idea that Floer theory can be applied to travelling wave problems in \eqref{eq:reaction_diffusion},
we have chosen to stick with the Hilbert space theory in order to reduce the technicality of the estimates,
which tend to complicate the spirit of the arguments.
\end{myremark}

\begin{myremark}
\label{remark:practical_restrictions_f}
  Conditions \Ggroup could seem restrictive,
but recall that the nonlinear term $g$ is only introduced to put the equation \eqref{eq:TWE}
into ``general position'' (i.e.\  to achieve transversality).
In the application of the theory we are eventually interested in the case where $g=0$,
hence conditions \Ggroup are not particularly restricting.

On the other hand, conditions \Fgroup clearly put restrictions
on the types of nonlinearities to which the theory is applicable.
Examples (but not exhausting all possibilities) of such nonlinearities are
$f = f_{\text{odd},\pm}$ and $f = f_{\text{even},\pm}$ which where introduced in equations \eqref{eq:f_odd} and \eqref{eq:f_even}.
Then clearly hypotheses \ref{H:growth_bound_f} and \ref{H:superlinear_growth_f} are satisfied.
The appropriate choice between hypotheses \ref{H:contact_condition_odd} and \ref{H:contact_condition_even}
and the corresponding value of $-1 < \theta < 1$ can be summarised as follows:
\begin{center}
  \begin{tabular}{ l || l | l | l }
    & & $\sigma = -$ & $\sigma = +$ \\
    \hline\hline
    $f_{\text{odd},\sigma}$ & \ref{H:contact_condition_odd}  & $\theta < - 2 / (p+1)$ & $\theta > 2/(p+1)$ \\
    $f_{\text{even},\sigma}$ & \ref{H:contact_condition_even} & $\theta < - 2/(p+1)$ & $\theta > 2/(p+1)$
  \end{tabular}
\end{center}
As was already indicated in the introduction, the difference between the various choices of $f$ is also reflected in the possible dynamics,
a fact which we will return to once we discuss applications of the theory in section \ref{sec:application}.
\end{myremark}

\subsubsection{Conditions on the nonautonomous equation} \label{sec:nonautonomous_conditions}
In order to develop continuation of the Floer homology groups, we need to allow $(f,g,c)$ to depend explicitly on $t$; 
i.e.\  consider a nonautonomous version of \eqref{eq:TWE}.
We shall assume that $t$-dependent 
$(f,g,c) \in C^m(\R\times\Omega\times\R) \times C^m(\R\times \Omega \times X^0) \times C^m(\R,(0,\infty))$ 
satisfy the following hypotheses.

\newcounter{C_counter}
\begin{enumerate}[label=\textbf{(n\arabic*)}, ref=(n\arabic*), series=Clist]

\item\label{C:pointwise_H}
For each $t \in \R$, the triple $(f(t,\cdot,\cdot), g(t,\cdot,\cdot), c(t))$ satisfy hypotheses \Hgroup,
with the constants $C_f$, $C_f'$, $\theta$, and $C_{k,g}$ uniform in $t \in \R$.

\item\label{C:asymptotic_constant}
There exists an $\ell > 0$ and $t$-independent triples $(f_-,g_-,c_-)$, $(f_+,g_+,c_+)$, such that 
\[
\begin{cases}
  \big( f(t,\cdot,\cdot) , g(t,\cdot,\cdot) , c(t) \big) = \big( f_- , g_- , c_- \big) & \text{for } t \leq -\ell, \\
  \big( f(t,\cdot,\cdot) , g(t,\cdot,\cdot) , c(t) \big) = \big( f_+ , g_+ , c_+ \big) & \text{for } t \geq +\ell.
\end{cases}
\]

\end{enumerate}
The final hypothesis makes use of a sufficiently small constant $\Theta$.
How small this $\Theta$ should really be depends on $\theta$, $\inf_{t\in\R} c(t)$, and $\ell$,
and will be dictated by theorems \ref{thm:solution_compactness_Dirichlet_nonautonomous}
and \ref{thm:solution_compactness_Neumann_nonautonomous}.
\begin{enumerate}[resume, label=\textbf{(n\arabic*)}, ref=(n\arabic*), series=Clist]

\item\label{C:growth_bound_f}
There exist $\Theta \geq 0$ sufficiently small, and $C_f'' \geq 0$, such that
\[
  |\partial_t F(t,x,u)| \leq C_f'' + \Theta | F(t,x,u) |.
\]


\setcounter{C_counter}{\value{enumi}}
\end{enumerate}

\newcommand{\Cgroup}{(n1)--(n\arabic{C_counter})\xspace}

\subsection{Auxiliary definitions}

\subsubsection{The energy/Lyapunov functional}
\label{subsec:Lyapunov}
Recall from hypothesis \ref{H:contact_condition_odd} that $F : \Omega \times \R^d \to \R$ is chosen such that $F(\cdot,0) = 0$, $\nabla_u F(x,u) = f(x,u)$.
By hypothesis \ref{H:growth_bound_f} and the continuous embedding $H^1(\Omega) \hookrightarrow L^{p+1}(\Omega)$,
it follows that the induced Nemytskii operator
\[
F : H^1(\Omega) \to L^1(\Omega), \qquad u(x) \mapsto F(x,u(x))
\]
is bounded and $C^1$, see e.g.\  \cite{appell1990nonlinear}.
We can therefore define energy functional ${\calE_f \in C^1(X^0, \R)}$ by
\begin{equation*}
    \calE_f(u,v) = \int_\Omega - \frac 1 2 |v(x)|^2 + \frac 1 2 |\nabla u(x)|^2 - F(x,u(x)) \d x.
\end{equation*}
Here $|\cdot|$ denotes the Euclidean norm.
When the choice of $f$ is clear from the context, we shall sometimes abbreviate $\calE_f$ by $\calE$.

Note that, in light of the continuous embedding $W^{1,2}_{\text{loc}}(\R;X^0,X^1) \hookrightarrow C^0(\R,X^0)$, the map
\[
  \calE_f : W^{1,2}_{\text{loc}}(\R ; X^0 , X^1) \to C^0(\R), \qquad U \mapsto \calE_f(U(\cdot))
\]
is $C^1$.
Let $\calA$ consist of all $U = (u,\partial_t u)$ that solve \eqref{eq:TWE}.
Elliptic regularity theory combined with hypothesis \ref{H:growth_bound_f} implies that in particular $\calA \subset C^m(\R;X^0,X^1)$.
Therefore $\calE_f$ restricts to a continuous map
\[
  \calE_f : \calA \to C^m(\R), \qquad U \mapsto \calE_f(U(\cdot)).
\]
\begin{myremark}
  Details on the regularity of $\calA$ can be found in the proof of theorem \ref{thm:solution_compactness_Dirichlet}.
Although step 1 of the proof depends on the regularity of the map $\calE_f$, steps 2 until 4 do not rely on such properties of $\calE_f$.
The argument given in steps 2 until 4 of the proof can be modified to show that $\calA \subset C^m(\R;X^0,\dots,X^m)$
whenever $f$ is of class $C^m$.
\end{myremark}

The derivative of $\calE_f$ along $U \in \calA$ has the form
\begin{align*}
  \ddx{\calE_f(U(t))}{t} &= \int_\Omega \bigg( - \partial_t^2 u(t,x) - \Delta u(t,x) - f(x,u(t,x)) \bigg) \partial_t u(t,x) \d x \\
  &= - c \| \partial_t u(t,\cdot) \|_{L^2(\Omega)}^2 + \langle g(\cdot,U(t)) , \partial_t u(t,\cdot) \rangle_{L^2(\Omega)}.
\end{align*}
Hypothesis \ref{H:cone_condition} then implies that
\begin{equation}
  \label{eq:Lyapunov_function_decay}
  \ddx{\calE_f(U(t))}{t} \leq - \frac c 2 \| \partial_t u(t,\cdot) \|_{L^2(\Omega)}^2,
\end{equation}
thus $\calE_f$ is a Lyapunov function for \eqref{eq:TWE}.

\begin{myremark}
\label{remark:strict_Lyapunov}
  In fact, from our regularity theory (section \ref{sec:compactness})
and the unique continuation theorem (theorem \ref{thm:unique_continuation}) it will follow that $\calE_f$ is a strict Lyapunov function.
That is, inequality \eqref{eq:Lyapunov_function_decay} is strict unless $\partial_t u(t) = 0$ for all $t$.
\end{myremark}

\subsubsection{Stationary solutions and hyperbolicity}
Denote by $\calS(f) \subset X^1$ the collection of stationary solutions (also referred to as rest points) of \eqref{eq:TWE}, i.e.\  $\calS = A_f^{-1}(0)$.
Given $-\infty \leq a \leq b \leq \infty$, we define $\calS_a^b(f) := \calS(f) \cap \calE_f^{-1}([a,b])$.
Whenever the choice of $f$ is clear from the context, it will be suppressed in the notation.

Given $Z = (z,0) \in \calS(f)$, we will see in section \ref{sec:compactness} (more specifically, see theorem \ref{thm:operator_regularity})
that the Nemytskii operator
\[
f : H_B^1(\Omega) \to L^2(\Omega), \qquad u(x) \mapsto f(x,u(x))
\]
is $m$ times continuously differentiable near $z$ 
whenever $f : \Omega \times \R \to \R$ is of class $C^m$.
In particular, the operator $A : X^1 \to X^0$ is differentiable near $Z$.
By hypothesis \ref{H:vanishing_condition} the linearised operator looks like
\begin{align*}
  &\d A(Z) : X^1 \to X^0, \\
  &\d A(Z) =
  \begin{pmatrix}
    0 & -1 \\
    \Delta + f_u(x,z) & -c
  \end{pmatrix}.
\end{align*}
To do spectral theory we shall consider the linear extension of this operator to the complexified Banach spaces $X^k_\C := X^k \times i X^k$.
We will say that $Z$ is hyperbolic if the linearised operator $\d A(Z)$, 
considered as an unbounded operator on $X^0_\C$ with domain $\calD(\d A(Z)) = X^1_\C$,
has its spectrum disjoint from the imaginary axis, i.e.\  $\sigma(\d A(Z)) \cap i\R = \emptyset$.
Denote by $\calS_{\text{hyp}}(f)$ the collection of all hyperbolic rest points.
A nonlinearity $f$ for which all rest points are hyperbolic shall be called \emph{regular}.

Note that $\d A(Z)$ is a compact perturbation of the operator $(u,v) \mapsto ( -v , \Delta u )$, hence it is Fredholm of index $0$.
Hence, if $Z$ is hyperbolic, the inverse function theorem can be applied, thus ensuring that hyperbolic rest points are isolated in $X^1$.
Later on, in section \ref{sec:genericity}, we will see that hyperbolicity can always be acquired by a slight perturbation of the nonlinearity $f$.

\subsubsection{Connecting orbits and transversality}
\label{subsec:transversality_definition}
A solution $U$ of \eqref{eq:TWE} is called a connecting orbit if there exist $Z_- \in \calS(f_-)$, $Z_+ \in \calS(f_+)$
such that $\| U(t) - Z_\pm \|_{X^0} \to 0$ as $t \to \pm\infty$.
We will later on see that any bounded solution of \eqref{eq:TWE} is in fact either an equilibrium or a connecting orbit.
Also note that, in light of the existence of the Lyapunov function $\calE_f$, connecting orbits in the autonomous equation are heteroclinic orbits, i.e.\  $Z_- \neq Z_+$.
Thus \eqref{eq:TWE} is a gradient-like system.

We need to introduce another technical condition.
A connecting orbit $U$ is said to be transversal provided that the linearised operator
$\partial_t + \d A(U(t))$ (which according to theorem \ref{thm:operator_regularity} is well-defined)
is surjective when considered as an operator from $W^{1,2}(\R;X^0,X^1)$ to $L^2(\R;X^0)$.

In section \ref{sec:fredholm} it is shown that there is a natural way to assign an index to 
connecting orbits between hyperbolic rest points.
Equation \eqref{eq:TWE} is said to satisfy the transversality condition up to order $k$ 
if all connecting orbits of index at most $k$ are transversal.
In section \ref{sec:genericity} we will see that, whenever $f$ is of class $C^m$, 
transversality up to order $m-1$ can always be obtained by choosing generic nonlocal perturbations $g$.
Such a $g$ shall then be called \emph{regular}.

\section{Regularity and compactness}\label{sec:compactness}

In this section we will see that the collection of solutions of \eqref{eq:TWE} is locally compact.
The relatively compact neighborhoods are determined by sub-superlevel sets of the energy functional.
These results, in a way, replace the Palais-Smale condition which appears in classical Morse theory,
and will form one of the cornerstones in defining the Floer boundary operator.

Throughout this section we let $J = (j_-,j_+)$, where $-\infty \leq j_- < j_+ \leq +\infty$.
Given numbers $a, b \in \R$, define
\[
  \calA_a^b(J,f,g,c) := \set{ U \in W^{1,2}_{\text{loc}}(J;X^0,X^1) }{ 
    \begin{matrix}
      \partial_t U + A_{f,g,c}(U) = 0 \text{ on } J, \\
      a \leq \liminf_{t\nearrow j_+} \calE_{f(t,\cdot,\cdot)}(U(t)), \\
      \limsup_{t\searrow j_-} \calE_{f(t,\cdot,\cdot)}(U(t)) \leq b
    \end{matrix}
    },
\]
where $\calE_f$ is as defined in section \ref{subsec:Lyapunov}.
Note that when $f$, $g$, and $c$ are independent of $t$ and hypothesis \ref{H:cone_condition} is satisfied, 
then the set $\calA_a^b(J,f,g,c)$ consists of all solutions to \eqref{eq:TWE} whose energy remains between $a$ and $b$.
Whenever the choices of $f$, $g$, and $c$ are clear we shall suppress them from the notation.
We will also write $\calA_a^b$ instead of $\calA_a^b(\R,f,g,c)$.

\subsection{Compactness of \texorpdfstring{$\calA_a^b(J,f,g,c)$}{A(J,f,g,c)} with Dirichlet boundary data}

\subsubsection{The autonomous case}
We have the following compactness result.
\begin{mythm}\label{thm:solution_compactness_Dirichlet}
  Consider Dirichlet boundary data.
Let $(f,g)$ be of class $C^m$, with $m \geq 1$.
Suppose hypotheses \ref{H:growth_bound_f}, either \ref{H:contact_condition_odd} or \ref{H:contact_condition_even}, 
 \ref{H:growth_bound_g}, and \ref{H:cone_condition} are satisfied.
  Then for any $J = (j_-,j_+) \subset \R$ and $J' \subset \cl{J'} \subset J$,
  the set $\rst{ \calA_a^b(J,f,g,c) }{ J' }$ is bounded in $C_b^m(J';X^0,\dots,X^m)$ and compact in $W^{m,2}_{\text{loc}}(J';X^0,\dots,X^m)$.
\end{mythm}
Here $\rst{ \calA_a^b(J,f,g,c) }{ J' } = \set{ \rst{U}{J'} }{ U \in \calA_a^b(J,f,g,c) }$.
Note that we cannot obtain compactness of $\calA_a^b(J,f,g,c)$ itself, since solutions may blow up as $t$ approaches $j_-$ or $j_+$.
However, if $J = \R$ we do retrieve compactness of the full space $\calA_a^b(\R,f,g,c)$. 

\begin{myproof}[of theorem \ref{thm:solution_compactness_Dirichlet}]
The proof is split into four steps.
In the first step we will use hypotheses \ref{H:growth_bound_g}, \ref{H:cone_condition} and 
either \ref{H:contact_condition_odd} or \ref{H:contact_condition_even}
to obtain a local $H^1$ bound on the solutions.
In the second step we combine these bounds together with hypothesis \ref{H:growth_bound_f} 
and a regularity argument to obtain a local $L^\infty$ bound.
This argument does not immediately extend to higher degrees of regularity, since the Nymetskii operator induced by $f$
potentially does not possess the required regularity.
To circumvent this problem, in the third step a new nonlinearity $\widetilde f$ is constructed which possesses the required regularity,
in such a way that solutions of the original equation \eqref{eq:TWE} are also solutions of the equation with this new nonlinearity.
In the fourth and final step, the desired compactness result are derived from the preceding steps.

\paragraph{Step 1.}
We will first construct a convenient family of test functions.
Let 
\[
\delta := \frac 1 2 \min\{ \inf J' - j_- , j_+ - \sup J' \}.
\]
Then let $\phi_0 : \R \to \R$ be a $C^2$ function such that $\supp(\phi_0) \subset [-\delta,\delta]$,
and $\phi_0(t) \geq 0$ for all $t \in \R$,
and $\phi_0(t) \geq C_{1,\phi} > 0$ for $t \in [-\delta/2,\delta/2]$,
and $|\phi'(t)| \leq C_{2,\phi} \phi(t)^{1/2}$ for all $t \in \R$, for some $C_{2,\phi} > 0$.
For any $\tau \in J'$ we then define $\phi_\tau(t) := \phi_0(t-\tau)$.
Note that the definition of $\delta$ ensures that $\supp(\phi_\tau) \subset J$.

Fix any $U = ( u , \partial_t u ) \in \calA_a^b(J,f,g,c)$.
We shall henceforth identify $u$ with the $\R$-valued function on $J \times \Omega$ given by $u(t,x) = u(t)(x)$.
Pick any $\tau \in J'$, and for the moment abbreviate $\phi_\tau$ by $\phi$.
Letting $Q = J \times \Omega$, we now list some estimates.
%
%
\newcounter{cpt_est_altcounter}
\begin{enumerate}[label=\textbf{(\alph*)}, ref=(\alph*), series=cpt_est_list]

\item
Observe that since $U = (u,\partial_t u)$ is a solution to \eqref{eq:TWE} and hypothesis \ref{H:cone_condition} is satisfied,
estimate \eqref{eq:Lyapunov_function_decay} holds, hence
\begin{equation}
  \label{eq:est_t_derivative}
  \int_Q | \partial_t u |^2 \d x \d t \leq - \frac 2 c \int_J \ddx{ \calE(U(t)) }{ t } \d t \leq \frac 2 c ( b - a ).
\end{equation}

\item 
Note that
\begin{equation}
  \label{eq:31}
  \int_Q \phi | \nabla u |^2 \d x \d t = 2 \int_J \phi(t) \calE(U(t)) \d t +  2 \int_Q \phi |\partial_t u|^2 \d x \d t + 2 \int_Q \phi F(x,u) \d x \d t.
\end{equation}
By hypothesis \ref{H:cone_condition} and since $J$ is connected, 
$t \mapsto \calE(U(t))$ is a monotone function, so that in particular $\calE(U(t)) \leq b$ for all $t \in J$.
Therefore the first term in \eqref{eq:31} is bounded from above.
By estimate \eqref{eq:est_t_derivative} the second term in \eqref{eq:31} is also bounded from above.
Hence
\[
\int_Q \phi | \nabla u |^2 \d x \d t \leq C + 2 \int_Q \phi | F(x,u) | \d x \d t,
\]
where $C \geq 0$ is independent of $U \in \calA_a^b(J,f,g,c)$.
We will now continue estimating the last term.

\item 
If $f$ satisfies hypothesis \ref{H:contact_condition_odd},
using the fact that $U = (u,\partial_t u)$ solves \eqref{eq:TWE} and partial integration we obtain
\begin{equation}
  \label{eq:intermediate_F_bound}
  \begin{split}
    2 & \int_Q \phi | F(x,u) | \d x \d t \leq 2 C_f' + \theta \int_Q \phi f(x,u) u \d x \d t \\
    &= 2 C_f' - \theta \int_Q \phi \big( \partial_t^2 u + \Delta u - c \partial_t u + g(x,U(t)) \big) u \d x \d t \\
    &= 2 C_f' + \theta \int_Q  \phi' u \partial_t u + \phi \big( |\partial_t u|^2 +  |\nabla u|^2 + c u \partial_t u - g(x,U(t)) u \big) \d x \d t \\
    &\leq 2 C_f + |\theta| \int_Q \phi |\partial_t u|^2 \d x \d t + |\theta| \int_Q \phi |\nabla u|^2 \d x \d t \\
    &\quad + |\theta| \int_Q ( |\phi'| + c \phi ) |u| |\partial_t u| \d x \d t + |\theta| \int_Q \phi | g(x,U(t)) | |u| \d x \d t.
  \end{split}
\end{equation}

We will now bound the last two terms in \eqref{eq:intermediate_F_bound}.
In light of Cauchy's inequality, for any $\mu > 0$ there exists a $0 < C_\mu < \infty$
\[
  \int_Q \phi |g(x,U(t))| |u| \d x \d t \leq \mu \int_Q \phi |u|^2 \d x \d t + C_\mu \int_Q \phi |g(x,U(t))|^2 \d x \d t.
\]
Using hypothesis \ref{H:growth_bound_g}, the last term can be estimated from above by some constant $\widetilde C_\mu$.
Now recall that $|\phi'(t)| \leq C_{2,\phi} \phi(t)^{1/2}$.
Hence, by again using Cauchy's inequality, for any $\nu > 0$ there exists a $0 < C_\nu < \infty$ such that
  \begin{align*}
    \int_Q ( |\phi'| + c \phi ) |u| |\partial_t u| \d x \d t 
    &\leq \int_Q \big( \phi^{1/2} |u| \big) \big( ( C_{2,\phi} + c \phi^{1/2} ) | \partial_t u | \big) \d x \d t \\
    &\leq \nu \int_Q \phi |u|^2 \d x \d t + C_\nu \int_Q ( C_{2,\phi} + c \phi^{1/2} )^2 |\partial_t u|^2 \d x \d t \\
    &\leq \nu \int_Q \phi |u|^2 \d x \d t + \widetilde C_\nu \int_Q |\partial_t u|^2 \d x \d t,
  \end{align*}
where $\widetilde C_\nu = C_\nu \| C_{2,\phi} + c \phi^{1/2} \|_{L^\infty(\R)}^2$.

Combining these estimates with \eqref{eq:intermediate_F_bound}, we obtain
\begin{equation}
  \label{eq:F_bound}
  \begin{split}
      2 \int_Q \phi |F(x,u)| \d x \d t &\leq 2 C_f + \widetilde C_\mu + |\theta| ( 1 + \widetilde C_\nu ) \int_Q |\partial_t u|^2 \d x \d t \\
&\quad + |\theta| \int_Q \phi |\nabla u|^2 + (\mu+\nu) \phi |u|^2 \d x \d t.
  \end{split}
\end{equation}

Using estimate \eqref{eq:est_t_derivative} the first integral is bounded from above by a constant which is independent of $U \in \calA_a^b(J,f,g,c)$.
Hence we have found that there exists a constant $C_{\theta,\mu,\nu}$ independent of $U \in \calA_a^b(J,f,g,c)$ such that
\begin{equation}
  \label{eq:12}
      2 \int_Q \phi |F(x,u)| \d x \d t \leq  C_{\theta,\mu,\nu} + |\theta| \int_Q \phi |\nabla u|^2 + (\mu+\nu) \phi |u|^2 \d x \d t.
\end{equation}

\setcounter{cpt_est_altcounter}{\value{enumi}}
\end{enumerate}
\begin{enumerate}[start=\value{cpt_est_altcounter}, label=\textbf{(\alph*$'$)}, ref=(\alph*$'$), series=cpt_est_list]

\item\label{est:primitive_bound_alt}
If on the other hand hypothesis \ref{H:contact_condition_even} holds, we still get the same estimate as above.
Care needs to be taken to avoid the non-differentiability of $u\mapsto|u|$ around $u = 0$, which prevents us from applying integration by parts.
As a workaround, we first select a function $\eta \in C^1(\R)$ such that 
$\eta(u) = |u|$ for $|u| \geq 1$, and $|\eta(u)| \leq 1$ for $|u| \leq 1$, and $|\eta'(u)| \leq 1$ for all $u$.
For example, one can consider
\[
\eta(u) =
\begin{cases}
  \frac 1 2 u^2 + \frac 1 2 & \text{if } |u| \leq 1, \\
  |u| & \text{if } |u| \geq 1.
\end{cases}
\]
Then, after updating the constant $C_f'$, we have
\[
|F(x,u)| \leq C_f' + \frac \theta 2 f(x,u) \eta(u).
\]
Then we estimate
\begin{align*}
  2 & \int_Q \phi |F(x,u)| \d x \d t \leq 2 C_f' + \theta \int_Q \phi f(x,u) \eta(u) \d x \d t \\
  &\leq 2 C_f' + |\theta| \int_Q \phi |\eta'(u)| |\partial_t u|^2 + \phi |\eta'(u)| |\nabla u|^2 \d x \d t \\
  &\quad + |\theta| \int_Q ( |\phi'| + c \phi ) |\eta(u)| |\partial_t u| \d x \d t + |\theta| \int_Q \phi | g(x,U(t)) | |\eta(u)| \d x \d t \\
  &\leq 2 C_f' + |\theta| \int_Q \phi |\partial_t u|^2 + \phi |\nabla u|^2 \d x \d t \\
  &\quad + |\theta| \int_Q ( |\phi'| + c \phi ) |u| |\partial_t u| \d x \d t + |\theta| \int_Q \phi | g(x,U(t)) | |u| \d x \d t \\
  &\quad + |\theta| \int_{Q_{|u|\leq 1}} ( |\phi'| + c \phi ) |\eta(u)| |\partial_t u| \d x \d t + |\theta| \int_{Q_{|u|\leq 1}} \phi | g(x,U(t)) | |\eta(u)| \d x \d t,
\end{align*}
where $Q_{|u|\leq 1} = \set{ (t,x) \in Q }{ |u(t,x)| \leq 1 }$.
Now by using that $|\eta(u)|\leq 1$ for $|u| \leq 1$, and the observation that the $L^2$-norm of $\partial_t u$ is bounded by a constant independent
of $U \in \calA_a^b(J,f,g,c)$, we see that the last two integrals in the above estimate are bounded by a constant independent of $U \in \calA_a^b(J,f,g,c)$.
Therefore, proceeding as before we again arrive at estimate \eqref{eq:F_bound}, and consequently \eqref{eq:12}.

\end{enumerate}
Combining all these estimates, we obtain the following gradient bound.
For any $U = ( u , \partial_t u ) \in \calA_a^b(J,f,g,c)$ we have
\begin{equation}
  \label{eq:gradient_bound}
    \int_Q \phi |\partial_t u |^2 + \phi |\nabla u|^2 \d x \d t \leq C_{\theta,\mu,\nu} + |\theta| \int_Q \phi |\nabla u|^2 + (\mu + \nu) \phi |u|^2 \d x \d t,
\end{equation}
where $C_{\theta,\mu,\nu} \geq$ is independent of $U \in \calA_a^b$.

Now we use the Dirichlet boundary data to apply Poincar\'e's inequality (with constant $C_P(\Omega)$), which yields
\[
\int_Q \phi |\partial_t u|^2 + \phi |\nabla u|^2 \d x \d t \leq C_{\theta,\mu,\nu} 
+ |\theta| \big( 1 + (\mu+\nu) C_P(\Omega) \big) \int_Q \phi |\nabla u|^2 \d x \d t.
\]
By choosing $\mu, \nu > 0$ sufficiently small and using the fact that $0 \leq |\theta| < 1$, the integral on the right-hand side can be absorbed
into the left-hand side.
Finally, we use that $\phi(t) = \phi_\tau(t) \geq C_{1,\phi} > 0$ for $t \in [\tau - \delta/2 , \tau + \delta/2]$ to obtain
\begin{equation}
  \label{eq:H1_u_bound}
  \int_{Q_\tau} |\partial_t u|^2 + |\nabla u|^2 \d x \d t \leq C \qquad \text{for all} \quad U = (u,\partial_t u) \in \calA_a^b(J,f,g,c).
\end{equation}
Here $Q_\tau = [\tau-\delta/2,\tau+\delta/2] \times \Omega$.
Note that the only way that this constant depends on $\tau$ is via $\| \phi_\tau \|_{W^{1,\infty}(\R)}$ and $\Vol(\supp(\phi_\tau))$, which are in fact independent of $\tau$.

\paragraph{Step 2.}
For any $U = (u,\partial_t u) \in \calA_a^b(J,f,g,c)$, note that
\[
\partial_t^2( \phi_\tau u ) + \Delta( \phi_\tau u ) = \phi_\tau'' u + 2 \phi_\tau' \partial_t u + c \phi_\tau \partial_t u - \phi_\tau f(x,u) - \phi_\tau g(x,U).
\]
By G{\aa}rding's inequality (see e.g.\  \cite{evans1998partial}) applied to the Laplacian $\Delta_{t,x} = \partial_t^2 + \Delta$,
and writing $H^k = H^k(Q)$ as a shorthand, we have
\begin{equation}
  \label{eq:8}
  \begin{split}
  \| \phi_\tau u \|_{H^{k+2}} &\leq C \big(  \| \phi_\tau'' u + 2 \phi_\tau' \partial_t u + c \phi_\tau \partial_t u 
  - \phi_\tau f(x,u) - \phi_\tau g(x,U) \|_{H^k} + \| \phi_\tau u \|_{H^k} \big) \\
  &\leq C \big( \| \phi_\tau'' u \|_{H^k} + \| \phi_\tau u \|_{H^k} + 2 \| \phi_\tau' \partial_t u \|_{H^k} + c  \| \phi_\tau \partial_t u \|_{H^k} \\
  &\quad + \|\phi_\tau f(x,u) \|_{H^k} + \|\phi_\tau g(x,U) \|_{H^k} \big).
  \end{split}
\end{equation}
We note here that by shift invariance of $\Delta_{t,x}$ the constant $C$ can in be chosen independent of $\tau \in J'$.

First we consider the case $k = 0$.
Using \eqref{eq:H1_u_bound} and hypothesis \ref{H:growth_bound_g} we obtain an upper bound for the first four terms and the last term.
Hypothesis \ref{H:growth_bound_f} implies that the Nemytskii operator $f : L^{2p}(Q_\tau) \to L^2(Q_\tau)$ is bounded and continuous,
hence by the Sobolev embedding $H^1(Q_\tau) \hookrightarrow L^{2p}(Q_\tau)$ 
the map $f : H^1(Q_\tau) \to L^2(Q_\tau)$ is bounded and continuous, see \cite{appell1990nonlinear}.
Again using \eqref{eq:H1_u_bound}, we see that also the fifth term in \eqref{eq:8} is bounded above by some constant.
Note that this upper bound is independent of $\tau \in J'$ and $U \in \calA_a^b(J,f,g,c)$.
Since $\phi_\tau(t) \geq C_{1,\phi} > 0$ for $t \in [\tau - \delta/2 , \tau + \delta/2]$, 
we find that there exists some constant $M$ such that
\[
\| u \|_{H^2( Q_\tau )} \leq M \qquad \text{for all} \quad \tau \in J',\; (u,\partial_t u) \in \calA_a^b(J,f,g,c).
\]
Using a Sobolev embedding it then follows that the set
\[
D := \set{ \rst{ u }{ J' \times \Omega } }{ (u,\partial_t u) \in \calA_a^b(J,f,g,c) }
\]
is bounded in $C_b^0(J' \times \Omega)$.

\paragraph{Step 3.}
Estimate \eqref{eq:8} cannot be directly employed for $k \geq 1$, since it is not clear whether the Nemytskii operator 
$f : H^{k+1}(Q_\tau) \to H^k(Q_\tau)$ is bounded (indeed, \ref{H:growth_bound_f} only ensures that the Nemytskii operator is bounded and $C^0$ for $k = 0$).
To circumvent this problem we will consider modified nonlinearities $\widetilde f$.
Set
\[
\rho := \sup_{ u \in D } \| u \|_{L^\infty(J' \times \Omega)},
\]
which in light of step 2 is a finite number.
The definition of $\rho$ implies that if $\widetilde f$ is another nonlinearity which coincides with $f$ on $\Omega \times [-\rho,\rho]$,
then clearly
\[
\partial_t U(t) + A_{\widetilde f}(U(t)) = 0 \quad \text{and} \quad \calE_{\widetilde f}(U(t)) = \calE_f(U(t)), 
\qquad \text{for} \quad U \in \calA_a^b(J,f,g,c),\; t \in J'
\]
hence $\rst{ \calA_a^b(J,f,g,c) }{ J' } \subset \calA_a^b(J',\widetilde f,g,c)$.

Let $\eta \in C^\infty(\R)$ be such that $\eta(u) = 1$ for $|u| \leq \rho$ and $\eta(u) = 0$ for $|u| \geq 2 \rho$.
Now set
\[
\widetilde f(x,u) := \eta(u) f(x,u) + (1 - \eta(u)) u^3.
\]
Then $\widetilde f$ satisfies hypotheses \Hgroup.
Furthermore, $\widetilde f$ induces a bounded $C^{m-k}$ Nemytskii operator from $H^{k+1}(Q_\tau)$ into $H^k(Q_\tau)$, see \cite{appell1990nonlinear}.
After choosing a further subinterval $J'' \subset \cl{ J'' } \subset \inter{ J' }$, the argument from steps 1 and 2 can be repeated
to obtain estimate \eqref{eq:8}, now with $f$ replaced by $\widetilde f$, and this time for $\tau \in J''$.
Inductively we can then obtain a bound for $\|\phi_\tau f(x,u) \|_{H^k}$ with $k \in \{0,\dots,m\}$.
Consequently, there exists some constant $M$ such that
\begin{equation}
  \label{eq:18}
  \| u \|_{H^{m+2}( Q_\tau )} \leq M \qquad 
  \text{for all} \quad \tau \in J'',\; (u,\partial_t u) \in \calA_a^b(J',\widetilde f,g,c),
\end{equation}
in particular this estimate holds for $(u,\partial_t u) \in \calA_a^b(J,f,g,c)$.

\paragraph{Step 4.}
Henceforth without loss of generality replace $J''$ by $J'$ in \eqref{eq:18}.
Note then that \eqref{eq:18} combined with a Sobolev embedding implies that 
$\rst{\calA_a^b(J,f,g,c)}{J'}$ is bounded in the topology of $C_b^m(J';X^0,\dots,X^m)$.
Since the embedding $H^{m+2}( Q_\tau ) \hookrightarrow H^{m+1}( Q_\tau )$ is compact
it follows that $\rst{\calA_a^b(J,f,g,c)}{J'}$ is relatively compact in $W^{m,2}_{\text{loc}}(J';X^0,\dots,X^m)$.
Moreover, hypotheses \ref{H:growth_bound_f} and \ref{H:growth_bound_g} imply that the nonlinear operator
\[
\partial_t + A(\cdot) : W^{m,2}_{\text{loc}}(J';X^0,\dots,X^m) \to L^2_{\text{loc}}(J';X^0)
\]
is continuous.
Hence the limit point $U$ of a sequence $(U_n)_n$ in $\rst{ \calA_a^b(J,f,g,c) }{ J' }$ is a solution of \eqref{eq:TWE} on $J'$.
To see that such a limit point has an extension to a solution of \eqref{eq:TWE} on $J$,
apply steps 1 through 3 with $J'$ replaced by $J''$, where $\cl{ J' } \subset J'' \subset \cl{J''} \subset J$.
We then find that $U_n$ converges over a subsequence to $U'$ in $W^{m,2}_{\text{loc}}(J'';X^0,\dots,X^m)$, and $U'$ solves \eqref{eq:TWE} on $J''$.
By uniqueness of the limits one has $\rst{ U' }{ J'' } = U$.
Since this holds for any such $J''$, we find that $U \in \rst{ \calA_{-\infty}^{+\infty}(J,f,g,c) }{ J' }$.
By continuity of the energy functional $\calE : W^{m,2}_{\text{loc}}(J ; X^0,X^1) \to C^0(J)$ we find that in fact $U \in \rst{ \calA_a^b(J,f,g,c) }{ J' }$.
Hence $\rst{ \calA_a^b(J,f,g,c) }{ J' }$ is compact in $W^{m,2}_{\text{loc}}(J';X^0,X^1)$.
This proves the theorem.
\end{myproof}

\begin{myremark}
\label{remark:Palais-Smale_type_compactness}
Note that the estimates in the proof of theorem \ref{thm:solution_compactness_Dirichlet}
do not depend explicitly on $(f,g,c)$, but only on the constants appearing in hypotheses \Hgroup.
Hence, if $( (f_n,g_n,c_n) )_n$ is a sequence of triplets satisfying hypotheses \Hgroup,
with contants uniform in $n$,
then step 1 through 3 of the proof shows that the sets $\rst{ \calA_a^b(J,f_n,g_n,c_n) }{ J' }$ are bounded in $C_b^m(J';X^0,\dots,X^m)$,
uniformly in $n$.
Suppose $(f,g,c)$ is another triples satisfying hypotheses \Hgroup,
and $(f_n,g_n,c_n) \to (f,g,c)$ as $n \to \infty$, where the convergence takes place in 
$C_{\text{loc}}^m(\Omega\times\R) \times C_b^m(\Omega \times X^0) \times (0,\infty)$.
For each $n$ select a solution $U_n \in \calA_a^b(J,f_n,g_n,c_n)$.
Then a slight adaption of step 4 of the proof shows that
there exists a subsequence $( U_{n_k} )_k$ of $(U_n)_n$ and a solution $U \in \calA_a^b(J,f,g,c)$
such that $U_{n_k} \to U$ as $k \to \infty$, with convergence in $W^{m,2}_{\text{loc}}(J';X^0,\dots,X^m)$.
The same result applies when considering different boundary conditions and/or the nonautonomous equations.
\end{myremark}

\subsubsection{The nonautonomous case}
For $t$-dependent $(f,g,c)$ we have the following compactness result.
\begin{mythm}\label{thm:solution_compactness_Dirichlet_nonautonomous}
Consider Dirichlet boundary data.
Given $\theta \in (-1,1)$, $\gamma > 0$, and $\ell > 0$, there exists $\Theta = \Theta(\theta,\gamma,\ell) > 0$ for which the following is true.
Let $(f,g,c)$ be of class $C^m$ ($m \geq 1$) for which hypotheses \Cgroup are satisfied with the chosen constants $\theta$, $\ell$, $\Theta$,
and $\inf_{t\in\R} c(t) \geq \gamma$.
Fix $J = (j_-,j_+) \subset \R$, and $J' \subset \cl{J'} \subset J$ with $[-\ell,\ell] \subset J'$.
Then the set $\rst{ \calA_a^b(J,f,g,c) }{ J' }$ is bounded in $C_b^m(J';X^0,\dots,X^m)$ and compact in $W^{m,2}_{\text{loc}}(J';X^0,\dots,X^m)$.
\end{mythm}
\begin{myproof}
The argument from the autonomous case cannot be directly applied to the nonautonomous case,
because along a solution $U$ of the nonautonomous equation the energy $t \mapsto \calE_{f(t,\cdot,\cdot)}(U(t))$ may increase.
This is the reason why we introduce hypothesis \ref{C:growth_bound_f}.
Conceptually, this condition allows us to extract an a priori bound for the amount the energy can increase along a solution $U$,
provided that the energy is asymptotically bounded as $t \to j_\pm$.

We will now explain in detail how the proof of theorem \ref{thm:solution_compactness_Dirichlet} can be adapted for the nonautonomous case.
First note that, since the problem is autonomous outside $(-\ell,\ell)$,
 there is no loss of generality in assuming that $J' = (-\ell-\epsilon,\ell+\epsilon)$, where $\epsilon > 0$ is small enough so that $\cl{ J' } \subset J$.
Indeed, suppose the conclusion of the theorem holds for this choice of $J'$, 
so that in particular $\rst{\calA_a^b(J,f,g,c)}{J'}$ is bounded in $C_b^m(J';X^0,\dots,X^m)$.
Hypothesis \ref{C:pointwise_H} then ensures that the map $t \mapsto \calE_{f(t,\cdot,\cdot)}(U(t))$ is bounded for $t \in J'$,
with a bound which is uniform in $U \in \rst{\calA_a^b(J,f,g,c)}{J'}$.
In particular, there exists an $M \geq 0$ such that
\[
| \calE_{f_+}(U(\ell)) - \calE_{f_-}(U(-\ell)) | \leq M \qquad \text{for all} \quad U \in \calA_a^b(J,f,g,c).
\]
Now consider $J' \subset \cl{J'} \subset J$ chosen arbitrarily.
Set $J_- := (j_-,-\ell)$ and $J_-' := J' \cap (j_-,-\ell-\epsilon)$; here we need the $\epsilon > 0$ in order to ensure that $\cl{J_-'} \subset J_-$.
Then
\[
\rst{ \calA_a^b(J,f,g,c) }{ J_-' } \subset \rst{ \calA_{a-M}^b(J_-,f_-,g_-,c_-) }{ J_-' }.
\]
Similarly, with $J_+ := (\ell,j_+)$ and $J_+' := J' \cap (\ell+\epsilon,j_+)$ it holds that
\[
\rst{ \calA_a^b(J,f,g,c) }{ J_+' } \subset \rst{ \calA_a^{b+M}(J_+,f_+,g_+,c_+) }{ J_+' }.
\]
For $J_\pm'$ the conclusion of the theorem thus follows from the analogous result for the autonomous case.
The result for general $J'$ then follows by combining the results for the autonomous and nonautonomous parts.

Henceforth assume $J' = (-\ell-\epsilon,\ell+\epsilon)$.
Compared to the proof of theorem \ref{thm:solution_compactness_Dirichlet} we will use a slightly modified test function $\phi \in C^2(\R)$,
namely, assume $0 \leq \phi(t) \leq 1$ for all $t \in \R$, 
and $\phi(t) = 1$ for $t \in J'$, 
and $\supp(\phi) \subset J$ is compact, 
and $|\phi'(t)| \leq C_{2,\phi} \phi(t)^{1/2}$ for all $t \in \R$, for some $C_{2,\phi} > 0$.
Then, using hypothesis \ref{C:growth_bound_f}, it follows that
\begin{equation}
  \label{eq:modified_growth_bound_f}
  |\partial_t F(t,x,u)| \leq C_f'' \1_{\supp(\phi)}(t) + \Theta \phi(t) |F(t,x,u)|.
\end{equation}
Here we used that $\phi(t) = 1$ for $t \in (-\ell,\ell)$ and that the left-hand side vanishes for $|t| \geq \ell$, in light of hypothesis \ref{C:asymptotic_constant}.

We will now point out how estimates (a)--(c) from step 1 of the proof of theorem \ref{thm:solution_compactness_Dirichlet} can be modified
to the nonautonomous case.
Throughout these estimates, let $U = (u,\partial_t u) \in \calA_a^b(J,f,g,c)$.
%
%
\newcounter{cpt_est_altcounter_nonauto}
\begin{enumerate}[label=\textbf{(\alph*)}, ref=(\alph*), series=cpt_est_list_nonauto]

\item
Note that by hypotheses \ref{C:pointwise_H} and \ref{C:growth_bound_f} the map $t \mapsto \calE_{f(t,\cdot,\cdot)}(U(t))$
is $C^1$, and
\begin{equation}
  \label{eq:33}
  \begin{split}
      \ddx{\calE_{f(t,\cdot,\cdot)}(U(t))}{t}
&= - \langle c(t) \partial_t u(t,\cdot) , \partial_t u(t,\cdot) \rangle_{L^2(\Omega)} + \rst{ \partial_s \calE_{f(s,\cdot,\cdot)}(U(t)) }{ s=t } \\
&\quad + \langle g(t,\cdot,U(t)) , \partial_t u(t,\cdot) \rangle_{L^2(\Omega)} \\
&\leq - \frac 1 2 \langle c(t) \partial_t u(t,\cdot) , \partial_t u(t,\cdot) \rangle_{L^2(\Omega)} + \rst{ \partial_s \calE_{f(s,\cdot,\cdot)}(U(t)) }{ s=t } \\
&\leq - \frac \gamma 2 \| \partial_t u(t,\cdot) \|_{L^2(\Omega)}^2 + \int_\Omega | \rst{ \partial_s F(s,x,u(t,x)) }{ s = t } | \d x.
  \end{split}
\end{equation}
Here the penultimate inequality holds since for each $t$ the pair $( g(t,\cdot,\cdot) , c(t) )$ satisfies hypothesis \ref{H:cone_condition}.
Therefore, using \eqref{eq:modified_growth_bound_f},
\begin{equation}
  \label{eq:nonautonomous_dt_bound}
  \begin{split}
      \int_Q |\partial_t u|^2 \d x \d t &\leq - \frac{2}{\gamma} \int_J \ddx{ \calE_{f(t,\cdot,\cdot)}(U(t))}{t} \d t 
+ \frac{2}{\gamma} \int_Q | \rst{ \partial_s F(s,x,u(t,x)) }{ s = t } | \d x \d t \\
  &\leq \frac{2}{\gamma} \big( b - a + C_f'' \Vol(\supp(\phi) \times \Omega) \big) + \Theta \frac{2}{\gamma} \int_Q \phi | F(t,x,u) | \d x \d t.
  \end{split}
\end{equation}

\item
Estimate \eqref{eq:33} implies that
\[
\ddx{\calE_{f(t,\cdot,\cdot)}(U(t))}{t} \leq  \int_\Omega | \rst{ \partial_s F(s,x,u(t,x)) }{ s = t } | \d x,
\]
so that
\begin{align*}
  \calE_{f(t,\cdot,\cdot)}(U(t)) &\leq b + \int_{j_-}^t \int_\Omega  | \rst{ \partial_s F(s,x,u(t,x)) }{ s = t } | \d x \d t  \\
&\leq b + \int_Q  | \rst{ \partial_s F(s,x,u(t,x)) }{ s = t } | \d x \d t \\
&\leq b + C_f'' \Vol( \supp(\phi) \times \Omega ) + \Theta \int_Q \phi | F(t,x,u) | \d x \d t,
\end{align*}
where we again used \eqref{eq:modified_growth_bound_f}.
Hence
\[
\int_J \phi \calE_{f(t,\cdot,\cdot)}(U(t)) \d t \leq C + \Theta \|\phi\|_{L^1(\R)} \int_Q \phi | F(t,x,u) | \d x \d t
\]
for some constant $C \geq 0$ independent of $U \in \calA_a^b(J,f,g,c)$.
Plugging this into \eqref{eq:31} and combining with \eqref{eq:nonautonomous_dt_bound}, after increasing the constant $C$ we find
\[
\int_Q \phi |\nabla u|^2 \d x \d t \leq C + 2 C_{1,\Theta} \int_Q \phi | F(t,x,u) | \d x \d t,
\]
where
\[
C_{1,\Theta} := \Theta \bigg( \frac 2 \gamma \|\phi\|_{L^\infty(\R)} + \|\phi\|_{L^1(\R)} \bigg) + 1.
\]

\item 
If for each $t$ the function $f(t,\cdot,\cdot)$ satisfies hypothesis \ref{H:contact_condition_odd},
the same computation as in the autonomous case leads to estimate \eqref{eq:F_bound}.
Combining this with estimate \eqref{eq:nonautonomous_dt_bound} results in
\begin{align*}
  2 \int_Q \phi | F(t,x,u) | \d x \d t &\leq C_{\theta,\mu,\nu} + \Theta \frac{2 |\theta| (1+\widetilde C_\nu)}{\gamma} \int_Q \phi |F(t,x,u)| \d x \d t \\
  &\quad + |\theta| \int_Q \phi |\nabla u|^2 + (\mu+\nu) \phi |u|^2 \d x \d t,
\end{align*}
for some $C_{\theta,\mu,\nu}$ independent of $U \in \calA_a^b(J,f,g,c)$.
Consequently,
\[
2 \int_Q \phi | F(t,x,u) | \d x \d t \leq C_{2,\Theta} C_{\theta,\mu,\nu} + |\theta| C_{2,\Theta} \int_Q \phi |\nabla u|^2 + (\mu+\nu) \phi |u|^2 \d x \d t,
\]
where
\[
C_{2,\Theta} = \bigg( 1 - \Theta \frac{|\theta| (1+\widetilde C_\nu)}{\gamma} \bigg)^{-1}.
\]
\end{enumerate}
The modifications needed in (c$'$) are similar to those made for (c).

Combining these estimates, we obtain
\begin{equation*}
  \label{eq:nonautonomous_gradient_bound}
  \begin{split}
  \int_Q \phi |\partial_t u |^2 + \phi | \nabla u |^2 \d x \d t &\leq C + 2 C_{3,\Theta} \int_Q \phi | F(t,x,u) | \d x \d t \\
  &\leq C_{\Theta,\theta,\mu,\nu} + |\theta| C_{2,\Theta} C_{3,\Theta} \int_Q \phi |\nabla u|^2 + (\mu+\nu) \phi |u|^2 \d x \d t,
  \end{split}
\end{equation*}
where $C_{3,\Theta} = \Theta \gamma^{-1} \| \phi \|_{L^\infty(\R)} + C_{1,\Theta}$, and $C, C_{\Theta,\theta,\mu,\nu} \geq 0$ are some constants independent of $U \in \calA_a^b(J,f,g,c)$.
Now observe that for $\Theta > 0$ sufficiently small 
(depending upon $\theta$, $\gamma$, and $\ell$ only)
we have  $|\theta| C_{2,\Theta} C_{3,\Theta} < 1$.
Therefore, if hypothesis \ref{C:growth_bound_f} is satisfied with a small enough $\Theta$,
one can proceed as in the autonomous case to arrive at estimate \eqref{eq:H1_u_bound}.
Steps 2, 3, and 4 in the proof of \ref{thm:solution_compactness_Dirichlet} are also valid for the nonautonomous case.
Hence the theorem is proven.
\end{myproof}

\subsection{Compactness of \texorpdfstring{$\calA_a^b(J,f,g,c)$}{A(J,f,g,c)} with Neumann or periodic boundary data}

\subsubsection{The autonomous case}

When dealing with Neumann or periodic boundary conditions, the question of compactness becomes more delicate.
We can no longer use Poincar\'e's inequality in order to bound the $L^2$-norm of $u$ in terms of the $L^2$-norm of $\nabla u$.
To obtain such a bound, we will need the superlinear growth condition \ref{H:superlinear_growth_f} on $f$.
\begin{mythm}\label{thm:solution_compactness_Neumann}
Consider Neumann or periodic boundary data.
Let $(f,g)$ be of class $C^m$, with $m \geq 1$.
Assume hypotheses \Fgroup, \ref{H:growth_bound_g}, and \ref{H:cone_condition} are satisfied.
  Then for any $J = (j_-,j_+) \subset \R$ and $J' \subset \cl{J'} \subset J$,
  the set $\rst{ \calA_a^b(J,f,g,c) }{ J' }$ is bounded in $C_b^m(J';X^0,\dots,X^m)$ and compact in $W^{m,2}_{\text{loc}}(J';X^0,\dots,X^m)$.
\end{mythm}
\begin{myproof}
First note that we can split hypothesis \ref{H:superlinear_growth_f} into four cases:
\begin{equation*}
  \begin{array}{l l}
  \displaystyle \liminf_{|u|\to\infty} \inf_{x\in\Omega} \frac{ f(x,u) }{|u|} > 0,  \quad & \displaystyle \liminf_{|u|\to\infty} \inf_{x\in\Omega} \frac{ f(x,u) }{u} > 0, \\
  \displaystyle \limsup_{|u|\to\infty} \sup_{x\in\Omega} \frac{ f(x,u) }{|u|} < 0,  \quad & \displaystyle \limsup_{|u|\to\infty} \sup_{x\in\Omega} \frac{ f(x,u) }{u} < 0.
  \end{array}
\end{equation*}
We will assume that the first case holds; the proof for the other three cases goes in a similar fashion.
There exist $\epsilon > 0$ and $K \geq 0$ such that
\[
\epsilon |u|^2 \leq f(x,u) |u| \qquad \text{for all} \quad x \in \Omega,\; |u| \geq K.
\]
Let $\eta \in C^1(\R)$ be as defined in estimate \ref{est:primitive_bound_alt} in the proof of theorem \ref{thm:solution_compactness_Dirichlet}.
Then there exists a constant $M$ such that
\[
\epsilon |u|^2 \leq M + f(x,u) \eta(u) \qquad \text{for all} \quad x \in \Omega,\; |u| \geq K.
\]

Fix any $U = ( u , \partial_t u ) \in \calA_a^b(J,f,g,c)$.
For a.e. $t$ we have
\begin{align*}
  \epsilon \| u(t,\cdot) \|_{L^2(\Omega)}^2 &= \epsilon \int_{\set{ x \in \Omega }{ |u(t,x)| < K }} |u(t,x)|^2 \d x + \epsilon \int_{\set{ x \in \Omega }{ |u(t,x)| \geq K }} |u(t,x)|^2 \d x \\
  &\leq ( \epsilon K^2 + M ) \Vol(\Omega) + \int_{\set{ x \in \Omega }{ |u(t,x)| \geq K }} f(x,u(t,x)) \eta(u(t,x)) \d x \\
  &\leq \left( \epsilon K^2 + M + \sup_{x \in \Omega,\; |v| \leq K} | f(x,v) \eta(v) | \right) \Vol(\Omega) + \langle f(\cdot,u(t,\cdot) , \eta(u(t,\cdot)) \rangle_{L^2(\Omega)}.
\end{align*}
Now multiply this inequality by the test function $\phi$ from theorem \ref{thm:solution_compactness_Dirichlet} and integrate over $t \in J$.
First we note that computations similar to the ones in the proof of theorem \ref{thm:solution_compactness_Dirichlet} yield
\[
\int_J \phi(t)  \langle f(\cdot,u(t,\cdot) , \eta(u(t,\cdot)) \rangle_{L^2(\Omega)} \d t
\leq C_\mu + \int_Q \phi | \nabla u |^2 \d x \d t + \mu \int_Q \phi | u |^2 \d x \d t
\]
for any $\mu > 0$.
We then find that
\[
  \epsilon \int_Q \phi |u|^2 \d x \d t
\leq C_{K,\mu} + \int_Q \phi |\nabla u |^2 \d x \d t + \mu \int_Q \phi |u|^2 \d x \d t,
\]
where $C_{K,\mu}$ is independent of $U = ( u , \partial_t u ) \in \calA_a^b(J,f,g,c)$.
Choosing $\mu$ sufficiently small, the last integral can be absorbed into the left hand side.
Thus we obtain the desired bound
\[
\int_Q \phi |u|^2 \d x \d t \leq C_1 + C_2 \int_Q \phi |\nabla u|^2 \d x \d t \qquad \text{for any} \quad U = (u,\partial_t u) \in \calA_a^b(J,f,g,c).
\]
This estimate can now replace the Poincar\'e inequality in the proof of theorem \ref{thm:solution_compactness_Dirichlet}.
The remainder of the proof of theorem \ref{thm:solution_compactness_Dirichlet} remains valid without further modifications.
\end{myproof}

\subsubsection{The nonautonomous case}
In the nonautonomous case we can now readily combine the ideas from the preceding paragraphs to conclude the following.
\begin{mythm}\label{thm:solution_compactness_Neumann_nonautonomous}
Consider Neumann or periodic boundary data.
Given $\theta \in (-1,1)$, $\gamma > 0$, and $\ell > 0$, there exists $\Theta = \Theta(\theta,\gamma,\ell) > 0$ for which the following is true.
Let $(f,g,c)$ be of class $C^m$ ($m \geq 1$) for which hypotheses \Cgroup are satisfied with the chosen constants $\theta$, $\ell$, $\Theta$,
and $\inf_{t\in\R} c(t) \geq \gamma$.
Fix $J = (j_-,j_+) \subset \R$, and $J' \subset \cl{J'} \subset J$ with $[-\ell,\ell] \subset J'$.
Then the set $\rst{ \calA_a^b(J,f,g,c) }{ J' }$ is bounded in $C_b^m(J';X^0,\dots,X^m)$ and compact in $W^{m,2}_{\text{loc}}(J';X^0,\dots,X^m)$.
\end{mythm}

\subsection{Regularity of the map \texorpdfstring{$\partial_t + A_{f,g,c}(\cdot)$}{dt + A}}

Consider the map
\[\
\partial_t + A_{f,g,c}(\cdot) : W^{1,2}_{\text{loc}}(J;X^0,X^1) \to L^2_{\text{loc}}(J;X^0).
\]
Smoothness of $g$ and hypothesis \ref{H:growth_bound_g} imply that the map $U \mapsto ( 0 , g(x,U(t)) )$ induces a smooth map
from $W^{1,2}_{\text{loc}}(J;X^0,X^1)$ into $L^2_{\text{loc}}(J;X^0)$.
Hence that the regularity class of this map is the same as that of the Nemytskii operator $f : H^1(J'' \times \Omega) \to L^2(J'' \times \Omega)$,
for bounded subsets $J'' \subset J$.
As we already remarked in the proofs of the compactness theorems, hypothesis \ref{H:growth_bound_f} implies that
the Nemytskii operator $f : H^1(J \times \Omega) \to L^2(J \times \Omega)$
is bounded and continuous, but in general does not possess additional regularity.
However, the failure to be more regular only stems from the behaviour of $f(t,x,u)$ for large $u$.
Since we have seen that $\rst{ \calA_a^b(J,f,g,c) }{ J' }$ is bounded in $C_b^1(J';X^0,X^1)$, the map does have additional regularity
near solutions of \eqref{eq:TWE}.
For clarity these observations are summarised in form of a theorem.
\begin{mythm}
\label{thm:operator_regularity}
  Suppose $f$ is of class $C^m$, with $m \geq 1$.
Given $J = (j_-,j_+) \subset \R$, and $J' \subset \cl{ J' } \subset J$, and $a , b \in \R$, the maps
\[
A_{f,g,c}(\cdot) : X^1 \to X^0
\]
and
\[
\partial_t + A_{f,g,c}(\cdot) : W^{1,2}_{\text{loc}}(J';X^0,X^1) \to L^2_{\text{loc}}(J';X^0)
\]
are bounded and of class $C^m$ in neighbourhoods of $\calS_a^b(f)$ and $\rst{ \calA_a^b(J,f,g,c) }{ J' }$, respectively.
\end{mythm}

\subsection{Energy bounds}

The following lemma shows that the energy is bounded from below on the collection $\calS$ of stationary solutions.
Observe that by the implicit function theorem hyperbolic stationary points are isolated in $\calS$ with respect to the topology of $X^1$.
Hence combining this energy bound with our compactness result shows that if \eqref{eq:TWE} is hyperbolic, 
then for arbitrary $a \in \R$ the collection $\calS_{-\infty}^a$ is a finite set.
We will need this fact later on when defining the Floer boundary operator.

\begin{mylemma}\label{lemma:energy_bounds_stationary}
Suppose either hypothesis \ref{H:contact_condition_odd} or \ref{H:contact_condition_even} is satisfied.
  There exists a constant $M \in \R$ such that $\calE(Z) \geq M$ whenever $Z \in \calS$.
In particular, $\calS_{-\infty}^a$ is finite for hyperbolic nonlinearities $f$.
\end{mylemma}
\begin{myproof}
  We will first prove the statement when hypothesis \ref{H:contact_condition_odd} is satisfied.
  For any $Z = (z,0) \in \calS$, we then have
  \begin{align*}
    \calE(Z) &\geq \int_\Omega \frac 1 2 | \nabla z(x) |^2 - | F(x,z(x)) | \d x \\
    &\geq \int_\Omega \frac 1 2 | \nabla z(x) |^2 - \frac{\theta}{2} f(x,z(x)) z(x) \d x - C_f' \Vol(\Omega) \\
    &\geq \int_\Omega \frac 1 2 | \nabla z(x) |^2 - \frac{|\theta|}{2} f(x,z(x)) z(x) \d x - C_f' \Vol(\Omega) \\
    &= \int_\Omega \frac 1 2 | \nabla z(x) |^2 + \frac{|\theta|}{2} \Delta z(x) z(x) \d x - C_f' \Vol(\Omega) \\
    &= \int_\Omega \frac{ 1 - |\theta| }{2} | \nabla z(x) |^2 \d x - C_f' \Vol(\Omega) \\
    &\geq - C_f' \Vol(\Omega),
  \end{align*}
where we used that $0 \leq |\theta| < 1$.

When on the other hand hypothesis \ref{H:contact_condition_even} holds, 
we let $\eta \in C^1(\R)$ be as defined in estimate \ref{est:primitive_bound_alt} in the proof of theorem \ref{thm:solution_compactness_Dirichlet}.
Then, after updating the constant $C_f'$, we have
\[
|F(x,u)| \leq C_f' + \frac \theta 2 f(x,u) \eta(u),
\]
and consequently it again holds that
\[
\calE(Z) \geq \int_\Omega \frac{ 1 - |\theta| |\eta'(z(x))| }{2} | \nabla z(x) |^2 \d x - C_f' \Vol(\Omega) \geq - C_f' \Vol(\Omega).
\]
\end{myproof}

%

\section{Unique continuation}
\label{sec:unique_continuation}

The initial value problem associated with \eqref{eq:TWE} is ill-posed.
However, in this section we will show that if a solution through some initial value exists, then it must be unique.
This implies that time shifting defines a dynamical system on the space of all heteroclinic solutions to \eqref{eq:TWE},
and this dynamical system posesses a strict Lyapunov function given by the energy functional.
For second order elliptic equations such a uniqueness result is relatively well known;
it follows for example from Aronszajn's unique continuation theorem (see \cite{aronszajn1956unique, aronszajn1962unique}).
However, the nonlocal term appearing in \eqref{eq:TWE} prohibits application of this theory.
Therefore we present a new continuation result, tailored towards \eqref{eq:TWE}.

\subsection{Carleman estimates}
Here we generalize the Carleman estimates (see \cite{carleman1939probleme}) for the Laplacian.
More precisely, instead of only considering compactly supported functions, we allow for Dirichlet, Neumann or periodic boundary conditions
in the variables which are not being controlled by the weight function.
We can do so at the expense of having the lower bound on the weight $\tau$ depend on $\|u(t,\cdot)\|_{L^2(\Omega)}$ and $\|\nabla u(t,\cdot)\|_{L^2(\Omega)}$.

For notational convenience we write $\nabla_{t,x} = (\partial_t , \nabla)$ and $\Delta_{t,x} = \partial_t^2 + \Delta$.

\begin{mylemma}
\label{lemma:carleman}
Let $\phi(t) := t + t^2 / 2$.
For each $0 < \epsilon < 1$ there exists $C > 0$ so that the following holds.
For each $u \in H^3( \R \times \Omega )$ with $B(u(t,\cdot)) = 0$ for each $t \in \R$,
and $\supp( u ) \subset (-\epsilon , \epsilon) \times \cl \Omega$
one has the following \emph{Carleman estimate}: 
there exists $\tau_0(u) > 0$ such that for all $\tau \geq \tau_0(u)$ one has
\[
\tau^4 \| e^{\tau \phi} u \|_{L^2(\R\times\Omega)}^2 + \tau \| e^{\tau \phi} \nabla_{t,x} u \|_{L^2(\R\times\Omega)}^2
\leq C \| e^{\tau \phi} \Delta_{t,x} u \|_{L^2(\R\times\Omega)}^2.
\]
Furthermore, if $\widetilde u \in H^3(\R\times\Omega)$ is another function satisfying the above mentioned hypotheses,
and for all $t \in \R$ one has
\[
\frac{ \| \nabla \widetilde u(t,\cdot) \|_{L^2(\Omega)} }{ \| \widetilde u(t,\cdot) \|_{L^2(\Omega)} }
\leq
\frac{ \| \nabla u(t,\cdot) \|_{L^2(\Omega)} }{ \| u(t,\cdot) \|_{L^2(\Omega)} },
\]
then $\tau_0(\widetilde u) \leq \tau_0(u)$.
\end{mylemma}
\begin{myproof}
  Let us abbreviate  $\| \cdot \| = \| \cdot \|_{L^2(\R\times\Omega)}$ and 
  $\langle\cdot,\cdot\rangle = \langle\cdot,\cdot\rangle_{L^2(\R\times\Omega)}$.
Set $v = e^{\tau \phi} u$ and observe that
\[
\| e^{\tau \phi} \partial_t u \|^2 = \| \partial_t v - \tau \dot \phi v \|^2 \leq 2 \| \partial_t v \|^2 + 2 \| \tau \dot \phi v \|^2
\leq 2 \|\partial_t v \|^2 + 2 (1+\epsilon)^2 \tau^2 \| v \|^2.
\]
Hence it suffices to see that for all $v \in H^3(\R\times\Omega)$ with $B(v) = 0$ and $\supp( v ) \subset (-\epsilon , \epsilon)\times\Omega$
\begin{equation}
  \label{eq:carleman_convexified}
  \tau^4 \| v \|^2 + \tau \| \nabla_{t,x} v \|^2 \leq C \| P v \|^2 \qquad \text{for all} \quad \tau \geq \tau_0(u),
\end{equation}
where $P = e^{\tau \phi} \Delta_{t,x} e^{-\tau \phi}$.

We now decompose $P$ into a symmetric part and an anti-symmetric part:
\begin{align*}
  P &= P^S + P^A, \\
  P^S &= \Delta_{t,x} + \tau^2 (1+t)^2 - \tau, \\
  P^A &= -2\tau (1+t) \partial_t.
\end{align*}
Then
\[
\| P v \|^2 = \| P^S v \|^2 + \| P^A v \|^2 + \langle [P^S , P^A] v , v \rangle
\]
where $[\cdot,\cdot]$ denotes the commutator bracket.
At this point we used that $v$ is of class $H^3$, so that it lies in the domain of definition of $[P^S,P^A]$.

We are now ready to make the estimates.
Note that
\begin{align*}
  \| P^S v \|^2 &= \| \Delta_{t,x} v \|^2 + 2 \langle \Delta_{t,x} v , ( \tau^2 (1+t)^2 - \tau ) v \rangle + \| ( \tau^2 (1+t)^2 - \tau ) v \|^2 \\
  &\geq 2 \langle \Delta_{t,x} v , ( \tau^2 (1+t)^2 - \tau ) v \rangle + \| ( \tau^2 (1+t)^2 - \tau ) v \|^2 \\
  &= 2 \tau \| \nabla_{t,x} v \|^2 - \frac 1 2 \| P^A v \|^2 - 2 \tau^2 \| (1+t) \nabla v \|^2  - 4 \langle \partial_t v , \tau^2 (1+t) v \rangle \\
  &\quad + \tau^4 \| (1+t)^2 v \|^2 - 2 \tau^3 \| (1+t) v \|^2 + \tau^2 \|v\|^2 \\
  &\geq \frac{(1-\epsilon)^4}{2} \tau^4 \| v \|^2 + 2 \tau \| \nabla_{t,x} v \|^2  - \frac 1 2 \| P^A v \|^2 \\
  &\quad - 2 \tau^2 (1+\epsilon)^2 \| \nabla v \|^2  - 4 \langle \partial_t v , \tau^2 (1+t) v \rangle
\end{align*}
for $\tau \geq 2 (1+\epsilon)^2 / (1-\epsilon)^4$.
Here we used that $\supp( v ) \subset (-\epsilon,\epsilon) \times \Omega$.
Hence
\begin{equation}
  \label{eq:37}
  \begin{split}
      \| P^S v \|^2 + \| P^A v \|^2 &\geq \frac{(1-\epsilon)^4}{2} \tau^4 \| v \|^2 + 2 \tau \| \nabla_{t,x} v \|^2  \\
      &\quad - 2 \tau^2 (1+\epsilon)^2 \| \nabla v \|^2  - 4 \langle \partial_t v , \tau^2 (1+t) v \rangle.
  \end{split}
\end{equation}

To estimate the last term in \eqref{eq:37}, note that
\[
[P^S,P^A] = -4 \tau \partial_t^2 + 4 \tau^3 (1+t)^2,
\]
Hence
\begin{align*}
  - 4 \langle \partial_t v , \tau^2 (1+t) v \rangle &= - 4 \langle \tau^{1/2} \partial_t v , \tau^{3/2} (1+t) v \rangle \\
  &\geq - 2 \tau \| \partial_t v \|^2 - 2 \tau^3 \| (1+t) v \|^2 \\
  &= 2 \tau \langle \partial_t^2 v , v \rangle - 2 \tau^3 \| (1+t) v \|^2 \\
  &= -\frac 1 2 \langle [P^S,P^A] v , v \rangle.
\end{align*}
Therefore, since 
\[
\langle [P^S,P^A] v , v \rangle = 4 \tau \| \partial_t v \|^2 + 4 \tau^3 \| (1+t) v \|^2 \geq 0,
\]
we find that
\begin{equation}
  \label{eq:38}
  \| P v \|^2 \geq \frac{(1-\epsilon)^4}{2} \tau^4 \| v \|^2 + 2 \tau \| \nabla_{t,x} v \|^2  - 2 \tau^2 (1+\epsilon)^2 \| \nabla v \|^2.
\end{equation}

To get rid of the last term in \eqref{eq:38} we need to take a $u$-dependent lower bound $\tau_0$.
Since $\| v \| = 0$ implies $\|\nabla v \| = 0$, for a fixed $v$ we can always find $\tau_0 \geq 2 (1+\epsilon)^2 / (1-\epsilon)^4$
such that
\begin{equation}
  \label{eq:28}
  \frac{(1-\epsilon)^4}{2} \tau^4 \| v \|^2 - 2 (1+\epsilon)^2 \tau^2 \| \nabla v \|^2 \geq \frac{ (1-\epsilon)^4 }{ 4 } \tau^4 \| v \|^2 
\qquad \text{for all} \quad \tau \geq \tau_0.
\end{equation}
If $\widetilde u$ is as in the hypotheses of the lemma, and $\widetilde v = e^{\tau\phi} \widetilde u$, 
then since the exponential factors through the inequalities, also
\[
\frac{ \| \nabla \widetilde v(t,\cdot) \| }{ \| \widetilde v(t,\cdot) \| }
\leq
\frac{ \| \nabla v(t,\cdot) \| }{ \| v(t,\cdot) \| }
\]
for all $t \in \R$.
Therefore, if $\tau_0(u)$ denotes the smallest constant $\tau_0$ for which \eqref{eq:28} holds, it is readily seen that $\tau_0(\widetilde u) \leq \tau_0(u)$.
Hence
\[
\| P v \|^2 \geq \frac{ (1-\epsilon)^4 }{ 4 } \tau^4 \|  v \|^2 + 2 \tau \| \nabla_{t,x} v \|^2 \qquad \text{for all} \quad \tau \geq \tau_0(u),
\]
from which \eqref{eq:carleman_convexified} follows.
\end{myproof}

\subsection{Continuation for an integro-differential inequality}
The following lemma is in a sense an integrated version of Aronszajn's continuation theorem.

\begin{mylemma}
\label{lemma:ineq_continuation}
Let $J \subset \R$ be an open interval,
and let $u \in H^3(J\times\Omega)$ with $B(u(t,\cdot)) = 0$ for all $t \in J$.
Assume that it satisfies the integro-differential inequality
  \begin{equation}
    \label{eq:integral_inequality}
    \int_\Omega | \Delta_{t,x} u(t,x) |^2 \d x \leq C \int_\Omega | u(t,x) |^2  + | \nabla_{t,x} u(t,x) |^2 \d x,
  \end{equation}
for almost every $t$ in some neighbourhood of $t_0$, for a certain $t_0 \in J$.
Assume furthermore that $u$ satisfies the following decay conditions around $t_0$:
\begin{equation}
  \label{eq:decay_conditions}
  \begin{split}
      &\int_{t_0-\delta}^{t_0+\delta} \int_\Omega |u(t,x)|^2 \d x \d t = O(\delta^5) \qquad \text{as} \quad \delta \downarrow 0, \\
      &\int_{t_0-\delta}^{t_0+\delta} \int_\Omega | \partial_t u(t,x) |^2 \d x \d t = O(\delta^3) \qquad \text{as} \quad \delta \downarrow 0.
  \end{split}
\end{equation}
Then $u \equiv 0$ in a neighbourhood of $\{t_0\}\times\Omega$.
\end{mylemma}
\begin{myproof}
  The strategy is as follows.
First we note there is no loss of generality in assuming that $t_0 = 0$.
We localize on the left side of the hyperplane $t=0$.
The Carleman estimates are still valid for these localized solutions, at the cost of some error terms.
One of the error terms stems from the behaviour of our localized solution away from $t=0$.
This error term can be made to decay at an exponential rate, 
precisely because the localized solution vanishes on the right side of the hyperplane $t=0$.
The other error terms stem from the behaviour of the localization near $t=0$.
To deal with these terms, we actually consider a sequence of localizations.
The decay conditions \eqref{eq:decay_conditions} allow us to take the limit in which the error terms disappear.
Combining this extension of the Carleman estimates with estimate \eqref{eq:integral_inequality},
we are left with a family of exponentially weighted inequalities, which forces the localization to be zero near $t=0$.
This just means that our original function must be zero for small negative $t$.
Then considering a time reversal, the same must be true for small positive $t$.

We will now first construct the sequence of localizations.
Let $\epsilon > 0$ be sufficiently small such that $(-\epsilon,\epsilon) \subset J$,
and let $\phi$ be as in lemma \ref{lemma:carleman}.
Let $0 < \ell < \epsilon / 2$ be such that $\phi$ is increasing on $[-2\ell,0]$,
and such that \eqref{eq:integral_inequality} holds for a.e. $t \in (-2\ell,2\ell)$.
Given $0 < \delta < \ell$, let $\chi_\delta \in C^3(\R)$ be such that 
$\chi_\delta(t) = 0$ for $t \not\in [-2\ell,0]$, and $\chi_\delta(t) = 1$ for $t \in [-\ell,-\delta]$.
Furthermore, we demand that $|\partial_t^k \chi_\delta(t)| = O(\delta^{-k})$ uniformly for $t \in [-\delta,0]$ and $k \in \{0,1,2\}$,
and also that $\partial_\delta \chi_\delta(t) = 0$ for $t \in [-2\ell,-\ell]$.
Now set $V_\delta(t,x) := \chi_\delta(t) u(t,x)$.
Then $V_\delta \in H^3(\R \times \Omega)$, and moreover $B(V_\delta(t,\cdot)) = 0$ for each $t \in \R$, 
and $\supp( V_\delta ) \subset (-\epsilon,\epsilon) \times \cl\Omega$, hence the Carleman estimates apply to $V_\delta$.

We will again abbreviate $\| \cdot \| = \| \cdot \|_{L^2(\R\times\Omega)}$.
Observe
\begin{align*}
  \| \chi_\delta e^{\tau\phi} \partial_t u \|^2 &= \| e^{\tau\phi} \partial_t V_\delta - \dot\chi_\delta e^{\tau\phi} u \|^2 \\
&\leq 2 \| e^{\tau\phi} \partial_t V_\delta \|^2 + 2 \| \dot\chi_\delta e^{\tau\phi} u \|^2 \\
&\leq 2 \| e^{\tau\phi} \partial_t V_\delta \|^2 
+ 2 \int_{-2\ell}^{-\ell} \int_\Omega | \dot\chi_\delta |^2 e^{2 \tau\phi} |u|^2 \d x \d t \\
&\quad + 2 \sup_{ -\delta \leq t \leq 0 } | \dot\chi_\delta(t) |^2 \int_{-\delta}^0 \int_\Omega |u|^2 \d x \d t \\
&\leq 2 \| e^{\tau\phi} \partial_t V_\delta \|^2 
+ 2 e^{2\tau\phi(-\ell)} \sup_{-2\ell \leq t \leq -\ell} |\dot\chi_\delta(t)| \int_{-2\ell}^{-\ell} \int_\Omega |u|^2 \d x \d t
+ O(\delta^3),
\end{align*}
where we used the decay condition \eqref{eq:decay_conditions} and the monotonicity of $\phi$.
Thus we have
\begin{equation}
  \label{eq:52}
  \sum_{|\alpha|\leq 1} \| \chi_\delta e^{\tau\phi} \partial^\alpha u \|^2 \leq 2 \sum_{|\alpha|\leq 1} \| e^{\tau\phi} \partial^\alpha V_\delta \|^2 
+ C_1 e^{2\tau\phi(-\ell)} + O(\delta^3).
\end{equation}
Here $C_1$ depends on $\|u\|$ but is independent of $\delta$ and $\tau$.
We use the multi-index notation $\alpha = (\alpha_0,\dots,\alpha_d) \in \N_0^{d+1}$, and $|\alpha|=\alpha_0 + \cdots + \alpha_d$, and
\[
\partial^\alpha = \partial_t^{\alpha_0} \partial_{x^1}^{\alpha_1} \cdots \partial_{x^d}^{\alpha_d}.
\]

In a similar fashion, we compute
\begin{equation}
  \label{eq:53}
  \begin{split}
    \| e^{\tau\phi} \Delta_{t,x} V_\delta \|^2 &\leq 3 \| \chi_\delta e^{\tau\phi} \Delta_{t,x} u \|^2 
    + 3 \| 2 \dot\chi_\delta e^{\tau\phi} \partial_t u \|^2 + 3 \| \ddot \chi_\delta e^{\tau\phi} u \|^2  \\
    &\leq 3 \| \chi_\delta e^{\tau\phi} \Delta_{t,x} u \|^2 
    + 3 \int_{-2\ell}^{-\ell} \int_\Omega e^{2\tau\phi} ( | 2 \dot\chi_\delta \partial_t u |^2 + | \ddot \chi_\delta u |^2 ) \d x \d t \\
    &\quad + 12\sup_{-\delta\leq t\leq 0} | \dot\chi_\delta(t) |^2 \int_{-\delta}^0 \int_\Omega | \partial_t u |^2 \d x \d t \\
    &\quad + 3 \sup_{-\delta\leq t \leq 0} |\ddot \chi_\delta(t)|^2 \int_{-\delta}^0 \int_\Omega |u|^2 \d x \d t \\
    &\leq 3 \| \chi_\delta e^{\tau\phi} \Delta_{t,x} u \|^2 + C_2 e^{2\tau \phi(-\ell)} + O(\delta),
  \end{split}
\end{equation}
where we again used the decay condition \eqref{eq:decay_conditions} and the monotonicity of $\phi$.
The constant $C_2$ depends on $\|u\|$ and $\|\partial_t u\|$ but is independent of $\delta$ and $\tau$.

Combining estimates \eqref{eq:52} and \eqref{eq:53} with the Carleman estimates from lemma \ref{lemma:carleman}, we find that
\begin{equation}
  \label{eq:54}
  \tau \sum_{|\alpha|\leq 1} \| \chi_\delta e^{\tau\phi} \partial^\alpha u \|^2 \leq C_3 \| \chi_\delta e^{\tau\phi} \Delta_{t,x} u \|^2 + C_3 (1 + \tau) e^{2\tau\phi(-\ell)} 
+ O(\delta) + O(\tau \delta^3)
\end{equation}
for all $\tau \geq \tau_0(V_\delta)$.
Since $|V_\delta(t,x)| / |\nabla V_\delta(t,x)| = |V_{\delta'}(t,x)| / |\nabla V_{\delta'}(t,x)|$ for any two $0 < \delta , \delta' < \ell$,
from lemma \ref{lemma:carleman} it follows that the lower bound $\tau_0 := \tau_0(V_\delta)$ on $\tau$ is independent of $\delta$.
Therefore, using the dominated convergence theorem, we can send $\delta$ to $0$ in equation \eqref{eq:54}, and obtain
\begin{equation}
  \label{eq:perturbed_carleman}
  \tau \sum_{|\alpha|\leq 1} \| \chi_0 e^{\tau\phi} \partial^\alpha u \|^2 \leq C_3 \| \chi_0 e^{\tau\phi} \Delta_{t,x} u \|^2 + C_3 (1 + \tau) e^{2\tau\phi(-\ell)} 
\qquad \text{for all} \quad \tau \geq \tau_0.
\end{equation}
Here
\[
\chi_0(t) =
\begin{cases}
  1 & \text{for } t \in [-\ell,0], \\
  \chi_1(t) & \text{for } t \not\in [-\ell,0].
\end{cases}
\]

Using \eqref{eq:integral_inequality} we have
\begin{align*}
    \| \chi_0 e^{\tau\phi} \Delta_{t,x} u \|^2 &\leq C \int_\R |\chi_0(t)|^2 e^{2\tau\phi(t)} \int_\Omega | u(t,x) |^2 + \| \nabla_{t,x} u(t,x) \|^2 \d x \d t \\
  &= C \| \chi_0 e^{\tau\phi} u \|^2 + C \| \chi_0 e^{\tau\phi} \nabla_{t,x} u \|^2.
\end{align*}
Combining this inequality with \eqref{eq:perturbed_carleman} yields
\[
\tau \sum_{|\alpha|\leq 1} \| \chi_0 e^{\tau\phi} \partial^\alpha u \|^2 
\leq C_4 \sum_{|\alpha|\leq 1} \| \chi_0 e^{\tau\phi} \partial^\alpha u \|^2 + C_3 (1 + \tau) e^{2\tau\phi(-\ell)} 
\qquad \text{for all} \quad \tau \geq \tau_0.
\]
After increasing $\tau_0$ if need be, the sum over $|\alpha|\leq 1$ on the right hand side can be absorbed by the left hand side, hence
\[
\tau \sum_{|\alpha|\leq 1} \| \chi_0 e^{\tau\phi} \partial^\alpha u \|^2 \leq C_5 (1 + \tau) e^{2\tau\phi(-\ell)} 
\qquad \text{for all} \quad \tau \geq \tau_0.
\]
Since $\chi_0 = 1$ on $[-\ell,0]$ and $\phi$ is increasing on $[-\ell,0]$, this implies that
\[
\| u \|_{L^2([t_*,0]\times\Omega)}^2 \leq C_5 \frac{1+\tau}{\tau} e^{2\tau (\phi(-\ell) - \phi(t_*)) } \to 0 \qquad \text{as} \quad \tau \to\infty,
\]
for any $t_* \in (-\ell,0]$.
Hence $u \equiv 0$ on $(-\ell,0] \times \Omega$.
\end{myproof}

\subsection{Uniqueness of the IVP}

\begin{mythm}\label{thm:unique_continuation}
  Let $J \subset \R$ be an open interval, assume $f$ is of class $C^3$,
and suppose $U, V \in W^{1,2}_{\text{loc}}(J;X^0,X^1)$ are both solutions of \eqref{eq:TWE}.
Suppose that $U(t_0) = V(t_0)$ for some $t_0 \in J$.
Then $U(t) = V(t)$ for $t \in J$.
\end{mythm}
\begin{myproof}
Let $J' \subset \cl{J'} \subset J$ be a bounded open interval such that $t_0 \in J'$.
Let us introduce the set
\[
\calZ(J') := \set{ t \in J' }{ \| U(t) - V(t) \|_{X^0} = 0 }.
\]
By assumption $\calZ(J') \neq \emptyset$.
By the regularity theory from section \ref{sec:compactness} we know that $U, V \in C^3(J';X^0,\dots,X^m)$, hence $\calZ(J')$ is closed in $J'$.
Thus by connectedness of $J'$ we can conclude that $\calZ(J') = J'$ if we are able to prove that $\calZ(J')$ is open in $J'$.
Since $J$ can be written as the union of bounded open intervals $J' \subset \cl{J'} \subset J$, the conclusion of the theorem will then follow.

Pick any $t_* \in \calZ(J')$.
Writing $U = (u,\partial_t u)$ and $V = (v , \partial_t v)$, and set $W = (w,\partial_t w) := U - V$.
We will prove that $w$ is zero in a neighbourhood of $\{t_*\} \times \Omega$.
To do so, we shall invoke lemma \ref{lemma:ineq_continuation}.
We thus have to check that $w$ satisfies the hypotheses of said lemma.

Note that by the regularity theory of section \ref{sec:compactness}, $w \in H^3(J'\times\Omega)$, and $B(w(t,\cdot)) = 0$ for each $t \in J'$, and $w$ satisfies
\begin{align*}
  \Delta_{t,x} w &= f(x,v) - f(x,u) + c \partial_t w + g(x,V) - g(x,U) \\
  &= \int_0^1 \partial_s f(x,u - sw) \d s + c \partial_t w + \int_0^1 \partial_s g(x,U - sW) \d s \\
  &= \alpha(t,x) w + c \partial_t w + L(t)[ W(t) ].
\end{align*}
Here $\alpha : J' \times \Omega \to \R$ is given by
\[
  \alpha(t,x) = - \int_0^1 f_u(x,u(t,x)-sw(t,x)) \d s,
\]
and note that by the regularity results from section \ref{sec:compactness} this $\alpha$ is a continuous function.
Furthermore, $L : J' \to \calL( X^0 , L^2(\Omega) )$ is given by
\[
\big( L(t) \xi \big)(x) = - \int_0^1 \dsub 2 g(x,U(t) - s W(t))\xi \d s.
\]
We note here that this integral is indeed well defined, since the map
\[
s \mapsto I(s) := \dsub 2 g(x,U(t) - s W(t))
\]
is continuous from $[0,1]$ to $\calL( X^0 , L^2(\Omega) )$ with its uniform operator topology.
Hence $s \mapsto I(s)$ is absolutely continuous and therefore strongly measurable, 
and since by hypothesis \ref{H:growth_bound_g} one has $\int_0^1 \| I(s) \|_{\calL} \d s < \infty$, it follows that $I$ is Bochner integrable.
By the regularity results from section \ref{sec:compactness} we know that 
$U$ and $V$ form continuous curves in $X^0$, hence $t \mapsto L(t)$ is continuous.
Hence $\| L(t) \|_{\calL(X^0,L^2(\Omega))} \leq C$ uniformly for $t$ in a neighbourhood of $t_*$,
and consequently $w$ satisfies an inequality of the form \eqref{eq:integral_inequality}.

All that is left is to check that $w$ satisfies the decay conditions \eqref{eq:decay_conditions} around $t_*$.
Let $\delta > 0$ be sufficiently small such that $[t_*-\delta , t_*+\delta] \subset J'$.
Since $f$ is of class $C^3$, it follows from section \ref{sec:compactness} 
that $W = (w,\partial_t w) \in C_b^3([t_*-\delta,t_*+\delta];X^0,\dots,X^3)$,
hence in particular $w \in C_b^4([t_*-\delta,t_*+\delta];L^2(\Omega))$.
Now consider the function $\eta : (t_*-\delta,t_*+\delta) \to \R$ given by
\[
\eta(t) = \int_\Omega | w(t,x) |^2 \d x.
\]
It is $C^4$ and the first three derivatives are given by
\begin{align*}
  \eta'(t) &= 2 \int_\Omega w(t,x) \partial_t w(t,x) \d x, \\
  \eta''(t) &= 2 \int_\Omega w(t,x) \partial_t^2 w(t,x) \d x + 2 \int_\Omega  |\partial_t w(t,x)|^2 \d x, \\
  \eta'''(t) &= 2 \int_\Omega w(t,x) \partial_t^3 w(t,x) \d x + 6 \int_\Omega \partial_t w(t,x) \partial_t^2 w(t,x) \d x.
\end{align*}
Since $w(t_*,\cdot) = \partial_t w(t_*,\cdot) = 0$, we have $\eta^{(k)}(t_*) = 0$ for $k \in \{0,\dots,3\}$.
By the mean value theorem it follows that $|\eta(t)| \leq C' \delta^4$, hence
\[
  \int_{t_*-\delta}^{t_*+\delta} \int_\Omega |w(t,x)|^2 \d x \d t \leq C \delta^5 \qquad \text{for} \quad |t-t_*| \leq \delta.
\]
A similar computation shows that
\[
\int_{t_*-\delta}^{t_*+\delta} \int_\Omega |\partial_t w(t,x)|^2 \d x \d t \leq C \delta^3 \qquad \text{for} \quad |t-t_*| \leq \delta.
\]
Hence lemma \ref{lemma:ineq_continuation} applies, proving that $w \equiv 0$ in a neighbourhood of $\{t_*\}\times\Omega$.
Therefore $U(t) = V(t)$ for $t$ in a neighbourhood $E \subset J'$ of $t_*$.
Hence $E \subset \calZ(J')$, and since this holds for any $t_* \in \calZ(J')$ it follows that $\calZ(J')$ is open in $J'$, thus proving the theorem.
\end{myproof}

\section{Fredholm theory}
\label{sec:fredholm}

Let $\frakL$ consist of all $L \in \calL(X^1,X^0)$ of the form
\[
L =
\begin{pmatrix}
  0 & -1 \\
  \Delta + L_1 & L_2,
\end{pmatrix}.
\]
where $L_1 \in \calL(H_B^1(\Omega),L^2(\Omega))$ and $L_2 \in \calL(L^2(\Omega))$.
Let $\frakL_{\text{hyp}}$ consist of those $L \in \frakL$ which are hyperbolic, i.e.\  $\sigma(L) \cap i\R = \emptyset$.

In this section we will study Fredholm properties of the linear operator
\begin{align*}
    &\calD_L : W^{1,2}(\R;X^0,X^1) \to L^2(\R;X^0), \\
    &\calD_L W = \partial_t W + L(t) W,
\end{align*}
where $L \in C^0(\R;\frakL)$ is such that the limits $L_\pm = \lim_{t\to\pm\infty} L(t)$ exist in the uniform operator topology on $\calL(X^1,X^0)$,
and $L_\pm \in \frakL_{\text{hyp}}$.
The study of Fredholm properties of this class of operators is motivated by the following lemma.
\begin{mylemma}
\label{lemma:abstract_Fredholm_application}
  Let $Z_- \in \calS(f_-)$ and $Z_+ \in \calS(f_+)$ be hyperbolic rest points of \eqref{eq:TWE},
and suppose $U$ is a path connecting $Z_-$ with $Z_+$.
Then the linearization $\partial_t + \d A_{f,g,c}(U)$ of \eqref{eq:TWE} along $U$ is Fredholm, with index given by $\ind(\calD_{L_{c,f}})$,
where
\[
L_{c,f}(t) =
\begin{pmatrix}
  0 & -1 \\
  \Delta + f_u(t,x,u(t,x)) & -c(t)
\end{pmatrix}.
\]
\end{mylemma}
\begin{myproof}
  
Note that the linearization of \eqref{eq:TWE} along $U$ takes the form
\begin{equation}
  \label{eq:26}
  \partial_t W + L_{c,f}(t) W + K(t) W,
\end{equation}
where
\[
K(t) =
\begin{pmatrix}
  0 & 0 \\
  \partial_1 g(x,U(t)) & \partial_2 g(x,U(t))
\end{pmatrix}.
\]
Here $\partial_1 g(x,(u,v)) = \ddi{ g(x,(u,v)) }{ u }$ and $\partial_2 g(x,(u,v)) = \ddi{ g(x,(u,v)) }{ v }$.
For each $t$, $K(t) \in \calL(X^0)$, hence $K(t) : X^1 \to X^0$ is compact.
Furthermore, $K \in C^0( \R , \calL(X^1,X^0) )$ and $\| K(t) \|_{\calL(X^1,X^0)} \to 0$ as $t \to \pm\infty$ by hypothesis \ref{H:vanishing_condition}.
This implies that the multiplication operator $K : W^{1,2}(\R;X^1,X^0) \to L^2(\R;X^0)$ is compact, see \cite{robbin1995spectral}.
Therefore the Fredholm properties of \eqref{eq:26} are the same as those of $\calD_{L_{c,f}}$.
\end{myproof}

The Fredholm properties which will be derived allow us to assign a (normalized) Morse index to hyperbolic rest points,
even though the operators $\d A_{c_-,f_-,g_-}(Z_-)$ and $\d A_{c_+,f_+,g_+}(Z_+)$ are strongly indefinite.

\subsection{Fredholm alternative for \texorpdfstring{$\calD_L$}{DL}}

Before discussing the Fredholm alternative for $\calD_L$, let us first consider a resolvent estimate for the operator $L(t)$.

\begin{mylemma}\label{lemma:resolvent_decay}
  Let $L \in C^0(\R;\frakL)$ such that the limits $L_\pm = \lim_{t\to\pm\infty} L(t)$ exist in the uniform operator topology on $\calL(X^1,X^0)$.
  Then one has the following resolvent estimate: there exist $M > 0$, $R_0 > 0$ such that
  \[
  \| \lambda ( L(t) - i\lambda )^{-1} \|_{\calL(X^0)} \leq M \qquad \text{for} \quad t \in \R,\; |\lambda| \geq R_0.
  \]
\end{mylemma}
\begin{myproof}
First consider the unbounded operator $P$ on $X^0$ with domain $\calD(P) = X^1$, given by
\[
P =
\begin{pmatrix}
  0 & -1 \\
  \Delta & 0
\end{pmatrix}.
\]
Note that $i \lambda \not\in \sigma(P)$ whenever $\lambda \in \R \setmin \{0\}$, and
\[
( P - i\lambda )^{-1} =
\begin{pmatrix}
  i \lambda^{-1} \big( 1 - (\Delta - \lambda^2)^{-1} \big) & (\Delta - \lambda^2)^{-1} \\
  - (\Delta - \lambda^2)^{-1} & -i \lambda (\Delta - \lambda^2)^{-1}
\end{pmatrix}.
\]
Now, since
\[
\| (\Delta - \mu)^{-1} \|_{\calL( L^2(\Omega) , H_B^k(\Omega))} \leq \frac{C}{\bigg( 1 + d\big( \mu , \sigma(\Delta) \big) \bigg)^{(2-k)/2}}
\]
for $k \in \{0,1,2\}$,
we find that
\[
\| ( P - i\lambda )^{-1} \|_{X^0} \leq \frac{C}{1 + |\lambda|} \qquad \text{for} \quad \lambda \in \R \setmin \{0\}.
\]

Now let $K(t)$ be defined by
\[
K(t) :=
\begin{pmatrix}
  0 & 0 \\
  L_1(t) & L_2(t)
\end{pmatrix},
\]
so that $L(t) = P + K(t)$.
Note that $K(t)$ is a bounded operator on $X^0$, and $K \in C^0(\R,\calL(X^0))$.
Using a perturbative argument (see e.g.\  \cite{kato2012perturbation})
one has that 
\begin{equation}
  \label{eq:perturbation_estimate}
  i\lambda \not\in \sigma(L(t)) \qquad \text{when} \quad \| K(t) \|_{\calL(X^0)} \| (P - i\lambda)^{-1} \|_{\calL(X^0)} < 1,
\end{equation}
and for such $\lambda$ one has
\[
\| ( L(t) - i \lambda )^{-1} \|_{\calL(X^0)} \leq \frac{1}{1 - \| K(t) \|_{\calL(X^0)} \|  ( P - i\lambda )^{-1} \|_{\calL(X^0)} } \|  ( P - i\lambda )^{-1} \|_{\calL(X^0)}.
\]
We have already argued that $\|  ( P - i\lambda )^{-1} \|_{\calL(X^0)} = O(|\lambda|^{-1})$ as $|\lambda| \to \infty$.
Hence for each fixed $t \in \R$, condition \eqref{eq:perturbation_estimate} is satisfied for $|\lambda|$ sufficiently big, say $|\lambda| \geq R_0(t)$.
Since $K$ depends continuously on $t$, and $K(t)$ converges as $t \to \pm\infty$, the constant $R_0$ can be chosen uniformly in $t$.
\end{myproof}

Combining this lemma with the results from \cite{rabier2004robbin}, we obtain the following theorem
(see appendix \ref{sec:appendix_fredholm} for details).
\begin{mythm}\label{thm:fredholm_alternative}
Let $L \in C^0(\R;\frakL)$ be such that $L(t) \to L_\pm$ as $t \to \pm\infty$ in the uniform operator topology on $\calL(X^1,X^0)$, where $L_\pm \in \frakL_{\text{hyp}}$.
  Then the map $\calD_L$ is Fredholm from $W^{1,2}(\R;X^0,X^1)$ to $L^2(\R;X^0)$, and its index depends on the endpoints $L_-$, $L_+$ only.
\end{mythm}
This allows us to define a relative index:
\begin{align*}
  & \nu : \frakL_{\text{hyp}} \times \frakL_{\text{hyp}} \to \Z, \\
  & \nu(L_-,L_+) = \ind( \calD_L ).
\end{align*}

The relative index has the following transitivity property.
\begin{mylemma}\label{lemma:cyclicity_index}
  Let $L_\alpha, L_\beta, L_\gamma \in \frakL_{\text{hyp}}$.
  Then
  \[
     \begin{array}{l l}
      \nu(L_\alpha,L_\beta) = - \nu(L_\beta,L_\alpha) & \quad \text{(antisymmetry)}, \\
      \nu(L_\alpha , L_\gamma) = \nu(L_\alpha, L_\beta) + \nu(L_\beta, L_\gamma) & \quad \text{(cyclicity)}.
    \end{array}
  \]
\end{mylemma}
The proof of the lemma uses an algebraic trick similar to the one employed in \cite{robbin1995spectral}.
See appendix \ref{sec:appendix_fredholm} for details.

\subsection{Computing the index}
Consider $L \in C_b^0(\R,\calL)$ with $\lim_{t\to\pm\infty} L(t) \in \calL_{\text{hyp}}$, given by
\begin{equation}
  \label{eq:L_class}
  L(t) = 
\begin{pmatrix}
  0 & -1 \\
  \Delta + L_1(t) & -c(t)
\end{pmatrix}.
\end{equation}
Here we assume $c \in C^\infty(\R,(0,\infty))$,
and $L_1 \in C^0(\R,\calL( H_B^1(\Omega),L^2(\Omega)))$,
and for each $t \in \R$ the operator $L_1(t)$ is symmetric when considered as an unbounded operator on $L^2(\Omega)$.
Consider the operator
\begin{align*}
    &\Psi_{L_1} : W^{1,2}(\R;H^1_B(\Omega),H^2_B(\Omega)) \to L^2(\R;L^2(\Omega)), \\
    &\Psi_{L_1} w = \partial_t w + \Delta w + L_1(t) w.
\end{align*}
The following lemma relates the Fredholm index of the elliptic operator $\calD_L$ with that of the parabolic operator $\Psi_{L_1}$.
\begin{mylemma}
\label{lemma:index_computation}
  Let $L$, $L_1$ be as in \eqref{eq:L_class}.
  Then $\calD_L$ and $\Psi_{L_1}$ are Fredholm operators, and
  \[
  \ind(\calD_L) = \ind(\Psi_{L_1}).
  \]
\end{mylemma}
\begin{myproof}
The Fredholm property of $\calD_L$ was discusses in theorem \ref{thm:fredholm_alternative}.
Moreover, theorem \ref{thm:rabier} is also applicable to the operator $\Psi_{L_1}$, thus establishing the Fredholm property of that operator as well.
To relate the two indices we will make use of spectral flows.

Loosely speaking, the spectral flow $\SF(A)$ of a curve of (densely defined unbounded) operators $t \mapsto A(t)$ is an algebraic count of the number of eigenvalues of $A(t)$ that
cross the imaginary axis as $t$ increases from $-\infty$ to $+\infty$.
More precisely, we define
\[
\SF(A) := - \sum_{t_*} \sum_i \sgn \real \lambda_i'(t_*).
\]
Here $t_*$ are those $t$ where a spectral crossing takes place, i.e.\  $\sigma(A(t_*)) \cap i\R \neq \emptyset$,
and $t \mapsto \lambda_i(t)$ are differentiable curves defined near $t_*$ which parametrize the spectrum,
i.e.\  $\sigma(A(t)) = \bigcup_i \{ \lambda_i(t) \}$ for $t$ near $t_*$.
This definition only makes sense if there is no ambiguity in counting the crossings of eigenvalues:
the operators $A(t)$ should have pure point spectrum near the imaginary axis,
there should be finitely many crossings, the crossings should be transverse to $i\R$, 
the crossing eigenvalues should be algebraically simple.
This is generically achieved by perturbing the curve $t \mapsto A(t)$ of operators.

Recall from \cite{robbin1995spectral} that we can choose a perturbation $S_1 \in C^0(\R,\calL( H_B^2(\Omega),L^2(\Omega)))$,
with $\lim_{t\to\pm\infty} \| S(t) \|_{\calL(X^1,X^0)} = 0$, such that $\SF(- (\Delta + L_1 + S_1))$ is well-defined, and
\[
\ind( \Psi_{L_1} ) = \SF(- (\Delta + L_1 + S_1)).
\]
Furthermore, one can ensure that for any $t \in \R$ the operator $S_1(t)$ is symmetric when considered as an unbounded operator on $L^2(\Omega)$,
and given $\epsilon > 0$ (to be specified in the next paragraph), $\sup_{t\in\R} \| S_1(t) \|_{\calL(X^1,X^0)} < \epsilon$.

Now let
\[
S(t) =
\begin{pmatrix}
  0 & 0 \\
  S_1(t) & 0
\end{pmatrix}.
\]
Suppose for the moment that $\SF(- (L + S))$ is well-defined.
It then follows from \cite{rabier2004robbin} that, provided that $\sup_{t\in\R} \| S(t) \|_{\calL(X^1,X^0)}$ is sufficiently small, one has
\[
\ind( \calD_L ) = \SF(- (L + S)).
\]
The proof of the lemma is then completed if we can show that $\SF(- (L+S)) = \SF(- (\Delta + L_1 + S_1))$.
We will do so by demonstrating that any spectral crossing of $-( \Delta + L_1 + S_1 )$ is in one-to-one correspondence with a spectral crossing of $- (L + S)$.
This then also shows that $\SF(-(L + S))$ is indeed well-defined.

Note that $\mu \in \sigma(-(L(t) + S(t)))$ if and only there exists a nonzero vector $(u,v) \in X^1$ such that
\[
\begin{cases}
  - \mu u + v = 0, \\
  -\Delta u - L_1(t) u - S_1(t) u + c(t) v - \mu v = 0.
\end{cases}
\]
This can only happen if $\mu^2 - c(t) \mu \in \sigma(-( \Delta + L_1(t) + S_1(t)))$.
Because the operator $-(\Delta + L_1(t) + S_1(t))$ is self-adjoint, it has real-valued spectrum.
Since $c(t) > 0$, the only $\mu \in i\R$ which can satisfy $\mu^2 - c(t) \mu \in \R$ is $\mu = 0$.
Collecting these observations, we conclude that eigenvalues of $-( L(t) + S(t) )$ which cross
the imaginary axis as $t$ increases must do so through the origin.
Furthermore, if an eigenvalue of $- ( L(t) + S(t))$ crosses the imaginary axis (i.e.\  the origin) from left to right as $t$ increases, 
then (since $c(t_*) > 0$) an eigenvalue of $- (\Delta + L_1(t) + S_1(t))$ passes the origin from left to right as $t$ increases, and vice versa.
From this we conclude that the spectral flow of $-(L + S)$ is well-defined for generic $S_1$, and
\[
\SF(-(L + S)) = \SF(-(\Delta + L_1 + S_1)).
\]
Combined with our previous observations this complete the proof.
\end{myproof}

\subsection{Normalized Morse indices}

Given a regular nonlinearity $f$ satisfying hypothesis \ref{H:growth_bound_f}, and $m_0 \in \Z$,
we define a normalized Morse index $\mu_f(Z)$ of $Z = (z,0) \in \calS(f)$ as
\begin{equation}
  \label{eq:normalized_index_def}
  \begin{split}
    \mu_f(Z) := m_0 - m_f(z), \qquad \text{where} \quad m_f(z) := \# \big( \sigma(\Delta + f_u(x,z) ) \cap (0,\infty) \big).
  \end{split}
\end{equation}
See remark \ref{remark:index_convention} after the next theorem for the rationale behind the sign convention.
Whenever the choice of nonlinearity $f$ is clear from the context we shall drop the subscript.
\begin{mythm}
    Let $Z_- \in \calS(f_-)$ and $Z_+ \in \calS(f_+)$ be hyperbolic rest points of \eqref{eq:TWE},
and suppose $U$ is a path connecting $Z_-$ with $Z_+$.
Then
\[
\ind( \partial_t + \d A_{f,g,c}(U) ) = \mu_{f_-}(Z_-) - \mu_{f_+}(Z_+).
\]
\end{mythm}
\begin{myproof}
Define the multiplication operator
\[
L_1(t) : H^1_B(\Omega) \to L^2(\Omega), \qquad w(x) \mapsto f_u(x,u(t,x)) w(x).
\]
Combining lemmata \ref{lemma:abstract_Fredholm_application} and \ref{lemma:index_computation}, we find
\[
\ind( \partial_t + \d A_{f,g,c}(U) ) = \ind( \Psi_{L_1} ).
\]
From \cite{robbin1995spectral} it follows that
\[
\ind( \Psi_{L_1} ) = m_{f_+}( z_+ ) - m_{f_-}(z_-).
\]
Combined with the definition of $\mu_f$ this concludes the proof.
\end{myproof}

\begin{myremark}
  \label{remark:index_convention}
The reason for choosing a minus sign in the definition of the normalized Morse index is to ensure that the index can only decrease along heteroclinic connections.
One can drop the minus sign in \eqref{eq:normalized_index_def} and arrive at a cohomology theory instead.
\end{myremark}

\section{Exponential decay}
\label{sec:exp_decay}

In this section we will show that a solution of \eqref{eq:TWE} which converges in forward time (or similarly, in backward time) towards some hyperbolic fixed point,
will in fact converge at an exponential rate.

\begin{mythm}\label{thm:exponential_decay}
Let $f$ be of class $C^4$.
Suppose hypotheses \Hgroup are satisfied.
Let $Z$ be a hyperbolic rest point of \eqref{eq:TWE}.
Then there exist constants $C, \gamma, \epsilon > 0$ such that the following holds.
Let $U$ be a solution of \eqref{eq:TWE} on $J = [\tau_0,\infty)$ for which it holds that $\lim_{t\to\infty} \| U(t) - Z \|_{X^0} = 0$.
Then
\[
\| U - Z \|_{W^{1,2}((T,T+1);X^0,X^1)} \leq C e^{-\gamma (T - T_0)} \qquad \text{for all} \quad T \geq T_0 + 1,
\]
where $T_0 \geq \tau_0 + 1$ is chosen such that $\| U(T_0) - Z \|_{X^0} < \epsilon$.
\end{mythm}
A similar statement holds for solutions $U$ of \eqref{eq:TWE} on $J = (-\infty,-\tau_0]$ which converge towards a hyperbolic rest point
in backward time.
\begin{myproof}[of theorem \ref{thm:exponential_decay}]
First note that from lemma \ref{lemma:resolvent_decay}, 
the unique continuation theorem \ref{thm:unique_continuation} (using that $f$ is $C^4$),
and \cite{peterhof1997exponential}, the following follows:
there exists $\gamma > 0$ and $K>0$ such that any $W \in C^0(J;X^1) \cap C^1(\inter J;X^0,X^1)$ which satisfies
\[
\partial_t W + \d A(U) W = 0
\]
and $\sup_{t\in J} \|W(t)\|_{X^0} < \infty$,
satisfies the exponential decay estimate
\begin{equation}
  \label{eq:58}
  \| W(t) \|_{X^0} \leq K e^{-\gamma |t - \tau|} \| W(\tau) \|_{X^0} \qquad \text{for} \quad t \geq \tau \geq \tau_0.
\end{equation}
See appendix \ref{sec:appendix_exp_decay} for details.

We know that $U \in C_b^4(J;X^0,X^1,\dots,X^4)$, hence $W := \partial_t U$ satisfies $\sup_{t\in J} \|W(t)\|_{X^0} < \infty$
and also $\partial_t W + \d A(U) W = 0$ on $J$.
From \eqref{eq:58} it then follows that $W \in L^1(J;X^0)$, and for any $t, T_0 \in J$ with $t \geq T_0$ we have
\begin{equation}
  \label{eq:20}
    \| U(t) - Z \|_{X^0} = \bigg\| - \int_t^\infty W(s) \d s \bigg\|_{X^0} \leq \int_t^\infty \| W(s) \|_{X^0} \d s 
\leq \frac{K}{\gamma} e^{-\gamma |t-T_0|} \| W(T_0) \|_{X^0}.
\end{equation}

Now let
\begin{align*}
  &\calL_Z : W^{1,2}(\R;X^0,X^1) \to L^2(\R;X^0), \\
  &\calL_Z W = \partial_t W + \d A(Z) W.
\end{align*}
By hyperbolicity of $Z$ this operator has a continuous inverse;
this follows using lemma \ref{lemma:resolvent_decay} and the results from \cite{rabier2003isomorphism}.
Let $\eta_0 \in C^\infty_c(\R)$ be such that $\eta_0(t) = 1$ for $t \in [0,1]$ and $\eta_0(t) = 0$ for $t \not\in [-1,2]$.
Then set $\eta_T(t) := \eta_0(t-T)$.
Now fix $T_0, T \in J$ such that $[T-1,T+2] \subset [T_0,\infty)$ (i.e.\  $T \geq T_0 + 1$).
Note that
\begin{equation}
  \label{eq:56}
  \begin{split}
  \calL_Z( \eta_T(t) ( U(t) - Z ) ) &= \dot\eta_T(t) \big( U(t) - Z \big) + \eta_T(t) \big( \partial_t U(t) + \d A(Z)( U(t) - Z ) \big) \\
  &= \dot\eta_T(t) \big( U(t) - Z \big) + \eta_T(t) \big( \d A(Z)( U(t) - Z ) - ( A(U(t)) - A(Z) ) \big) \\
  &\quad + \eta_T(t) \big( \partial_t U(t) + A(U(t)) \big) - \eta_T(t) A(Z) \\
  &= \dot\eta_T(t) \big( U(t) - Z \big) + \eta_T(t) \big( \d A(Z)( U(t) - Z ) - ( A(U(t)) - A(Z) ) \big) \\
  &= \dot\eta_T(t) \big( U(t) - Z \big) + \eta_T(t)
  \begin{pmatrix}
    0 \\
    f_u(x,z) ( u(t) - z )  - f(x,u) + f(x,z) 
  \end{pmatrix} \\
  &\quad + \eta_T(t)
  \begin{pmatrix}
    0 \\
    g(x,U(t))
  \end{pmatrix},
  \end{split}
\end{equation}
where $U = (u,\partial_t u)$ and $Z = (z,0)$.

Since the map $f : H^1(\Omega) \to L^2(\Omega)$ is differentiable near $z$ we have
\[
\bigg\|
 \begin{pmatrix}
    0 \\
    f_u(x,z) ( u(t) - z )  - f(x,u) + f(x,z) 
  \end{pmatrix}
\bigg\|_{X^0}
= o( \| U(t) - Z \|_{X^0} ) \qquad \text{as} \quad t \to \infty.
\]
By hypotheses \ref{H:growth_bound_g}, \ref{H:cone_condition}, and a mean value estimate, 
we find that for all $t \in J$ one has
\begin{equation*}
    \bigg\| 
\begin{pmatrix}
    0 \\
    g(x,U(t))
\end{pmatrix}
\bigg\|_{X^0}
= \bigg\| \int_0^1 \ddxe{s} g(x, Z + s(U(t)-Z)) \d s \bigg\|_{L^2(\Omega)}
\leq C_{1,g} \sqrt{ \Vol(\Omega) } \| U(t) - Z \|_{X^0}.
\end{equation*}
Combining these observations with estimate \eqref{eq:56}, we find that
\[
\| \calL_Z( \eta_T(t) ( U(t) - Z ) ) \|_{X^0} = O(\|U(t)-Z\|_{X^0}) \qquad \text{for} \quad t \in J, \text{ uniform in } T.
\]
Since $\supp( \eta_T ) \subset [T_0,\infty)$ we can now apply estimate \eqref{eq:20} to obtain
\[
\| \calL_Z( \eta_T(t) ( U(t) - Z ) ) \|_{X^0} \leq C e^{-\gamma |t-T_0|} \| W(T_0) \|_{X^0},
\]
for some constant $C$ which is independent of $T$.
By invertibility of $\calL_Z : W^{1,2}(\R;X^0,X^1) \to L^2(\R;X^0)$ we then find
\begin{align*}
  \| U - Z \|_{W^{1,2}((T,T+1);X^0,X^1)} &\leq \| \eta_T \cdot ( U - Z ) \|_{W^{1,2}(\R;X^0,X^1)}  \\
  &\leq \| \calL_Z^{-1} \|\| \calL_Z( \eta_T \cdot ( U - Z ) ) \|_{L^2(\R;X^0)} \\
  &\leq \frac{ C \| \calL_Z^{-1} \| }{ \gamma } e^{-\gamma |T-T_0|} \| W(T_0) \|_{X^0}.
\end{align*}

The above estimate already proves the exponential decay of the single solution $U$ towards $Z$.
Next, we will argue why the constant $\| W(T_0) \|_{X^0} = \| \partial_t U(T_0) \|_{X^0}$ can be bounded independently of $U$,
as long as $T_0$ is chosen big enough such that $U(T_0)$ lies in some given $X^0$-neighbourhood of $Z$.
Let $a := \calE(Z)$, and choose any $b > a$.
From section \ref{sec:compactness} it follows that $\rst{ \calA_a^b(J,f,g,c) }{ [\tau_0+1,\infty) }$ is bounded in
the norm topology on $C_b^1([\tau_0+1,\infty);X^0,X^1)$.
Denote this upper bound by $M$.
By continuity of $\calE : X^0 \to \R$, there exists some $\epsilon > 0$ such that if $T_0 \geq \tau_0+1$ 
is such that $\| U(T_0) - Z \|_{X^0} < \epsilon$, then $a \leq \calE(U(T_0)) \leq b$.
If we define $\widetilde U(t) := U(t-\tau_0-1+T_0)$, then $\widetilde U \in \rst{ \calA_a^b(J,f,g,c) }{ [\tau_0+1,\infty) }$, hence
\[
\| \partial_t U(T_0) \|_{X^0} = \| \partial_t \widetilde U(\tau_0+1) \|_{X^0} \leq M.
\]
\end{myproof}

\section{Moduli spaces}
\label{sec:moduli}

At this point we are ready to define the moduli spaces of heteroclinic orbits.
Given two hyperbolic rest points $Z_- , Z_+ \in \calS_{\text{hyp}}$, we define
\[
\calM(Z_-,Z_+) = \set{ U \in W^{1,2}_{\text{loc}}(\R;X^0,X^1) }{ U \text{ solves \eqref{eq:TWE}, } \lim_{t\to\pm\infty} \| U(t) - Z_\pm \|_{X^0} = 0 }.
\]
We will now discuss two useful modes of convergence on this space, and the interplay between these two.
For conciseness we will do this only for the autonomous case.
We refer to remark \ref{remark:nonautonomous_moduli} at the end of this section for details on how to adapt the arguments to the nonautonomous setting.

\subsection{The manifold structure}

Let $Z \in C^1(\R;X^0,X^1)$ be such that $Z(t) = Z_-$ for $t \leq -1$ and $Z(t) = Z_+$ for $t \geq 1$.
We then define the affine space of paths between $Z_-$ and $Z_+$
\[
\calP(Z_-,Z_+) := Z + W^{1,2}(\R;X^0,X^1).
\]
By the exponential decay of solutions towards hyperbolic rest points we have $\calM(Z_-,Z_+) \subset \calP(Z_-,Z_+)$.
Therefore, if we introduce the nonlinear operator
\begin{align*}
  &\Phi_{Z_-,Z_+} : \calP(Z_-,Z_+) \to L^2(\R;X^0), \\
  &\Phi_{Z_-,Z_+}(U) = \partial_t U + A(U),
\end{align*}
we find that
\[
\calM(Z_-,Z_+) = \Phi_{Z_-,Z_+}^{-1}(0).
\]
The linearization of $\Phi_{Z_-,Z_+}$ around any $U \in \calP(Z_-,Z_+)$ is Fredholm with $\ind( \d\Phi_{Z_-,Z_+}(U) ) = \mu(Z_-) - \mu(Z_+)$.
This index is independent of the chosen $U \in \calP(Z_-,Z_+)$, 
thus $\Phi_{Z_-,Z_+}$ is a Fredholm map with $\ind( \Phi_{Z_-,Z_+} ) = \mu(Z_-) - \mu(Z_+)$.

Recall from section \ref{subsec:transversality_definition} the definition of transversality up to order $k$.
We now assume that the transversality condition up to order $\mu(Z_-)-\mu(Z_+)$ is satisfied.
If the nonlinearity $f$ is of class $C^m$, 
then the map $\Phi_{Z_-,Z_+}$ is of class $C^m$ in a neighbourhood of $\calM(Z_-,Z_+)$ (recall theorem \ref{thm:operator_regularity}), 
hence the implicit function theorem implies that $\calM(Z_-,Z_+)$ with the topology inherited from $\calP(Z_-,Z_+)$ 
is again a $C^m$ manifold of finite dimension $\mu(Z_-) - \mu(Z_+)$.

Since we study the autonomous case, time translations $s \mapsto U(\cdot + s)$ induce an $\R$-action on $\calM(Z_-,Z_+)$.
We will denote the equivalence class of $U \in \calM(Z_-,Z_+)$ with respect to this $\R$-action by $[ U ]$, and write
\[
\widehat \calM(Z_-,Z_+) := \calM(Z_-,Z_+) / \R
\]
for the quotient space.
\begin{mylemma}
  Assume the transversality condition up to order $\mu(Z_-)-\mu(Z_+)$ is satisfied.
  Then the space $\widehat \calM(Z_-,Z_+)$ is a $C^m$ manifold of dimension $\mu(Z_-) - \mu(Z_+) - 1$.
\end{mylemma}
\begin{myproof}
  We verify that the $\R$-action of time translations on $\calM(Z_-,Z_+)$ is $C^m$, free, and proper.
  The lemma then follows from the quotient manifold theorem \cite{Lee2003}.
  The $C^m$ smoothness of the action, where we consider $\calM(Z_-,Z_+)$ to be endowed with the topology inherited from $\calP(Z_-,Z_+)$,
  follows from the regularity results stated in section \ref{sec:compactness}.
  Since the energy functional $\calE$ is strictly decreasing with time translations, see also remark \ref{remark:strict_Lyapunov}, the action is free.
  Finally, properness follows from the compactness estimates in section \ref{sec:compactness}.
\end{myproof}


\subsection{Geometric convergence}

\begin{figure}
  \centering
  
  \begin{subfigure}{0.48\textwidth}
   \centering
  \makeatletter
\tikzset{nomorepostaction/.code={\let\tikz@postactions\pgfutil@empty}}
\makeatother

\begin{tikzpicture}[scale=1]
\coordinate (A1) at (0,-2);
\coordinate (A2) at (1.4,-1);
\coordinate (A3) at (1.8,0);
\coordinate (A4) at (1.4,1);
\coordinate (A5) at (0,2);

\coordinate (T1) at (2.2,-1.5);
\coordinate (T2) at (2.2,1.5);

\begin{scope}[every path/.style={color=black!60!white,
    postaction={ nomorepostaction, decorate, decoration={ markings,mark=at position 0.5 with \arrow{>} } }
  }]
\draw[thick] (A2)--(A1);
\draw[thick] (A3)--(A2);
\draw[thick] (A4)--(A3);
\draw[thick] (A5)--(A4);
\end{scope}

\begin{scope}[every path/.style={
    postaction={ nomorepostaction, decorate, decoration={ markings,mark=at position 0.5 with \arrow{>} } }
  }]
\draw[thick] (A5) to[in=40,out=-40] (A1);
\draw[thick] (A5) to[in=70,out=-70] (A1);
\end{scope}

\begin{scope}
\draw[thick,dashed] (0.7,1.4) to[in=60,out=-60] (0.7,-1.4);
\draw[thick,dashed] (1,1.2) to[in=60,out=-60] (1,-1.2);
\draw[thick,dashed] (1.3,0.8) to[in=60,out=-60] (1.3,-0.8);
\end{scope}

\draw[fill=black] (A1) circle (0.1em); 
\draw[fill=black!60!white] (A2) circle (0.1em); 
\draw[fill=black!60!white] (A3) circle (0.1em); 
\draw[fill=black!60!white] (A4) circle (0.1em); 
\draw[fill=black] (A5) circle (0.1em);

\draw [->] (T1) -- (T2) node[below right] {$\mathcal{E}$}; 

\end{tikzpicture}
  \caption{Convergence towards a broken trajectory.}
  \label{fig:geometric_convergence}
  \end{subfigure}
  \hspace*{\fill}
  \begin{subfigure}{0.48\textwidth}
  \centering

\makeatletter
\tikzset{nomorepostaction/.code={\let\tikz@postactions\pgfutil@empty}}
\makeatother

\begin{tikzpicture}[scale=1.5]
\coordinate (A1) at (0,1);
\coordinate (A2) at (-1,0);
\coordinate (A3) at (1,0);
\coordinate (A4) at (0,-1);

\fill[fill=black!10!white] (A1)--(A2)--(A4)--(A3);

\draw (0,0) node {$2$};

\begin{scope}[every path/.style={postaction={ nomorepostaction, decorate, decoration={ markings,mark=at position 0.5 with \arrow{>} } } }]
\draw[thick] (A1)-- node[above left] {$1$} (A2);
\draw[thick] (A1)-- node[above right] {$1$} (A3);
\draw[thick] (A2)-- node[below left] {$1$} (A4);
\draw[thick] (A3)-- node[below right] {$1$} (A4);
\end{scope}

\draw[fill=black] (A1) circle (0.1em) node[above] {$Z_-$};
\draw[fill=black] (A2) circle (0.1em) node[left] {$Y$};
\draw[fill=black] (A3) circle (0.1em) node[right] {$Y'$};
\draw[fill=black] (A4) circle (0.1em) node[below] {$Z_+$};
\end{tikzpicture}
  \caption{Noncompact connected component of dimension $2$ of $\calM(Z_-,Z_+)$.}
  \label{fig:moduli}
  \end{subfigure}

  \caption{Analysis of the moduli spaces $\calM(Z_-,Z_+)$ in the autonomous setting.}
\end{figure}
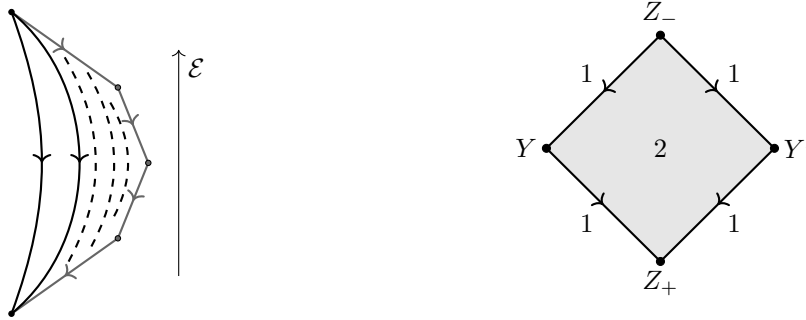

In general $\calM(Z_-,Z_+)$ will not be compact in the topology of $\calP(Z_-,Z_+)$;
just consider a sequence of time translations of a nonconstant solution, this cannot have a convergent subsequence in this topology.
In order to better understand exactly how $\calM(Z_-,Z_+)$ fails to be compact, we introduce the notion of geometric convergence.
\begin{mydef}
Let $Z_0,\dots,Z_{k+1}$ be hyperbolic stationary solutions of \eqref{eq:TWE}.
Pick any
\[
( [U_0] , [U_1], \dots , [U_{k-1}] , [U_k] ) \in \widehat\calM(Z_0,Z_1) \times \widehat\calM(Z_1,Z_2) \times \cdots \times 
\widehat\calM(Z_{k-1},Z_k) \times \widehat\calM(Z_k , Z_{k+1}).
\]
We will call such a $k$-tuple $([V_0],\dots,[V_k])$ a $k$-fold broken orbit.

Given $([U_n])_n \subset \widehat\calM(Z_0,Z_{k+1})$, we will say that $([U_n])_n$ converges geometrically to the $k$-fold broken orbit $([V_0],\dots,[V_k])$ if the following holds.
For each $j \in \{0,\dots,k\}$ there exists a sequence $(s_{j,n})_n \subset \R$ so that
\[
U_n(\cdot + s_{j,n}) \to V_j \qquad \text{in} \quad W^{1,2}_{\text{loc}}(\R;X^0,X^1), \qquad \text{as} \quad n \to \infty,
\]
and this $k$-fold broken orbit is maximal in the sense that for each sequence $(s_n')_n \subset \R$ and $V'$ such that 
$U_n(\cdot + s_n') \to V'$ in $W^{1,2}_{\text{loc}}(\R;X^0,X^1)$ as $n\to\infty$ it holds that $[V'] \in \{[V_0],\dots,[V_n]\}$.
In this case we will also write
\[
[U_n] \leadsto ([V_0],\dots,[V_k]) \qquad \text{as} \quad n \to \infty.
\]
See also figure \ref{fig:geometric_convergence}.
\end{mydef}

To arrive at the main result regarding geometric convergence, 
we will first discuss a property which was already hinted at in the introduction.
First, given a solution $U$ of \eqref{eq:TWE}, let $\alpha(U)$ denote the set of accumulation points in the topology of $W^{1,2}_{\text{loc}}(\R;X^0,X^1)$
of the sequence $( U(\cdot + \tau) )_\tau$ where $\tau \to -\infty$.
Similarly, let $\omega(U)$ denote the set of accumulation points of the sequence $( U(\cdot + \tau) )_\tau$ where $\tau \to \infty$.
In other words, these are the $\alpha$- and $\omega$-limit sets of the shift dynamics on the set of solutions to \eqref{eq:TWE}.
\begin{mylemma}[Gradient-like behaviour]
\label{lemma:gradient_like}
  Consider a $t$-independent triple $(f,g,c)$ satisfying \Hgroup.
  Let $U$ be a bounded solution of \eqref{eq:TWE}.
  Then the limit sets $\alpha(U)$ and $\omega(U)$ are both nonempty and connected, and $\alpha(U) , \omega(U) \subset \calS$.
  Moreover, if $\alpha(U) \cap \omega(U) \neq \emptyset$, then $U$ is a stationary solution of \eqref{eq:TWE}.
  If in addition we assume that $f$ is regular, then $\alpha(U)$ and $\omega(U)$ each consist of a single point,
  hence any bounded solution of \eqref{eq:TWE} is either a rest point or a heteroclinic orbit.
\end{mylemma}
\begin{myproof}
Let $U \in W^{1,2}_{\text{loc}}(\R;X^0,X^1)$ be a solution of \eqref{eq:TWE} with $\sup_{t\in\R} \| U(t) \|_{X^0} < \infty$.
Then there exist $a , b \in \R$ such that $a \leq \calE(U(t)) \leq b$ for all $t \in \R$.
Endow $\calA_a^b$ with the topology inherited from $W^{1,2}_{\text{loc}}(\R;X^0,X^1)$.
Then note that time translation $s \mapsto U(\cdot + s)$ defines a continuous dynamical system on $\calA_a^b$.
The compactness results from section \ref{sec:compactness} imply that 
the $\alpha$- and $\omega$-limit sets of any $U \in \calA_a^b$ are nonempty.
The dynamical system also posesses a Lyapunov function given by $\calL(U) = \calE(U(0))$, with $\calE$ as defined in section \ref{subsec:Lyapunov}.
In light of theorem \ref{thm:unique_continuation} this Lyapunov function is stricly decreasing along nonstationary trajectories.
Since the Lyapunov function must be constant on the limit sets, this implies that $\alpha(U) , \omega(U) \subset \calS_a^b$.
Futhermore, this implies that $\alpha(U) \cap \omega(U) = \emptyset$ when $U$ is a nonstationary solution of \eqref{eq:TWE}.

If we in addition assume that $f$ is regular, the hyperbolicity implies that $\calS_a^b$ is totally disconnected.
On the other hand, since $\alpha(U)$ and $\omega(U)$ are the limit sets of a continuous dynamical system, they are both connected.
Hence each of these limit sets must consist of a single point.
\end{myproof}

With this settled, a standard argument shows that $\widehat \calM(Z_-,Z_+)$ is compact up to broken orbits.
Roughly speaking, any sequence $(U_n)_n \subset \calM(Z_-,Z_+)$ must have some convergent subsequence by the results from section \ref{sec:compactness},
and by lemma \ref{lemma:gradient_like} the limit point must then again belong to some moduli space $\calM(Z_i,Z_j)$.
If $Z_i \neq Z_-$ or $Z_j \neq Z_+$, using the energy functional $\calE$ defined in section \ref{subsec:Lyapunov} we can then find a sequence $(t_n)_n \subset \R$ such that
the time-translated sequence $(U_n(t_n + \cdot))_n$ converges over a subsequence towards a limit point in yet another moduli space $\calM(Z_k,Z_l)$.
This iterative procedure must terminate after a finite number of steps, 
since hyperbolicity implies there are only a finite number of rest points with bounded energy.
See figure \ref{fig:geometric_convergence}, and also \cite{schwarz1993morse} for additional details.
Note that in this reference only gradient flows are considered, but the particular argument only relies on the existence of a strict Lyapunov function.
Summarizing, we have the following result.
\begin{mythm}
\label{thm:compactness_broken_orbit}
Consider \eqref{eq:TWE} with $t$-independent $(f,g,c)$.
  Suppose all rest points of \eqref{eq:TWE} are hyperbolic.
  Then the space $\widehat \calM(Z_-,Z_+)$ is compact up to broken orbits, i.e., for each $([U_n])_n \subset \widehat \calM(Z_-,Z_+)$ there exists
a $k \in \N_0$, intermediate points $Z_1,\dots,Z_k \in \calS$, a broken orbit
\[
( [V_0] , [V_1], \dots , [V_{k-1}] , [V_k] ) \in \widehat\calM(Z_-,Z_1) \times \widehat\calM(Z_1,Z_2) \times \cdots \times \widehat\calM(Z_{k-1},Z_k) \times \widehat\calM(Z_k , Z_+),
\]
and a subsequence $([U_{n_j}])_j$ such that
\[
[U_{n_j}] \leadsto ([V_0],\dots,[V_k]) \qquad \text{as} \quad j \to \infty.
\]
See also figure \ref{fig:geometric_convergence}.
\end{mythm}

\subsection{Relating the modes of convergence}

The following theorem highlights an important link between geometric convergence and the topology inherited from $\calP(Z_-,Z_+)$.
An important consequence will be that for generic choices of $f$ and $g$ one can count index $1$ orbits.
\begin{mythm}\label{thm:unparameterized_trajectories_compactness}
Consider \eqref{eq:TWE} with $t$-independent $(f,g,c)$.
Let $Z_-$ and $Z_+$ be hyperbolic rest points of \eqref{eq:TWE}, and suppose $\widehat \calM(Z_-,Z_+)$ is compact up to $0$-fold broken orbits.
Then $\widehat \calM(Z_-,Z_+)$ is sequentially compact in the quotient topology inherited from $\calP(Z_-,Z_+)$.
In particular, if $\mu(Z_-) - \mu(Z_+) = 1$, then $\widehat \calM(Z_-,Z_+)$ is sequentially compact.
\end{mythm}
\begin{myproof}
Without loss of generality we can assume $\calE(Z_+) < \calE(Z_-)$; 
because otherwise either $\calM(Z_-,Z_+) = \emptyset$ (if $Z_- \neq Z_+$) or $\calM(Z_-,Z_+) = \{ Z_- \}$ (if $Z_- = Z_+$).
Select any $([U_n])_n \subset \widehat \calM(Z_-,Z_+)$.
For each $n$ fix a representative $U_n$ of $[U_n]$.
By assumption we may find a subsequence $(U_n')_n$, a sequence $(t_n)_n \subset \R$, 
and a limit point $V_0 \in \calM(z_-,z_+)$ such that $U_n'(\cdot + t_n) \to V_0$ in $W^{1,2}_{\text{loc}}(\R;X^0,X^1)$ as $n\to\infty$.
After replacing $U_n'$ by $U_n'(\cdot + t_n)$, i.e.\  choosing a different representative for $[U_n']$, 
we may as well assume that $U_n' \to V_0$ in $W^{1,2}_{\text{loc}}(\R;X^0,X^1)$ as $n\to\infty$.
To ease on notation, we shall henceforth denote the subsequence $(U_n')_n$ by just $(U_n)_n$.

We claim that $U_n(t)$ converges uniformly in $n$ towards $Z_\pm$ as $t \to \pm\infty$, i.e.\ 
\begin{equation}
  \label{eq:10}
  \forall \epsilon > 0\; \exists T_0 > 0\; \forall n \in \N\; \forall t \geq T_0 : \qquad \| U_n(t) - Z_\pm \|_{X^0} \leq \epsilon.
\end{equation}
Suppose for the moment that this is true.
Choosing $\epsilon$ as small as needed, we may now apply the exponential decay theorem \ref{thm:exponential_decay}
to find $\delta_0 > 0$, $C > 0$, and $T_0 > 0$ such that
\begin{align*}
  \| U_n - V_0 \|_{W^{1,2}((T,T+1);X^0,X^1)} &\leq \| U_n - Z_\pm \|_{W^{1,2}((T,T+1);X^0,X^1)} \\
  &\quad + \| V_0 - Z_\pm \|_{W^{1,2}((T,T+1);X^0,X^1)} \\
  &\leq C e^{-\delta_0 T}
\end{align*}
for all $T \geq T_0$.
Consequently,
\[
\| U_n - V_0 \|_{W^{1,2}(\R \setmin [-T,T] ; X^0,X^1)} \to 0 \qquad \text{as} \quad T \to \infty, \text{ uniformly in } n \in \N.
\]
Since we also have that $U_n \to V_0$ in $W^{1,2}_{\text{loc}}(\R;X^0,X^1)$ as $n \to \infty$, we in particular find that for each $T > 0$ it holds that
\[
\| U_n - V_0 \|_{W^{1,2}((-T,T) ; X^0,X^1)} \to 0 \qquad \text{as} \quad n \to \infty.
\]
Together these estimates imply that $\|U_n - V_0\|_{W^{1,2}(\R;X^0,X^1)} \to 0$ as $n\to\infty$, so that $U_n \to V_0$ in $\calP(Z_-,Z_+)$ as $n\to\infty$.
Hence $\widehat \calM(Z_-,Z_+)$ is compact in the quotient topology inherited from $\calP(Z_-,Z_+)$.

It remains to prove the claim \eqref{eq:10}.
We will only discuss the uniform convergence towards $Z_+ = (z_+,0)$ as $t\to\infty$; the case for $t \to -\infty$ is obtained analogously.
We shall first show that $\| \partial_t u_n(t) \|_{L^2(\Omega)} \to 0$ uniformly in $n$ as $t \to \infty$, where $U_n = (u_n,\partial_t u_n)$.
Assume on the contrary that we can find a sequence $(t_k)_k \subset \R$ with $t_k \to \infty$, $(n_k)_k \subset \N$ with $n_k \to \infty$, and $\epsilon > 0$ such that
\begin{equation}
  \label{eq:39}
  \| \partial_t u_{n_k}(t_k) \|_{L^2(\Omega)} > \epsilon \qquad \text{for all} \quad k \in \N.
\end{equation}
Recall from section \ref{sec:compactness} that $\calM(Z_-,Z_+)$ is uniformly bounded in $C_b^1(\R;X^0,X^1)$. 
From this we obtain the following equicontinuity condition:
there exists an $M \in \R$ such that
\begin{align*}
  \bigg| \| \partial_t w(s) \|_{L^2(\Omega)}^2 - \|\partial_t w(t) \|_{L^2(\Omega)}^2 \bigg| &\leq
\int_t^s \big| \partial_\tau \| \partial_t w(\tau) \|_{L^2(\Omega)}^2 \big| \d\tau \\
&= 2 \int_s^t | \langle \partial_t^2 w(\tau) , \partial_t w(\tau) \rangle_{L^2(\Omega)} |\d\tau \\
&\leq 2 \int_s^t \| \partial_t^2 w(\tau) \|_{L^2(\Omega)} \|\partial_t w(\tau) \|_{L^2(\Omega)} \d\tau \\
&\leq M |s-t|
\end{align*}
for all $s, t \in \R$ and $W = (w,\partial_t w) \in \calM(Z_-,Z_+)$.
Combining this with \eqref{eq:39}, we may find $\delta > 0$ such that
\[
\| \partial_t u_{n_k}(s) \|_{L^2(\Omega)}^2 \geq \epsilon^2 / 2 \qquad \text{whenever} \quad |t_{n_k} - s| \leq \delta. 
\]
In turn, this gives
\begin{align*}
  \calE(U_{n_k}(t_k)) - \calE(u_{n_k}(t_k + \delta)) &= - \int_{t_k}^{t_k+\delta} \partial_\tau \calE(U_{n_k}(\tau)) \d\tau \\
  &\geq \frac c 2 \int_{t_k}^{t_k+\delta} \| \partial_t u_{n_k}(\tau) \|_{L^2(\Omega)}^2 \d \tau \\
  &\geq c \delta \epsilon^2 / 4 =: \epsilon_0.
\end{align*}
Since $t \mapsto \calE(U_{n_k}(t))$ is a decreasing function, this finally implies that for all $k$
\[
\calE(U_{n_k}(t_k)) - \calE(Z_+) \geq \epsilon_0.
\]
But we can find $s_0 \in\R$ such that $\calE(V_0(s_0)) - \calE(Z_+) = \epsilon_0 / 2$, and since $t_k \to \infty$ we then find that for $k$ large enough
\[
\calE(U_{n_k}(s_0)) - \calE(Z_+) \geq \calE(U_{n_k}(t_k)) - \calE(Z_+) \geq \epsilon_0.
\]
However, $U_{n_k} \to V_0$ in $W^{1,2}_{\text{loc}}(\R;X^0,X^1)$, so we in particular have $\calE(U_{n_k}(s_0)) \to \calE(V_0(s_0))$.
Hence the left hand side of the last inequality tends to $\epsilon_0 / 2$ as $k \to \infty$ over a subsequence, which is impossible since $\epsilon_0 > 0$.
This contradiction shows that $\| \partial_t u_n(t) \|_{L^2(\Omega)} \to 0$ uniformly in $n$ as $t \to \infty$.

Next, we will show that from this it also follows that $\| u_n(t) - z_+ \|_{H^1(\Omega)} \to 0$ uniformly in $n$ as $t \to \infty$,
thus completing the proof of \eqref{eq:10}.
Suppose on the contrary that we can find a sequence $(t_k)_k \subset \R$ with $t_k \to \infty$, $(n_k)_k \subset \N$ with $n_k \to \infty$, and $\epsilon > 0$ such that
\[
    \| u_{n_k}(t_k,\cdot) - z_+ \|_{H^1(\Omega)} > \epsilon \qquad \text{for all} \quad k \in \N.
\]
From section \ref{sec:compactness} we know that, after selecting a further subsequence,
 $U_{n_k}(t_k + \cdot) \to W = (w,\partial w)$ in $W^{1,2}_{\text{loc}}(\R;X^0,X^1)$ as $k \to \infty$.
The assumed inequality ensures that $w \neq z_+$.
Yet since $\| \partial_t u_n(t) \|_{L^2(\Omega)} \to 0$ uniformly in $n$ as $t\to\infty$, it holds that $\partial_t w = 0$.
Since $\calM(Z_-,Z_+)$ is assumed compact up to $0$-fold broken orbits, we thus must have $W = Z_-$.
Hence for any $\delta > 0$ and $s_0 > 0$ we can choose $k$ large enough so that
\[
\calE(Z_+) < \calE(Z_-) \leq \calE(U_{n_k}(t_k)) + \frac \delta 2 \leq \calE(U_{n_k}(s_0)) + \frac \delta 2 \leq \calE(V_0(s_0)) + \delta.
\]
However, $\calE(V_0(s_0)) \to \calE(Z_+)$ as $s_0 \to \infty$.
So we find that for any $\delta > 0$ we have $\calE(Z_+) < \calE(Z_-) \leq \calE(Z_+) + \delta$.
Hence $\calE(Z_-) = \calE(Z_+)$, in contradiction with the assumption that $\calE(Z_+) < \calE(Z_-)$.
\end{myproof}

\subsection{The glueing map}

The following ``glueing theorem'' allows us to understand the structure of the boundary of the two-dimensional moduli spaces.

\begin{mythm} \label{thm:glueing}
Suppose $Z_0$, $Z_1$, $Z_2$ are hyperbolic stationary solutions of \eqref{eq:TWE}, where $\mu(Z_0) = \mu(Z_1)+1 = \mu(Z_2)+2$.
Assume the transversality condition up to order $2$ is satisfied.
Let $(U,V) \in \calM(Z_0,Z_1) \times \calM(Z_1,Z_2)$.
Then there exists an immersion
\begin{align*}
  & \# : [T_0,\infty) \to \calM(Z_0,Z_2), \\
  & T \mapsto U \#_T V,
\end{align*}
such that $[U \#_T V] \leadsto ([U] , [V])$ as $T \to \infty$.
Moreover, any sequence in $\widehat\calM(Z_0,Z_2)$ which converges geometrically towards $([U],[V])$ eventually lies within the range of $[U \# V]$.
\end{mythm}

The ideas in this construction are fairly standard, 
see e.g.\  \cite{floer1989symplectic, schwarz1995cohomology, schwarz1993morse, audin2014morse}.
We will only give a sketch here.
First define a pre-glueing map
\[
\big( U \#^0_T V \big)(t) := \bigg( 1 - \eta\bigg( \frac t T \bigg) \bigg) U(t+2T) + \eta\bigg( \frac t T \bigg) V(t - 2T),
\]
where $\eta \in C^\infty(\R)$ is such that $0 \leq \eta \leq 1$, $\eta(t) = 0$ for $t \leq -1$, and $\eta(t) = 1$ for $t \geq 1$.
Note that $U \#^0_T V  \in \calP(Z_0,Z_2)$,
and $U \#^0_T V$ converges geometrically towards $([U] , [V])$ as $T \to \infty$.
However, the pre-glueing is in general not a solution of \eqref{eq:TWE},
but a brief computation yields
\[
  \big\| \Phi_{Z_0,Z_2}( U \#_T^0 V ) \big\|_{L^2(\R;X^0)} \to 0 \qquad \text{as} \quad T \to \infty,
\]
which suggest there must be a true solution nearby the pre-glueing.
The aim is to find this true solution using a contraction mapping argument.

Define
\begin{align*}
  & F_{T} : W^{1,2}(\R;X^0,X^1) \to L^2(\R;X^0), \\
  & F_{T}(\gamma) = \Phi_{Z_0,Z_2}(U \#^0_T V + \gamma).
\end{align*}
By hyperbolicity and transversality, the maps $\d \Phi_{Z_0,Z_1}(U)$ and $\d \Phi_{Z_1,Z_2}(V)$ are surjective Fredholm operators,
hence they have bounded right inverses $M_{01}$ and $M_{12}$ respectively.
We then patch these operators together to obtain an ``approximate right inverse'' for $\d F_{T}(0)$:
\[
M_{T} := \zeta_T^- \tau_{2T} M_{01} \tau_{-2T} \zeta_T^- + \zeta_T^+ \tau_{-2T} M_{12} \tau_{2T} \zeta_T^+.
\]
Here $\tau_a$ denotes the operator induced by time translation $t \mapsto t + a$, and $\zeta_T^\pm \in C^\infty(\R)$ is such that 
$\zeta_T^-(t)^2 + \zeta_T^+(t)^2 = 1$, $\zeta_T^+(t) = 0$ for $t \leq -T$, $\zeta_T^-(t) = \zeta_T^+(-t)$, $\zeta_T^\pm(t) = \zeta_1^\pm(t/T)$.
This is an approximate right inverse in the sense that $\d F_{T}(0) \circ M_{T} \to I$ in $\calL(L^2(\R;X^0))$ as $T \to \infty$.
In turn, this implies the existence of a true right inverse $G_{T}$ of $\d F_{T}(0)$;
in addition the operator norm of $G_{T}$ can be bounded independent of $T$.
This allows us to define a Newton-like operator
\begin{align*}
  & N_{T} : W^{1,2}(\R;X^0,X^1) \to \img(G_{T}), \\
  & N_{T} = G_{T} \circ \big( \d F_{T}(0) - F_{T} \big).
\end{align*}
Since $F_{T}(0) \to 0$ as $T \to \infty$ and the norm of $G_T$ can be bounded independent of $T$, 
this operator turns out to be a contraction for $T$ large enough.
Consequently, there exists an $\epsilon > 0$ such that, for large $T$,
there exists a unique $\gamma_{T} \in B_\epsilon(0) \cap \img(G_{T})$ such that $F_{T}(\gamma_{T}) = 0$.
Furthermore, one can show that $\| \gamma_{T} \|_{W^{1,2}(\R;X^0,X^1)} \to 0$ as $T \to \infty$.
Hence, if we set
\begin{equation}
  \label{eq:40}
  U \#_T V := U \#^0_T V + \gamma_{T},
\end{equation}
then $U \#_T V \in \calM(Z_0,Z_2)$ and $[ U \#_T V ] \leadsto ( [U] , [V] )$ as $T \to \infty$.

\subsection{The geometric picture}

Now assume that both $f$ and $g$ are regular.
Pick any $Z_- , Z_+ \in \calS$ with $\mu(Z_-) - \mu(Z_+) \leq 2$.
Suppose $\calM(Z_-,Z_+) \neq \emptyset$, and pick $U \in \calM(Z_-,Z_+)$.
Then
\[
\mu(Z_-) - \mu(Z_+) = \dim \ker( \d \Phi_{Z_-,Z_+}(U) ) \geq 0.
\]
Hence transversality implies that the Morse index can never increase along orbits.

We now combine the various results from this chapter to get a detailed picture of geometric properties of the trajectory spaces.
We present this in the following list.
\begin{description}[font=\normalfont]
\item[$\mu(Z_-) = \mu(Z_+)$.]
Assume that $\calM(Z_-,Z_+) \neq \emptyset$.
In this case $\calM(Z_-,Z_+)$ is a $0$-dimensional manifold, i.e.\  a discrete set.
If $U \in \calM(Z_-,Z_+)$, then $s \mapsto U(s+\cdot)$ defines a continuous curve in $\calM(Z_-,Z_+)$.
This curve has to be constant since $\calM(Z_-,Z_+)$ is discrete.
Hence $U$ is $t$-independent, and since $U(t) \to Z_\pm$ as $t \to\pm\infty$, we find that $Z_- = Z_+$.
We conclude that \emph{the space of index $0$ trajectories $\calM(Z_-,Z_+)$ is a finite set.}

\item[$\mu(Z_-)=\mu(Z_+)+1$.]
Let $([U_n])_n \subset \widehat\calM(Z_-,Z_+)$ and suppose that $[U_n] \leadsto ([V_0],\dots,[V_k])$, a $k$-fold broken orbit.
Since the index can never increase along orbits, we find that all except one of the $V_j$'s is of index $0$.
But as we just saw, index $0$ orbits are stationary solutions, hence all but one of the $V_j$'s equal either $Z_-$ or $Z_+$.
So $([U_n])_n$ converges geometrically to a $0$-fold broken orbit, i.e.\  an element of $\widehat\calM(Z_-,Z_+)$ itself.
From this and theorem \ref{thm:compactness_broken_orbit} we deduce that $\widehat\calM(Z_-,Z_+)$ is compact up to $0$-fold broken orbits.
We then conclude from theorem \ref{thm:unparameterized_trajectories_compactness} that 
$\widehat \calM(Z_-,Z_+)$ is sequentially compact in the quotient topology.
We also know that $\widehat \calM(Z_-,Z_+)$ is a $0$-dimensional manifold with this quotient topology.
Consequently $\widehat \calM(Z_-,Z_+)$ is a finite set.
This means that \emph{modulo time shifts, $\calM(Z_-,Z_+)$ consists of finitely many trajectories.}

\item[$\mu(Z_-) = \mu(Z_+)+2$.]
Let $O$ be a connected component of $\calM(Z_-,Z_+)$.
Arguing as above, we find that $O$ is either compact up to $0$-fold broken orbits, or compact up to $1$-fold broken orbits.

We can also study $\widehat O := O / \R$, the $1$-dimensional manifold obtained by dividing out the time shifts.
Since $\calM(Z_-,Z_+)$ is obtained as a regular level $0$ set via the implicit function theorem, it is a manifold without boundary.
Therefore also $O$ and hence $\widehat O$ are manifolds without boundary.
It follows from the classification of $1$-dimensional spaces that $\widehat O$ is homeomorphic to either $S^1$ or $(0,1)$.
The former corresponds to the case where $\widehat O$ is compact up to $0$-fold broken orbits.
In the latter case, we obtain a $1$-parameter family $([U_s])_{s\in(0,1)} \subset \widehat O$
such that $[U_s] \leadsto ([V_0^-],[V_1^-])$ as $s \downarrow 0$ and $[U_s] \leadsto ([V_0^+],[V_1^+])$ as $s \uparrow 1$.
These broken orbits $([V_0^-],[V_1^-])$ and $([V_0^+],[V_1^+])$ are distinct,
otherwise we would have two sequences which are separated by open sets, yet converging geometrically to the same $1$-fold broken trajectory.
However, the latter is impossible in light of the glueing theorem.
See figure \ref{fig:moduli} for a schematic depiction of this situation.
\end{description}

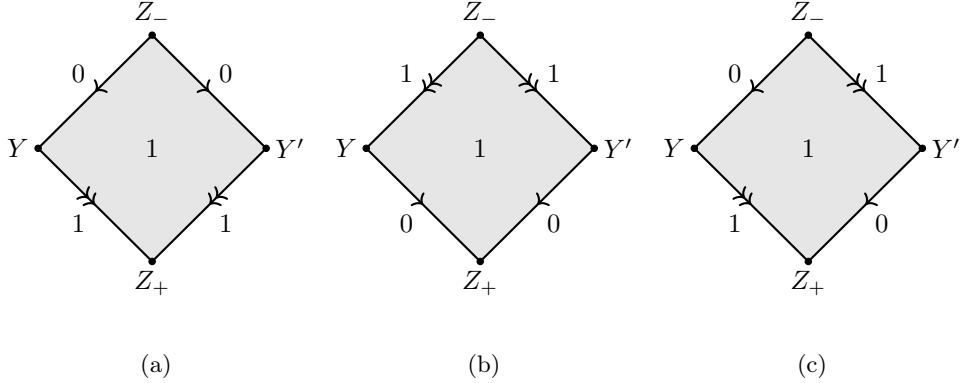
\begin{figure}[t]
  \centering

\begin{subfigure}[b]{.30\linewidth}
  \centering

\makeatletter
\tikzset{nomorepostaction/.code={\let\tikz@postactions\pgfutil@empty}}
\makeatother

\begin{tikzpicture}[scale=1.5]
\coordinate (A1) at (0,1);
\coordinate (A2) at (-1,0);
\coordinate (A3) at (1,0);
\coordinate (A4) at (0,-1);

\fill[fill=black!10!white] (A1)--(A2)--(A4)--(A3);

\draw (0,0) node {$1$};

\begin{scope}[every path/.style={postaction={ nomorepostaction, decorate, decoration={ markings,mark=at position 0.5 with \arrow{>} } } }]
\draw[thick] (A1)-- node[above left] {$0$} (A2);
\draw[thick] (A1)-- node[above right] {$0$} (A3);
\end{scope}

\begin{scope}[every path/.style={postaction={ nomorepostaction, decorate, decoration={ markings,mark=at position 0.5 with \arrow{>>} } } }]
\draw[thick] (A2)-- node[below left] {$1$} (A4);
\draw[thick] (A3)-- node[below right] {$1$} (A4);
\end{scope}

\draw[fill=black] (A1) circle (0.08em) node[above] {$Z_-$};
\draw[fill=black] (A2) circle (0.08em) node[left] {$Y$};
\draw[fill=black] (A3) circle (0.08em) node[right] {$Y'$};
\draw[fill=black] (A4) circle (0.08em) node[below] {$Z_+$};
\end{tikzpicture}
  \caption{}
  \label{fig:traj_n_n0a1}
\end{subfigure}
\begin{subfigure}[b]{.30\linewidth}
  \centering

\makeatletter
\tikzset{nomorepostaction/.code={\let\tikz@postactions\pgfutil@empty}}
\makeatother

\begin{tikzpicture}[scale=1.5]
\coordinate (A1) at (0,1);
\coordinate (A2) at (-1,0);
\coordinate (A3) at (1,0);
\coordinate (A4) at (0,-1);

\fill[fill=black!10!white] (A1)--(A2)--(A4)--(A3);

\draw (0,0) node {$1$};

\begin{scope}[every path/.style={postaction={ nomorepostaction, decorate, decoration={ markings,mark=at position 0.5 with \arrow{>>} } } }]
\draw[thick] (A1)-- node[above left] {$1$} (A2);
\draw[thick] (A1)-- node[above right] {$1$} (A3);
\end{scope}

\begin{scope}[every path/.style={postaction={ nomorepostaction, decorate, decoration={ markings,mark=at position 0.5 with \arrow{>} } } }]
\draw[thick] (A2)-- node[below left] {$0$} (A4);
\draw[thick] (A3)-- node[below right] {$0$} (A4);
\end{scope}

\draw[fill=black] (A1) circle (0.08em) node[above] {$Z_-$};
\draw[fill=black] (A2) circle (0.08em) node[left] {$Y$};
\draw[fill=black] (A3) circle (0.08em) node[right] {$Y'$};
\draw[fill=black] (A4) circle (0.08em) node[below] {$Z_+$};
\end{tikzpicture}
  \caption{}
  \label{fig:traj_n_a1n0}
\end{subfigure}
\begin{subfigure}[b]{.30\linewidth}
  \centering

\makeatletter
\tikzset{nomorepostaction/.code={\let\tikz@postactions\pgfutil@empty}}
\makeatother

\begin{tikzpicture}[scale=1.5]
\coordinate (A1) at (0,1);
\coordinate (A2) at (-1,0);
\coordinate (A3) at (1,0);
\coordinate (A4) at (0,-1);

\fill[fill=black!10!white] (A1)--(A2)--(A4)--(A3);

\draw (0,0) node {$1$};

\begin{scope}[every path/.style={postaction={ nomorepostaction, decorate, decoration={ markings,mark=at position 0.5 with \arrow{>} } } }]
\draw[thick] (A1)-- node[above left] {$0$} (A2);
\draw[thick] (A3)-- node[below right] {$0$} (A4);
\end{scope}

\begin{scope}[every path/.style={postaction={ nomorepostaction, decorate, decoration={ markings,mark=at position 0.5 with \arrow{>>} } } }]
\draw[thick] (A1)-- node[above right] {$1$} (A3);
\draw[thick] (A2)-- node[below left] {$1$} (A4);
\end{scope}

\draw[fill=black] (A1) circle (0.08em) node[above] {$Z_-$};
\draw[fill=black] (A2) circle (0.08em) node[left] {$Y$};
\draw[fill=black] (A3) circle (0.08em) node[right] {$Y'$};
\draw[fill=black] (A4) circle (0.08em) node[below] {$Z_+$};
\end{tikzpicture}
  \caption{}
  \label{fig:traj_n_mixed}
\end{subfigure}

  \caption{Geometric closures of index $1$ connected component of $\calM(Z_-,Z_+)$ in the nonautonomous setting.
  Double arrow heads indicate autonomous trajectories, single arrow heads indicate nonautonomous trajectories,
numbers indicate their indices.}
\label{fig:boundaries_nonautonomous}
\end{figure}

\begin{myremark}
  \label{remark:nonautonomous_moduli}

In this section we analyzed the geometry of (the compactification of) the moduli spaces in the autonomous case.
We now indicate how this analysis can be adapted to the nonautonomous case.
The main technical difference is that one needs the compactness estimates for the nonautonomous case, as given in section \ref{sec:compactness}.
The other difference is the lack of translational invariance of \eqref{eq:TWE}.
In this case we can therefore only use the classification of $1$-dimensional manifolds when the index difference is $\mu(Z_-)-\mu(Z_+) = 1$.
What one obtains is, much the same as in the autonomous setting, that any connected component of $\calM(Z_-,Z_+)$ 
is either compact, or can be compactified by two pairs of broken trajectories.
Here, a nonautonomous trajectory can break into either one of the following (using notation as in hypothesis \ref{C:asymptotic_constant})
\begin{itemize}
\item a concatenation of an index 0 nonautonomous trajectory corresponding to $(f,g,c)$, and an index 1 autonomous trajectory corresponding to $(f_+,g_+,c_+)$, or
\item a concatenation of an index 1 autonomous trajectory corresponding to $(f_-,g_-,c_-)$, and an index 0 nonautonomous trajectory corresponding to $(f,g,c)$.
\end{itemize}
Hence, accounting for multiplicity, there are three possible boundaries for any noncompact connected component of $\calM(Z_-,Z_+)$ when $\mu(Z_-)-\mu(Z_+) = 1$,
see figure \ref{fig:boundaries_nonautonomous}, and \cite{schwarz1993morse} for more detail.
\end{myremark}

\section{Generic properties}\label{sec:genericity}

In this section we show that for generic choices of $f$ and $g$, 
all rest points are hyperbolic and connecting orbits up to order $m-1$ are transversal.
In other words, we will show that regular $f$ and $g$ are generic.
Thus, the results from the preceding sections apply in generic cases.

\subsection{Hyperbolicity}
\label{subsec:generic_hyperbolicity}
As already indicated, hyperbolicity of all rest points of \eqref{eq:TWE} can always be achieved by perturbing the nonlinearity $f$.
We shall now first construct a space from which our generic perturbations of the nonlinearity can be chosen.
Since we will apply the Sard-Smale theorem to a map defined on this space, 
it must be constructed in such a way that it is a Banach manifold.

Given $m \in \N$, let $C_0^m(\Omega \times \R)$ consist of those functions $\phi_0 \in C_b^m(\cl\Omega\times\R)$
such that $\lim_{u\to\pm\infty} \phi_0(x,u) = 0$ uniformly in $x \in \Omega$, and $\rst{ \phi_0 }{ \bdy \Omega \times \R } = 0$.
Equipped with the norm inherited from $C_b^m(\cl \Omega\times\R)$, this becomes a separable Banach space.
Now let $\calF^m$ consist of those functions $\phi$ of the form
\[
\phi(x,u) = e^{-|u|^2} \phi_0(x,u), \qquad \text{where} \quad \phi_0 \in C_0^m(\Omega \times \R).
\]
Define a norm on $\calF^m$ by setting $\|\phi\|_{\calF^m} := \| \phi_0 \|_{C_b^m}$.
As such, $\calF^m$ is isometric to $C_0^m(\Omega\times\R)$, hence $\calF^m$ is a separable Banach space.
Note that the rapid decay of $\phi \in \calF^m$ implies that $f+\phi$ satisfies hypotheses \Fgroup whenever $f$ does.

Recall that a subset of a topological space is called residual if it can be written as the countable intersection of open and dense subsets.
Since $\calF^m$ is a Baire space, any residual subset of $\calF^m$ is also dense in $\calF^m$.
\begin{mythm}\label{thm:genericity_hyperbolicity}
  Fix a nonlinearity $f \in C^m(\Omega \times \R)$, with $m \geq 1$, satisfying hypotheses \Fgroup.
Then there exists a residual set $\calF^m_\text{reg} \subset \calF^m$ such that for each $\phi \in \calF^m_{\text{reg}}$
equation \eqref{eq:TWE} with the perturbed nonlinearity $f + \phi$ has only hyperbolic rest points.
\end{mythm}
\begin{myproof}
  Given a nonlinearity $f$, for notational convenience we will write $A_f$ instead of $A_{f,g,c}$ 
(the construction of $\calF^m_{\text{reg}}$ will in fact be independent of $g$ and $c$).
Consider the map
\begin{align*}
  & \Psi : \calF^m \times X^1 \to X^0, \\
  & \Psi(\phi,Z) = A_{f + \phi}(Z).
\end{align*}
Note that this $\Psi$ is of class $C^m$.

We claim that $0$ is a regular value of $\Psi$, 
i.e.\  for any $(\psi,Z) \in \Psi^{-1}(0)$ the linear operator $\d \Psi(\psi,Z) \in \calL(\calF^m\times X^1,X^0)$ has a continuous right inverse.
Assume without loss of generality that $\Psi^{-1}(0) \neq \emptyset$, and pick any $(\phi,Z) \in \Psi^{-1}(0)$.
For any $(\psi,W) \in \calF^m \times X^1$, one has
\[
\d \Psi(\phi,Z)(\psi,W) = \d A_{f+\phi}(Z)W + B(Z)\psi,
\]
where
\[
B(Z) \psi(x) =
\begin{pmatrix}
  0 \\
  \psi(x,z(x))
\end{pmatrix}, \qquad \text{where} \quad Z = (z,0).
\]

Note that by the Rellich-Kondrachov theorem the operator $\d A_{f+\phi}(Z)$ 
is a compact perturbation of $(u,v) \mapsto (-v , \Delta u) \in \calL(X^1_\C,X^0_\C)$.
The latter is Fredholm of index $0$, hence $\d A_{f+\phi}(Z)$ is also Fredholm of index $0$. 
Hence $\d \Psi(\phi,Z)$ is the direct sum of a Fredholm map and a bounded map, 
therefore the existence of a bounded right inverse will follow from surjectivity of $\d \Psi(\phi,Z)$.

Note that $\img \d \Psi(\phi,Z) \supset \img \d A_{f+\phi}(Z)$ is finite codimensional, and consequently $\d \Phi(\psi,z)$ has closed range.
It thus suffices to see that $\img \d \Psi(\phi,Z)$ is dense in $X^0$.
For this, we note that by the compactness results from section \ref{sec:compactness} the function $z$ is continuous and uniformly bounded.
Thus with $\psi$ defined by $\psi(x,u) = \xi(x) \eta(u)$, 
where we choose $\xi \in C^\infty_c(\Omega)$ and $\eta \in C_c^\infty(\R)$ with $\eta(u) = 1$ on $\supp( z )$,
it follows that $B(Z) \psi = (0,\xi)$.
Hence $\{0\} \times C^\infty_c(\Omega) \subset \img \d \Psi(\phi,Z)$.
Also, because of the shape of $\d A_{f+\phi}(Z)$, there exists a subspace $E \subset L^2(\Omega)$ such that $H^1_B(\Omega) \times E\subset \img \d \Psi(\phi,Z)$.
From this we see that $\img \d \Psi(\phi,Z)$ is dense in $X^0$, thus showing that $0$ is a regular value of $\Psi$.

Let $Z \in \calS(f+\phi)$ be a stationary solution of \eqref{eq:TWE} with nonlinearity $f + \phi$.
We claim that $\d A_{f+\phi}(Z)$ is surjective precisely when $Z$ is a hyperbolic solution of \eqref{eq:TWE}.
Indeed, for any $\mu \in \C$ it follows by the Rellich-Kondrachov theorem that $\d A_{f+\phi}(Z) - \mu$ is a compact perturbation of $\d A_{f+\phi}(Z)$,
hence Fredholm of index $0$.
Therefore the spectrum of $\d A_{f+\phi}(Z)$ consists solely of eigenvalues.
Consequently, $\mu \in \sigma(\d A_{f+\phi}(Z))$ if and only if $\mu^2 - c \mu \in \sigma_p(P)$, 
where $P = \Delta + f_u(x,z) + \phi_u(x,z)$ as an unbounded operator on $L^2(\Omega)_\C$ with domain $\calD(P) = H^2_B(\Omega)_\C$.
As such $P$ is symmetric, hence it has real-valued point spectrum; i.e.\  $\sigma_p(P) \subset \R$.
If $\mu = i\lambda \in i\R \setmin \{0\}$, then $\mu^2 - c \mu \not\in \R$,
thus proving that $i \R \setmin \{0\} \cap \sigma(\d A(Z)) = \emptyset$.
From this the claim follows.

Now let $\calF^m_{\text{reg}}$ consist of those $\phi \in \calF^m$ for which $0$ is a regular value of $\Psi(\phi,\cdot)$.
We will now argue that $\calF^m_{\text{reg}}$ is residual in $\calF^m$.
By the implicit function theorem, $\calZ := \Psi^{-1}(0)$ is a $C^m$ manifold.
Let $\pi : \calZ \to \calF^m$ be the restriction to $\calZ$ of the projection $\Pr_1 : \calF^m \times X^1 \to \calF^m$,
i.e.\  $\pi = \rst{ \Pr_1 }{ \Psi^{-1}(0) }$, and note that this is a $C^m$ map.
Moreover note that $\phi$ is a regular value of $\pi = \rst{ \Pr_1 }{ \Psi^{-1}(0) }$ 
if and only if $0$ is a regular value of $\Psi(\phi, \cdot) = \rst{ \Psi }{ \Pr_1^{-1}(\phi) }$.
Pick any $\phi \in \calF^m$, and note that $\phi$ is a regular value of $\Pr_1$.
The linearization of $\rst{ \Psi }{ \Pr_1^{-1}(\phi) } = \Psi(\phi,\cdot)$ around $Z$ equals $\d A_{f+\phi}(Z)$, 
hence the former map is Fredholm of index $0$.
Then also $\pi$ is Fredholm of index $0$.
Consequently, the Sard-Smale theorem \cite{smale1965infinite} 
implies that the regular values of $\pi$ are residual in $\calF^m$, thus proving the claim.
\end{myproof}

\subsection{Transversality}

We will show here that there is an abundance of $g$ for which the transversality condition up to order $m-1$ is satisfied.
To do so we must first introduce a separable Banach manifold $\calG^m(c)$ from which our perturbations $g$ can be chosen.

The following lemma indicates how we can build localized perturbations $g$.
We stress that the proof of the lemma relies heavily on the unique continuation theorem (theorem \ref{thm:unique_continuation}).
\begin{mylemma}
\label{lemma:psi_localization}
Fix any $\psi \in C^\infty_c(\R)$ with $\int_0^\infty \psi(t) \d t \neq 0$.
Suppose $f$ is of class $C^4$ and suppose hypothesis \Fgroup are satisfied.
Let $U$ be a non-stationary bounded solution of \eqref{eq:TWE}.
For any $\epsilon > 0$ and $t_* \in \R$, define $\Lambda_\epsilon \in C_b^1(\R)$ by
\[
\Lambda_\epsilon(t) := \epsilon^{-1} \psi\left( \epsilon^{-1} \| U(t) - U(t_*) \|_{X^0} \right).
\]
Set $C := 2 \|\partial_t U(t_*) \|_{X^0}^{-1} \int_0^\infty \psi(t) \d t$.
Then for any $\chi \in C^0(\R)$ one has
\[
\lim_{\epsilon \to 0} C^{-1} \int_\R \chi(t) \Lambda_\epsilon(t) \d t  = \chi(t_*).
\]
\end{mylemma}
\begin{myproof}
  Let $a > 0$ be such that $\supp(\psi) \subset [-a,a]$.
  Given $\epsilon > 0$, consider the function $\xi^+_\epsilon : [t_*,\infty) \to [0,\infty)$, 
  defined by $\xi^+_\epsilon(t) := \epsilon^{-1} \| U(t) - U(t_*) \|_{X^0}$.
  Note that, for $t \in (t_*,\infty)$,
  \[
  \ddxe{t} \xi^+_\epsilon(t) = \frac{ 1 }{ \epsilon } \frac{ \langle \partial_t U(t) , U(t) - U(t_*) \rangle_{X^0} }{ \| U(t) - U(t_*) \|_{X^0} },
  \]
  where $\langle \cdot , \cdot \rangle_{X^0}$ denotes the Hilbert space inner product on $X^0$.
  We claim that there exists $\delta > 0$ and $\epsilon_0 > 0$ such that 
  \begin{enumerate}
  \item $\ddxe{t} \xi^+_\epsilon(t) > 0$ for $t \in (t_*,t_*+\delta)$ and $0 < \epsilon \leq \epsilon_0$, and
  \item $\xi^+_\epsilon(t) \geq a$ for any $t \in [t_*+\delta,\infty)$ and $0 < \epsilon \leq \epsilon_0$.
  \end{enumerate}

  To see why (1) holds, note that
  \[
    \ddxe{t} \xi^+_\epsilon(t) = 
    \frac{ 1 }{ \epsilon } \frac{ \left\langle \partial_t U(t) , 
                                 \dfrac{ U(t) - U(t_*) }{ t - t_* } \right\rangle_{X^0} }{ \left\| \dfrac{ U(t) - U(t_*) }{ t - t_* } \right\|_{X^0} } 
    \to \frac{ 1 }{ \epsilon } \| \partial_t U(t_*) \|_{X^0} \qquad \text{as} \quad t \downarrow t_*.
  \]
  Since $U$ is assumed to be non-stationary, the unique continuation theorem \ref{thm:unique_continuation}
  (here we use that $f$ is $C^4$) ensures that $\|\partial_t U(t_*)\|_{X^0} > 0$.
  Hence, by continuity, there must exist a $\delta > 0$ such that claim (1) holds (with any $\epsilon > 0$).

  We will now consider claim (2).
  Let $\delta > 0$ be as in the preceding paragraph, and suppose an $\epsilon_0$ such that claim (2) holds does not exist.
  Then we would be able to find sequences $(t_n)_n \subset [t_*+\delta,\infty)$ and $(\epsilon_n)_n \subset (0,\infty)$
  with $\epsilon_n \to 0$ as $n\to\infty$, such that $\xi^+_{\epsilon_n}(t_n) < a$ for all $n$.
  Hence $\| U(t_n) - U(t_*) \|_{X^0} < a \epsilon_n \to 0$ as $n \to \infty$.

  We claim that the sequence $(t_n)_n$ is convergent, with $\lim_{n\to\infty} t_n = t_*$.
  This is obviously in contradiction with the construction of the sequence, hence this will prove claim (2).
  To see why $t_n \to t_*$ as $n \to \infty$, we first note that $(t_n)_n$ must be bounded.
  Indeed, since we assume $U$ to be bounded, lemma \ref{lemma:gradient_like} implies that it is a connecting orbit.
  Hence if we could find an unbounded subsequence $(t'_n)_n$ of $(t_n)_n$,
  then $U(t'_n) \to Z \in \calS(f)$.
  But by assumption $U(t_n) \to U(t_*)$, hence $U(t_*) = Z$.
  The unique continuation theorem \ref{thm:unique_continuation} would then imply that $U(t) = Z$ for all $t \in \R$,
  which contradicts our assumption that $U$ is non-stationary.
  We have thus proved that $(t_n)_n$ is bounded.
  Given any subsequence, select a further subsequence $(t'_n)_n$ which is convergent, say $t'_n \to t_\infty \in \R$.
  Then $U(t_\infty) = \lim_{n\to\infty} U(t'_n) = U(t_*)$,
  hence if $t_\infty \neq t_*$ the unique continuation theorem \ref{thm:unique_continuation} would imply that
  $U$ is periodic.
  However, such behaviour is excluded by the gradient-like structure of \eqref{eq:TWE}.
  Hence $t_\infty = t_*$, which proves that $(t_n)_n$ is convergent, with $\lim_{n\to\infty} t_n = t_*$.
  This proves claim (2).

  The implicit function theorem now ensures the existence of a family of maps $t^+_\epsilon : [0,a) \to [t_*,t_*+\delta)$ (with $0 < \epsilon \leq \epsilon_0$),
  which restrict to $C^1$ diffeomorphisms from $(0,a)$ onto their image, such that
  \[
    \xi^+_\epsilon \circ t^+_\epsilon = \id_{[0,a)},
  \]
  \[
  \big( t^+_\epsilon \big)^{-1}(t_*) = 0,
  \]
  and
  \[
    \xi^+_\epsilon\big( [t_*,\infty) \setmin t^+_\epsilon([0,a)) \big) \cap [0,a) = \emptyset.
  \]
  Furthermore, since $\xi^+_\epsilon(t) \to \infty$ as $\epsilon \to 0$ when $t \neq t_*$, and $\xi^+_\epsilon(t_*) = 0$,
  it follows that $t^+_\epsilon(s) \downarrow t_*$ as $\epsilon \to 0$.

  Similarly, we define the family of maps $\xi^-_\epsilon : (-\infty,t_*] \to [0,\infty)$ by $\xi^-_\epsilon(t) := \epsilon^{-1} \| U(t) - U(t_*) \|_{X^0}$.
  The same argument as above then proves the existence of a family $t^-_\epsilon : [0,a) \to (-\infty,t_*]$ with the same properties as the maps $t^+_\epsilon$.

  At this point we are prepared to compute the limit of $\Lambda_\epsilon$ as $\epsilon \to 0$.
  Fix any $\chi \in C^0(\R)$.
  Then, since $\Lambda_\epsilon(t) = 0$ for $t \not\in t^-_\epsilon([0,a)) \cup t^+_\epsilon([0,a))$, one has
  \begin{align*}
    \int_\R \chi(t) \Lambda_\epsilon(t) \d t &= \int_{\R\setmin\{0\}} \chi(t) \Lambda_{\epsilon}(t) \d t 
                                               = \int_{t^-_\epsilon((0,a))} \chi(t) \Lambda_\epsilon(t) \d t + \int_{t^+_\epsilon((0,a))} \chi(t) \Lambda_\epsilon(t) \d t \\
    &= \int_0^a \chi(t^-_\epsilon(s)) \Lambda_\epsilon(t^-_\epsilon(s)) \ddx{t^-_\epsilon(s)}{s} \d s 
      + \int_0^a \chi(t^+_\epsilon(s)) \Lambda_\epsilon(t^+_\epsilon(s)) \ddx{t^+_\epsilon(s)}{s} \d s.
  \end{align*}
  Consider for example the last integral.
  Filling in the definition of $\Lambda_\epsilon$ gives
  \[
    \int_0^a \chi(t^+_\epsilon(s)) \Lambda_\epsilon(t^+_\epsilon(s)) \ddx{t^+_\epsilon(s)}{s} \d s =
      \int_0^a  \chi(t^+_\epsilon(s)) \psi(s)
        \frac{ \| U(t^+_\epsilon(s)) - U(t_*) \|_{X^0} }{ \langle \partial_t U(t_*) , U(t^+_\epsilon(s)) - U(t_*) \rangle_{X^0} } \d s.
  \]
  Since $t^+_\epsilon(s) \downarrow t_*$ as $\epsilon \to 0$, we have
  \[
  \chi(t^+_\epsilon(s)) \to \chi(t_*)
  \]
  as $\epsilon \to 0$, and
  \[
  \dfrac{ \| U(t^+_\epsilon(s)) - U(t_*) \|_{X^0} }{ \langle \partial_t U(t_*) , U(t^+_\epsilon(s)) - U(t_*) \rangle_{X^0} }
  = \frac{ \bigg\| \dfrac{ U(t^+_\epsilon(s)) - U(t_*) }{ t^+_\epsilon(s) - t_* } \bigg\|_{X^0} }
  { \bigg\langle \partial_t U(t_*) , \dfrac{ U(t^+_\epsilon(s)) - U(t_*) }{ t^+_\epsilon(s) - t_* } \bigg\rangle_{X^0} }
  \to \| \partial_t U(t_*) \|_{X^0}^{-1}
  \]
  as $\epsilon \to 0$.
  Using the dominated convergence theorem it follows that
  \[
  \int_0^a \chi(t^+_\epsilon(s)) \Lambda_\epsilon(t^+_\epsilon(s)) \ddx{t^+_\epsilon(s)}{s} \d s \to \frac{C}{2} \chi(t_*) \qquad \text{as} \quad \epsilon \to 0.
  \]
  Similarly, we find that
  \[
  \int_0^a \chi(t^-_\epsilon(s)) \Lambda_\epsilon(t^-_\epsilon(s)) \ddx{t^-_\epsilon(s)}{s} \d s \to \frac{C}{2} \chi(t_*) \qquad \text{as} \quad \epsilon \to 0.
  \]
  This concludes the proof of the lemma.
\end{myproof}

Let $\widetilde \calG^m$ consist of all $g \in C_b^m(\Omega \times X^0)$ 
which vanish on $X^0_S := \Omega \times ( H^1_B(\Omega) \times \{0\} )$.
Equipped with the $C_b^m(\Omega\times X^0)$-norm $\widetilde \calG^m$ becomes a Banach space.
We now construct a suitable separable subspace of $\widetilde \calG^m$.
Fix $\psi \in C^\infty_c(\R)$ as in lemma \ref{lemma:psi_localization}.
Consider $g \in \widetilde \calG^m$ of the form
\begin{equation}
  \label{eq:g_representative}
  g(x,U) = \phi(x) \psi\left(\frac{1}{\epsilon} \| U - U_0 \|_{X^0} \right),
\end{equation}
where $\phi \in C^\infty(\cl \Omega)$,
$U_0 \in X^0 \setmin X^0_S$, and $\epsilon > 0$ is sufficiently small (depending on $U_0$).
Denote by $\calG^m$ the smallest closed subspace of $\widetilde \calG^m$ which contains all maps of this form.
Since $C_b^m(\Omega)$, $\R$, and $X^0$ are all separable, and the map
\begin{align*}
  C_b^m(\Omega) \times \R \times X^0 &\to \widetilde \calG^m , \\
  (\phi , \epsilon , U_0 ) &\mapsto g
\end{align*}
is continous, the space $\calG^m$ is a separable Banach space.

Now let $\calG^m(c)$ consist of those $g \in \calG^m$ for which 
\[
\sup_{x\in\Omega,\; u \in H^1_B(\Omega)} | g(x,u,v) | < \frac{ c }{ 2 \sqrt{ \Vol(\Omega) } } \| v \|_{L^2(\Omega)} \qquad \text{whenever} \quad v \neq 0.
\]
Note that $\calG^m(c)$ is an open subspace of $\calG^m$.
Also remark that any $g \in \calG^m(c)$ satisfies hypotheses \Ggroup.

\begin{mythm}\label{thm:genericity_transversality}
  Let $f$ be of class $C^m$ with $m \geq 4$ and suppose hypotheses \Fgroup are satisfied.
  Furthermore, assume that each stationary solution of \eqref{eq:TWE} is hyperbolic.
  Then there exists a residual set $\calG^m_\text{reg} \subset \calG^m(c)$ such that for any $g \in \calG^m_\text{reg}$
  the transversality condition up to order $m-1$ is satisfied.
\end{mythm}
\begin{myproof}
  Fix $Z_-,Z_+ \in \calS$ with $\mu(Z_-) - \mu(Z_+) \leq m-1$, and consider the map
  \begin{align*}
    & \Psi_{Z_-,Z_+} : \calG^m(c) \times \calP(Z_-,Z_+) \to L^2(\R;X^0), \\
    & \Psi_{Z_-,Z_+}(g,U) = \partial_t U + A_{f,g,c}(U).
  \end{align*}
  Here $\calP(Z_-,Z_+)$ is the affine space defined in section \ref{sec:moduli}.
  With the aid of theorem \ref{thm:operator_regularity} we see that $\Psi_{Z_-,Z_+}$ 
  is of class $C^m$ in a neighbourhood of $\Psi_{Z_-,Z_+}^{-1}(0)$.

We will argue that $0$ is a regular value of $\Psi_{Z_-,Z_+}$.
Assume without loss of generality that $\Psi_{Z_-,Z_+}^{-1}(0) \neq \emptyset$, and pick any $(g,U) \in \Psi_{Z_-,Z_+}^{-1}(0)$.
Then
\[
\d \Psi_{Z_-,Z_+}(g,U)(\gamma,V) = \dsub 2 \Psi_{Z_-,Z_+}(g,U) V + B(U) \gamma,
\]
where $B(U) \in \calL(\calG^m,L^2(\R;X^0))$ is given by
\[
B(U) \gamma =
\begin{pmatrix}
  0 \\
  \gamma(\cdot,U)
\end{pmatrix}.
\]
We know from section \ref{sec:fredholm} that $\dsub 2 \Psi_{Z_-,Z_+}(g,U) \in \calL(W^{1,2}(\R;X^0,X^1),L^2(\R;X^0))$ 
is Fredholm of index $\mu(Z_-) - \mu(Z_+)$.
Therefore $\d \Psi_{Z_-,Z_+}(g,U)$ has a bounded right inverse as soon as it is surjective.
Furthermore, $\img{ \d \Psi_{Z_-,Z_+}(g,U) }$ is finite codimensional hence closed in $L^2(\R;X^0)$.
So it suffices to check that $\img{ \d \Psi_{Z_-,Z_+}(g,U) }$ is dense in $L^2(\R;X^0)$.

Select any $(\xi_1,\xi_2) \in ( \img{ \d \Psi_{Z_-,Z_+}(g,U) } )^\perp$.
Then
\begin{align}
  \label{eq:perp1}
  \left\langle
  \begin{pmatrix}
    \xi_1 \\
    \xi_2
  \end{pmatrix}
  ,
  \dsub 2 \Psi_{Z_-,Z_+}(g,U) V
\right\rangle_{L^2(\R;X^0)} &= 0 \qquad \text{for all} \quad V \in W^{1,2}(\R;X^0,X^1), \\
\label{eq:perp2}
\left\langle
  \begin{pmatrix}
    \xi_1 \\
    \xi_2
  \end{pmatrix}
  ,
  B(U)\gamma
\right\rangle_{L^2(\R;X^0)} &= 0 \qquad \text{for all} \quad \gamma \in \calG^m.
\end{align}

Set $\xi_1^*(t) := \langle \xi_1(t) , \cdot \rangle_{H^1_B(\Omega)}$ and $\xi_2^*(t) := \langle \xi_2(t) , \cdot \rangle_{L^2(\Omega)}$,
so that $\xi_1^* \in L^2(\R;H^1_B(\Omega)^*)$ and $\xi_2^* \in L^2(\R;L^2(\Omega)^*)$.
Combining equation \eqref{eq:perp1} with the regularity results from \cite{rabier2004robbin} show that
\[
(\xi_1^*,\xi_2^*) \in W^{1,2}(\R;H^2_B(\Omega)^* \times H^1_B(\Omega)^* , H^1_B(\Omega)^* \times L^2(\Omega)^* )
\]
and
\[
(\xi_1^*,\xi_2^*) \in \ker \big( -\partial_t + \d A_{f,g,c}(U)^* \big),
\]
which means that
\begin{equation}
  \label{eq:coker_eqs}
\begin{cases}
  - \partial_t \xi_1^* + \Delta^* \xi_2^* + f'(\cdot,u)^* \xi_2^* + \partial_1 g(\cdot,U)^* \xi_2^* = 0 \\
  - \partial_t \xi_2^* - \xi_1^* - c \xi_2^* + \partial_2 g(\cdot,U)^* \xi_2^* = 0.
\end{cases}
\end{equation}
Here the adjoints are to be considered as the dual operators of bounded operators between Banach spaces,
where $\Delta : H^2_B(\Omega) \to L^2(\Omega)$, $f'(\cdot,u) : H^1_B(\Omega) \to L^2(\Omega)$,
 $\partial_1 g(\cdot,(u,v)) = \ddi{ g(\cdot,(u,v)) }{ u } : H^1_B(\Omega) \to L^2(\Omega)$, and $\partial_2 g(\cdot,(u,v)) = \ddi{ g(\cdot,(u,v)) }{ v } : L^2(\Omega) \to L^2(\Omega)$.
We shall be using this observation later on.

Equation \eqref{eq:perp2} implies that for all $\gamma \in \calG^m$ it holds that
\[
\int_\R \int_\Omega \gamma(x,U(t)) \xi_2(t,x) \d x \d t = 0.
\]
In particular, consider $\gamma = \gamma_\epsilon$ of the form
\[
\gamma_{\epsilon}(x,V) := \phi(x) \psi\left(\frac{1}{\epsilon} \| V - U(t_0) \|_{X^0} \right),
\]
where $t_0 \in \R$, and $\phi \in C_c^\infty(\Omega)$, and $\psi$ is as in \eqref{eq:g_representative}.
Since $\xi_2^* \in W^{1,2}(\R;H^1_B(\Omega)^* , L^2(\Omega)^*)$ it follows that $\xi_2^* \in C^0_b(\R;L^2(\Omega)^*)$, hence the map
\[
t \mapsto \int_\Omega \phi(x) \xi_2(t,x) \d x
\]
is continuous.
By lemma \ref{lemma:psi_localization}
\[
0 = \lim_{\epsilon \to 0} \epsilon^{-1} \int_\R \int_\Omega \gamma_\epsilon(x,U(t)) \xi_2(t,x) \d x \d t
= \frac{ \displaystyle 2 \int_0^\infty \psi(s) \d s }{ \| \partial_t U(t_0) \|_{X^0} } \int_\Omega \phi(x) \xi_2(t_0,x) \d x.
\]
This holds for all $t_0 \in \R$ and $\phi \in C_c^\infty(\Omega)$, hence $\xi_2 = 0$.
Consequently, equation \eqref{eq:coker_eqs} implies that $\xi_1 = 0$ as well.
This shows that $( \img{ \d \Psi_{Z_-,Z_+}(g,U) } )^\perp = \{0\}$, hence $E$ is dense in $L^2(\R,X^0)$, as desired.

By the implicit function theorem the regular level set 
$\calZ := \Psi_{Z_-,Z_+}^{-1}(0)$ is a $C^m$ smooth submanifold of $\calG^m(c) \times \calP(Z_-,Z_+)$.
The projection $\pi : \calZ \to \calG^m(c)$ is $C^m$ smooth and Fredholm of index $\ind(\pi) = \mu(Z_-) - \mu(Z_+) \leq m-1$.
Applying the Sard-Smale theorem \cite{smale1965infinite} to this map (here we use that $m \geq \max\{0,\ind(\pi)\} + 1$) 
and using a transversality argument similar to the one in the proof of theorem \ref{thm:genericity_hyperbolicity} 
we find a residual subset $\calG^m_{\text{reg}}(c ; Z_-,Z_+) \subset \calG^m(c)$ such that 
for each $g \in \calG^m_{\text{reg}}(c ; Z_-,Z_+)$ orbits connecting $Z_-$ and $Z_+$ are transversal.
Now set
\[
\calG^m_{\text{reg}}(c) := \bigcap_{\substack{ Z_-,Z_+ \in \calS \\ \mu(Z_-) - \mu(Z_+) \leq m-1}} \calG^m_{\text{reg}}(c ; Z_-,Z_+).
\]
By the compactness results from section \ref{sec:compactness} and hyperbolicity of the rest points,
this is a countable intersection, hence $\calG^m_{\text{reg}}(c)$ is residual.
\end{myproof}

\section{The travelling wave homology}\label{sec:homology}

\subsection{The homology for generic perturbations}
Given $N \subset X^0$, denote by $\BInv(N;f,g,c)$ the set of all points $U(t) \in X^0$, with $t \in \R$, 
where $U$ is a solution of \eqref{eq:TWE} for which $U(t) \in N$ for all $t \in \R$,
and for which $\sup_{t\in\R} \| U(t) \|_{X^0} < \infty$.
We will call $\BInv(N;f,g,c)$ the bounded invariant set of $N$.
In light of lemma \ref{lemma:gradient_like}, if $f$ is a regular nonlinearity
the set $\BInv(N;f,g,c)$ will consists solely of stationary solutions and connecting orbits.
The set $N$ shall be called an \emph{isolating neighbourhood} if
\begin{enumerate}
\item $N$ is closed in $X^0$,
\item $\BInv(N;f,g,c) \subset \inter N$, i.e.\  $N$ isolates the rest points and connecting orbits.
\end{enumerate}
Note that our definition of an isolating neighbourhood differs from the usual one since we allow $N$ to be unbounded.
Let $\calS(N,f) := \calS(f) \cap N$, and given a normalized Morse index $\mu$, define
\[
\calS_k(N,f) := \set{ Z \in \calS(N,f) }{ \mu_f(Z) = k }.
\]
An isolating neighbourhood $N$ shall be called \emph{finitely generating} provided that for each $k \in \Z$, the set $\calS_k(N,f)$ is finite.
Note that this notion is independent of the chosen normalized Morse index $\mu$.

The chain complex will depend on the choices of a finitely generating isolating neighbourhood $N$, 
a triplet $(f,g,c)$ satisfying hypotheses \Hgroup and for which $f$ and $g$ are both regular 
(henceforth $(f,g,c)$ shall also be called a \emph{regular triplet}),
the chosen boundary data $B$, and the chosen normalized Morse index $\mu$.
Assume $f$ and $g$ are at least $C^4$ smooth; this ensures that all the results from preceding sections can be applied in this setting.

Define the graded $\Z_2$-module
\[
C := \bigoplus_{n \in \Z} C_n, \qquad \text{where} \quad C_n := \bigoplus_{ X \in \calS_n(N,f) } \Z_2 \langle X \rangle.
\]
Since the chain groups $C_n$ are finitely generated,
we can define a homomorphism $\partial_n : C_n \to C_{n-1}$ by declaring its action on the generators of $C_n$ as follows: 
for each $X \in \calS_n(N,f)$, we set
\[
  \partial_n X := \sum_{ Y \in \calS_{n-1}(N,f)} i(X,Y) Y,
\]
where
\[
i(X,Y) := \# \set{ [U] \in \widehat \calM(X,Y) }{ U(t) \in N \text{ for all } t \in \R } \modn 2.
\]
From section \ref{sec:moduli} we know that $\widehat \calM(X,Y)$ is a finite collection, hence the number $i(X,Y)$ is well-defined.
The sum in the definition of $\partial_n$ is always finite since $N$ is finitely generating.

We now arrive at one of our main theorem.
\begin{mythm}
\label{thm:homology_definition}
One has $\partial_{n} \circ \partial_{n+1} = 0$, and consequently,
\[
\HTW_n(N,f,g,c) := \homology_n(C_*,\partial_*) = \frac{ \ker{ \partial_n } }{ \img{ \partial_{n+1} } }
\]
is well-defined.
\end{mythm}
\begin{myproof}
The homomorphism $\partial_n \circ \partial_{n+1}$ counts (modulo 2) the 1-fold broken orbits, consisting of two index $1$ orbits,
which are entirely contained in $N$.
We have seen that 1-fold broken orbits consisting of two index $1$ orbits are always the limit of a $1$-parameter family of index 2 orbits,
and therefore always appear in pairs.
Since $N$ is an isolating neighbourhood, $1$-parameter families of solutions with fixed endpoints are trapped by $N$,
hence pairs of 1-fold broken orbits are also trapped by $N$.
Since all counting is done modulo 2, it follows that $\partial_n \circ \partial_{n+1} \equiv 0$.
\end{myproof}

We will call $\HTW_*(N,f,g,c)$ the travelling wave homology.

\subsection{Invariance}
We now want to study what happens to the chain complex and the homology if we perturb $N$, $f$, $g$, or $c$.
First we introduce conditions under which perturbations of $N$ do not change the homology.
We will then show that the homology is independent of the choice of $g$,
thus allowing the definition of $\HTW_*(N,f,c)$.
Finally, we give criteria under which a homotopy $t \mapsto (f_t,c_t)$ induces an isomorphism on the level of homology.

\subsubsection{Perturbing \texorpdfstring{$N$}{N}}
Let $N$, $f$, $g$, and $c$ satisfy the same conditions as in the previous section.
Clearly the construction of the homology does not depend explicitly on $N$, but only on $\BInv(N;f,g,c)$.
Hence if $\widetilde N$ is another isolating neighbourhood
such that $\BInv(N;f,g,c) = \BInv(\widetilde N;f,g,c)$, then $\HTW_*(\widetilde N,f,g,c)$ is well-defined
and $\HTW_*(N,f,g,c) = \HTW_*(\widetilde N,f,g,c)$.

\subsubsection{Independence of \texorpdfstring{$g$}{g}}
Begin by fixing a regular nonlinearity $f$ which satisfies hypotheses \Fgroup.
Also fix $c > 0$.
Assume $N \subset X^0$ is a finitely generating isolating neighbourhood for $(f,0,c)$.
Since $N$ is finitely generating, the set
\[
\calE_f( \calS(N,f) ) = \bigcup_{k\in\Z} \calE_f( \calS_k(N,f) )
\]
is countable, hence the set of regular energy levels
\[
E_{\text{reg}}(N,f) := \R \setmin \calE_f( \calS(N,f) )
\]
is dense in $\R$.
Then, for each $E \in E_{\text{reg}}(N,f)$, the set
\[
N^E := N \cap \calE_f^{-1}((-\infty,E])
\]
is also a finitely generating isolating neighbourhood for $(f,0,c)$.

Recall from section \ref{sec:genericity} the definition of the space $\calG^m(c)$.
Denote by $\calB(\delta)$ the set of those $g \in \calG^m(c)$ for which $\| g \|_{C_b^m(\Omega \times X^0)} < \delta$.
The following lemma now guarantees that $N^E$ remains an isolating neighbourhood when we consider small perturbations in $g$.

\begin{mylemma}
\label{lemma:stable_perturbations}
Let $N$, $f$, and $c$ be as described above.
Then for each $E \in E_{\text{reg}}(N,f)$ there exists $\delta_E > 0$ such that the following is true:
for any curve $t \mapsto g_t$ with values in $\calB(\delta_E)$ which is constant on $(-\infty,-1)$ and on $(1,\infty)$,
the set $N^E$ is an isolating neighbourhood for the nonautonomous equation \eqref{eq:TWE} associated with $(f,g_t,c)$.
\end{mylemma}
\begin{myproof}
Suppose the claim is false.
Then one could find a $E \in E_{\text{reg}}(N,f)$, a sequence $(\delta_n)_n$ with $\lim_{n\to\infty} \delta_n = 0$,
and a sequence of curves $t \mapsto g_t^n$ with values in $\calB(\delta_n)$ 
which are constant for $t \in (-\infty,-1)$ and for $t \in (1,\infty)$,
such that the following holds:
for each $n$ there is a bounded solution $U_n$ to the nonautonomous equation associated to $(f,g_t^n,c)$
such that $U_n(t) \in N^E$ for all $t \in \R$, 
but $U_n(t_n) \in \bdy{ N^E }$ for a certain $t_n$.

Combining lemmata \ref{lemma:energy_bounds_stationary} and \ref{lemma:gradient_like}
it follows that there exists an $M \geq E$ such that $| \calE_f(U_n(t)) | \leq M$ for all $n \in \N$ and $t \in \R$.
As discussed in remark \ref{remark:Palais-Smale_type_compactness}
one then finds that $U_n(\cdot + t_n)$ converges in $W^{1,2}_{\text{loc}}(\R;X^0,X^1)$ over a subsequence to a bounded solution $U$
of the autonomous equation associated with $(f,0,c)$.
This solution is entirely contained in $N^E$, but $U(0) \in \bdy{ N^E }$.
Since $E$ was chosen to be a regular energy level,
this contradicts the assumption that $N$ is an isolating neighbourhood for \eqref{eq:TWE} with $(f,0,c)$.
\end{myproof}

Let us fix, for the moment, a value $E \in E_{\text{reg}}(N,f)$ and a corresponding $\delta_E > 0$ 
as dictated by lemma \ref{lemma:stable_perturbations}.
Since $\calB(\delta_E)$ is open in $\calG^m(c)$, it follows from theorem \ref{thm:genericity_transversality}
that the collection $\calB_{\text{reg}}(\delta_E)$ of regular $g \in \calB(\delta_E)$ is dense in $\calB(\delta_E)$.
For each $g \in \calB_{\text{reg}}(\delta_E)$, the triplet $(f,g,c)$ is regular,
and $N^E$ is a finitely generating isolating neighbourhood for $(f,g,c)$.
Hence, the homology $\HTW_*(N^E,f,g,c)$ is well-defined.

\paragraph{The isomorphism induced by homotopies of $g$.}

Let $\lambda \mapsto g_\lambda$ be a smooth homotopy between regular endpoints in $\calB(\delta_E)$,
i.e.\  $\lambda \mapsto g_\lambda \in C^\infty([0,1],\calB(\delta_E))$ and $g_0 , g_1 \in \calB_{\text{reg}}(\delta_E)$.
After choosing a suitable reparameterization $t \mapsto \lambda(t)$, we obtain a curve $t \mapsto g_{\lambda(t)}$
which satisfies hypotheses \Cgroup.
Henceforth we shall write $g_t$ instead of $g_{\lambda(t)}$.

Denote by $C_k$ the $k$-th chain group associated with $f$ and the isolating neighbourhood $N^E$.
Define a homomorphism
\[
\psi_k^{1,0} : C_k \to C_k,
\]
by defining it on a generator $X_0$ of $C_k$ to be
\[
\psi_k^{1,0}( X_0 ) = \sum_{ X_1 \in \calS_k(N,f) } i^{1,0}(X_0,X_1) X_1.
\]
Here $i^{1,0}(X_0,X_1)$ counts (modulo $2$) the number of heteroclinic orbits 
$U$ of the nonautonomous equation \eqref{eq:TWE} associated with $(f,g_t,c)$,
with $U(t) \in N$ for all $t\in\R$, and $\lim_{t\to-\infty} U(t) = X_0$, and $\lim_{t\to\infty} U(t) = X_1$.

\begin{mylemma}
\label{lemma:chain_homomorphism}
  The map $\psi_k^{1,0}$ is well-defined and satisfies
  \begin{equation}
  \label{eq:chain_homomorphism}
  \partial_k( N^E , f , g_1 , c ) \circ \psi_k^{1,0} = \psi_{k-1}^{1,0} \circ \partial_k( N^E , f , g_0 , c),
\end{equation}
where $\partial_k(N^E , f , g_i , c) : C_k \to C_{k-1}$ denotes the boundary operator associated with the chosen quadruple $(N^E,f,g_i,c)$.
\end{mylemma}
\begin{myproof}
  By lemma \ref{lemma:energy_bounds_stationary} and since $N$ is finitely generating 
we know that the sum appearing in the definition of $\psi_k^{1,0}$ is finite.
The fact that $i^{1,0}(X_0,X_1)$ itself is well-defined relies on the compactness results we have obtained for the nonautonomous equation,
together with a detailed analysis of the manifold structure of the moduli space.
For this one first has to repeat the analysis from section \ref{sec:moduli} for the nonautonomous equation.
We refer to remark \ref{remark:nonautonomous_moduli} and note that the essential ingredients are compactness, transversality, and glueing for the nonautonomous equation,
together with the fact that $N^E$ is an isolating neighbourhood for the dynamics.

The fact that \eqref{eq:chain_homomorphism} holds follows entirely out of the geometry of the moduli space of index 1 nonautonomous solutions.
See figure \ref{fig:boundaries_nonautonomous} towards the end of section \ref{sec:moduli}, but also \cite{schwarz1993morse}.
The case depicted in figure \ref{fig:traj_n_n0a1}, corresponding solely to the term $\partial_k( N^E , f , g_1 , c ) \circ \psi_k^{1,0}$
in the left hand side of \eqref{eq:chain_homomorphism}, results in no net contribution as we count with $\Z_2$ coefficients.
Likewise, the situation described by figure \ref{fig:traj_n_a1n0} corresponds solely to the right hand side of \eqref{eq:chain_homomorphism} and has no net contribution.
The last case, depicted in figure \ref{fig:traj_n_mixed}, yields the identity \eqref{eq:chain_homomorphism}.
\end{myproof}

Note that lemma \ref{lemma:chain_homomorphism} implies that $\psi_k^{1,0}$ induces a homomorphism $\Psi_k^{1,0}$ on the homology groups.
To verify that $\Psi_k^{1,0}$ is indeed an isomorphism, independence of the chosen homotopy $g_\lambda$ has to be checked.
To do so, one choses a 2-parameter family $g_{\lambda,\mu} \in \calB(\delta_E)$, and considers the parameter dependent moduli space $\calM_\mu(Z_-,Z_+)$,
consisting of pairs
\[
( \mu_* , U_{\mu_*} ) \in \calM_\mu(Z_-,Z_+),
\]
where $U_{\mu_*}$ belongs to the moduli space $\calM(Z_-,Z_+)$ for the nonautonomous equation \eqref{eq:TWE} corresponding to $(f,g_{{\lambda(t)},\mu_*},c)$.
The analysis of this parameter dependent moduli space is much the same as the work in section \ref{sec:moduli},
with the exception that there are more cases to distinguish in the compactification.
We refer to \cite{schwarz1993morse} for a concise description, the conclusion of which is that on the level of homology the homomorphism $\Psi_k^{1,0}$
is independent of the chosen homotopy $g_\lambda$.
Combined with a glueing argument, this yields the relation
\[
\Psi_k^{2,0} = \Psi_k^{2,1} \circ \Psi_k^{1,0},
\]
which in particular shows that $\Psi_k^{1,0}$ is an isomorphism.
We summarize this discussion in the following theorem.
\begin{mythm}
  Given $g_0 , g_1 \in \calB_{\text{reg}}(\delta_E)$, the homomorphism $\Psi_*^{1,0} : \HTW_*(N^E,f,g_0,c) \to \HTW_*(N^E,f,g_1,c)$ is 
  independent of the choice of homotopy $\lambda \mapsto g_\lambda$ between $g_0$ and $g_1$.
  Furthermore,
  \begin{enumerate}
  \item if $g_0 = g_1$, then $\Psi_*^{1,0}$ is the identity, and
  \item for any three $g_0 , g_1 , g_2 \in \calB_{\text{reg}}(\delta_E)$
    one has $\Psi_*^{2,0} = \Psi_*^{2,1} \circ \Psi_*^{1,0}$.
  \end{enumerate}
In particular, $\Psi_*^{1,0}$ is an isomorphism.
\end{mythm}
We will refer to this isomorphism as the canonical isomorphism between $\HTW_*(N^E,f,g_0,c)$ and $\HTW_*(N^E,f,g_1,c)$.

\paragraph{The limit as $E \to \infty$ and $g \to 0$.}

For each $k \in \Z$, define
\[
E_{\text{crit}}(k) := \max\set{ \calE_f(X) }{ X \in \calS_k(N,f) \cup \calS_{k+1}(N,f) }.
\]
The maximum is attained since $N$ is finitely generating.
Pick regular energy levels $E_0 > E_{\text{crit}}(k)$ and $E_1 > E_{\text{crit}}(k)$,
and let $g \in \calB_{\text{reg}}( \min\{ \delta_{E_0} , \delta_{E_1} \} )$.
Let $U$ be a connecting orbit of \eqref{eq:TWE} associated with $(f,g,c)$,
with $\lim_{t\to-\infty} U(t) = X$ where $X$ is a rest point with index $\mu_f(X) = k+1$ or $\mu_f(X) = k$.
Then $\calE_f(U(t)) \leq E_{\text{crit}}(k)$ for all $t \in \R$, hence if $U$ is trapped by $N^{E_0}$ then it is also trapped by $N^{E_1}$,
and vice versa.
Thus $\partial_k(N^{E_0},f,g,c) = \partial_k(N^{E_1},f,g,c)$ and $\partial_{k+1}(N^{E_0},f,g,c) = \partial_{k+1}(N^{E_1},f,g,c)$.
Hence there is a canonical isomorphism
\[
\HTW_k(N^{E_0},f,g,c) \iso \HTW_k(N^{E_1},f,g,c), \qquad \text{for} \quad g \in \calB_{\text{reg}}( \min\{ \delta_{E_0} , \delta_{E_1} \} ).
\]

For any two $\widetilde g_0 , \widetilde g_1 \in \calB_{\text{reg}}(\min\{ \delta_{E_0} , \delta_{E_1} \})$,
the following diagram commutes:
\begin{center}
\begin{tikzcd}
  & \HTW_k(N^{E_0},f,g_0,c) \arrow[ld] \arrow[rd] & \\
  \HTW_k(N^{E_0},f,\widetilde g_0 ,c) \arrow[d] \arrow[rr] & & \HTW_k(N^{E_0},f,\widetilde g_1 ,c) \arrow[d] \\
  \HTW_k(N^{E_1},f,\widetilde g_0 ,c) \arrow[rr] \arrow[rd] & & \HTW_k(N^{E_1},f,\widetilde g_1 ,c) \arrow[ld] \\
  & \HTW_k(N^{E_1},f,g_1,c) &
\end{tikzcd}
\end{center}
Here each of the arrows denote one of the canonical isomorphisms.
It follows that the isomorphism between $\HTW_k(N^{E_0},f,g_0,c)$ and $\HTW_k(N^{E_1},f,g_1,c)$,
which is defined via the commuting diagram, is independent of the intermediate point $\widetilde g_0$.
Denote this isomorphism by $\Phi^{(E_1,g_1) , (E_0,g_1)}_k$.

Thus $\HTW_k(N^E,f,g,c)$ is independent, up to a canonical isomorphism, of $E$ and $g$,
as $E \to \infty$ and $g \to 0$.
We then define $\HTW_k(N,f,c)$ as the isomorphism class of $\HTW_k(N,f,g,c)$, for small generic $g$.
To formalize the notion of defining $\HTW_k(N,f,c)$ up to natural isomorphism,
we make use of an inverse limit of the isomorphisms $\Phi^{(E_1,g_1) , (E_0,g_1)}_k$ and set
\[
\HTW_k(N,f,c) := \varprojlim \HTW_k(N^E,f,g,c).
\]

\subsubsection{Continuation in \texorpdfstring{$f$}{f} and \texorpdfstring{$c$}{c}}

Suppose a curve $t \mapsto (f_t,c_t)$ with regular endpoints $(f_-,c_-)$ and $(f_+,c_+)$ satisfies hypotheses \Cgroup.
In addition, assume $N \subset X^0$ is a finitely generating isolating neighbourhood for the dynamics associated
with the autonomous equations \eqref{eq:TWE} corresponding to $(f_-,0,c_-)$ and $(f_+,0,c_+)$,
as well as being an isolating neighbourhood the nonautonomous equation \eqref{eq:TWE} corresponding to $(f_t,0,c_t)$.
We will call such an $N$ \emph{stable with respect to $t \mapsto (f_t,c_t)$}.
By repeating the construction from the preceding section, 
but with the constant $f$ and $c$ replaced by their $t$-dependent analogues, 
one finds the following theorem.
\begin{mythm}\label{thm:f_isomorphism}
  A curve $t \mapsto (f_t,c_t)$ which satisfies hypotheses \Cgroup,
  which has regular endpoints $(f_-,c_-)$ and $(f_+,c_+)$,
  and for which $N$ is stable,
  induces an isomorphism of homologies:
  \[
  \HTW_*(N,f_-,c_-) \iso \HTW_*(N,f_+,c_+).
  \]
\end{mythm}

\subsubsection{Classes of isomorphic homologies when \texorpdfstring{$N = X^0$}{N=X0}}

Let us now consider the special case where $N = X^0$.
Clearly this means that $N$ is stable with respect to any homotopy between regular endpoints $(f_-,c_-)$ and $(f_+,c_+)$
for which $X^0$ is finitely generating.
Thus $(X^0,f_-,c_-)$ and $(X^0,f_+,c_+)$ will have isomorphic homologies whenever there exists a curve $t \mapsto (f_t,c_t)$
connecting $(f_-,c_-)$ with $(f_+,c_+)$ and satisfying hypotheses \Cgroup.
A first thing to note is that any two $c_- , c_+ > 0$ can be connected via such a path (keeping $f$ fixed).
In fact, the induced isomorphism will then be independent of the chosen homotopy between $c_-$ and $c_+$
(this can be verified by considering two-parameter families of $c > 0$, similar to how independence of the chosen
path $t \mapsto g_t$ is verified).
Thus, for $f$ for which $X^0$ is finitely generating, we can define
\[
\HTW_*(f) := \varprojlim \HTW_*(X^0,f,c),
\]
where $\varprojlim$ is the inverse limit over the isomorphisms induced by homotopies of $c$.
In the remainder of this subsection we give a concrete description of a class of nonlinearities $f$ for which
the homology remains unchanged.

Let $f_*$, which will function as a reference point for our perturbations,
be of class $C^m$ (with $m \geq 4$) and satisfy hypotheses \Fgroup.
Recall that hypotheses \ref{H:superlinear_growth_f} was only needed when considering Neumann or periodic boundary data,
in order to arrive at the compactness results from section \ref{sec:compactness}.
In contrast to this, even when considering Dirichlet boundary data
we will now demand that $f_*$ satisfies the superlinear growth condition \ref{H:superlinear_growth_f}.
This will help in the construction of allowed perturbations from $f_*$, as we will see shortly hereafter.
Consider $f$ of the form
\begin{equation}
  \label{eq:perturbed_f}
  f(x,u) = \alpha(x) f_*(x,u) + h(x,u),
\end{equation}
where
$ \alpha \in C_b^m(\Omega)$, and $h \in C^m(\Omega\times\R)$,
\[
\inf_{x\in\Omega} \alpha(x) > 0,
\]
and
\[
\limsup_{|u| \to \infty} \sup_{x\in\Omega} \left| \frac{ h(x,u) }{ f_*(x,u) } \right| = 0.
\]

\begin{mylemma}
\label{lemma:Sigma_hypotheses}
Any $f$ of the form \eqref{eq:perturbed_f} satisfies hypotheses \Fgroup.
\end{mylemma}
\begin{myproof}
  It is obvious that $f$ will satisfy hypotheses \ref{H:growth_bound_f} and \ref{H:superlinear_growth_f}.
  Checking whether \ref{H:contact_condition_odd} or \ref{H:contact_condition_even} are satisfied takes slightly more effort.
  Let us assume that $f_*$ satisfies \ref{H:contact_condition_odd}, the argument for the other case is completely similar.
  We shall prove that $f$ then also satisfies \ref{H:contact_condition_odd},
  i.e., we need to show that
  \begin{equation}
    \label{eq:perturbed_contact_condition}
    | \alpha(x) F_*(x,u) + H(x,u) | \leq C + \frac{ \theta }{ 2 } ( \alpha(x) f_*(x,u) + h(x,u) ) u,
  \end{equation}
  for some constants $C \geq 0$, and $-1 < \theta < 1$.
  Here $H(x,u) := \int_0^u h(x,u) \d s$.
  First note that, by dividing by $\alpha(x)$ and updating the values of $C$ and $h(x,u)$,
  it suffices to prove \eqref{eq:perturbed_contact_condition} for $\alpha \equiv 1$.

Suppose for the moment that
  \begin{equation}
    \label{eq:47}
    \limsup_{|u|\to\infty} \sup_{x\in\Omega} \left| \frac{ H(x,u) }{ F_*(x,u) } \right| = 0.
  \end{equation}
  Then for each $\epsilon > 0$ we can find $C_\epsilon > 0$ such that
  \[
  | F_*(x,u) + H(x,u) | \leq C_\epsilon + (1 + \epsilon) \frac{\theta_*}{2} f_*(x,u) u,
  \]
  where $\theta_*$ is the value of $\theta$ for which $f_*$ satisfies hypothesis \ref{H:contact_condition_odd}.
  Then note that hypotheses \ref{H:contact_condition_odd} and \ref{H:superlinear_growth_f} taken together
  imply that $\theta_* f_*(x,u) u$ is strictly positive for $|u|$ large.
  Combining this observation with the assumption that $h = o(f_*)$ as $|u| \to \infty$, uniformly in $x$, gives
  \[
  \frac{ \theta_* }{ 2 } f_*(x,u) u \leq \widetilde C_\epsilon + (1+\epsilon) \frac{ \theta_* }{ 2 } ( f_*(x,u) + h(x,u) ) u.
  \]
  Combined with the penultimate estimate this shows that $f$ satisfies hypothesis \ref{H:contact_condition_odd},
  with $\theta$ arbitrarily close to $\theta_*$.
  
  We still need to see why \eqref{eq:47} holds, i.e. that $H = o(f_*)$ as $|u| \to \infty$, uniformly in $x$.
As we have already seen that hypotheses \ref{H:contact_condition_odd} and \ref{H:superlinear_growth_f} together imply that
$f_*(x,u)$ does not change signs for $|u|$ large.
Therefore, to prove \eqref{eq:47} we can replace $f_*$ by $|f_*|$, and can thus assume that $f_*(x,u) > 0$ for all $x$ and $u$.
  Since $h = o(f_*)$ as $|u|\to\infty$, uniformly in $x$,
  for each $\epsilon > 0$ we can find $K \geq 0$ such that for all $|u| \geq K$ we have
  \[
      |H(x,u)| \leq \epsilon F_*(x,u) + \frac{1}{2} \int_{-K}^K |h(x,s)| \d s - \frac{\epsilon}{2} \int_{-K}^K f_*(x,s) \d s.
  \]
  Hypothesis \ref{H:superlinear_growth_f} ensures that $\inf_{x\in\Omega} F_*(x,u) \to \infty$ as $|u| \to \infty$.
  Hence given any $\epsilon > 0$ and $K \geq 0$ we can find $L \geq K$ such that for any $x \in \Omega$ and $|u| \geq L$ one has
  \[
  \frac{1}{2} \int_{-K}^K |h(x,s)| \d s - \frac{\epsilon}{2} \int_{-K}^K f_*(x,s) \d s \leq \epsilon f_*(x,u).
  \]
  Hence $|H(x,u)| \leq 2 \epsilon F_*(x,u)$ for all $|u| \geq L$.
  Since $\epsilon > 0$ was chosen arbitrarily, we arrive at \eqref{eq:47}.
\end{myproof}

We will call a curve $t \mapsto f_t \in C^m(\R,C^m(\Omega\times\R))$,
where $f_t$ is of the form 
\[
f_t(x,u) = \alpha_t(x) f(x,u) + h_t(x,u),
\]
an \emph{$\epsilon$-perturbation of $f_*$} provided that 
\begin{enumerate}
\item it is constant for $t \in (-\infty,1)$ and for $t \in (1,\infty)$,
\item $\sup_{(t,x) \in \R\times\Omega} | 1 - \alpha_t(x) | < \epsilon$ and $\sup_{(t,x) \in\R\times\Omega} | \partial_t \alpha_t(x) | < \epsilon$, and
\item $|h_t(x,u)| \leq \epsilon ( 1 + |f_*(x,u)| )$ for all $t \in \R$, $(x,u) \in \Omega\times\R$.
\end{enumerate}

\begin{mylemma}
\label{lemma:Sigma_hypotheses_nonautonomous}
  There exists a sufficiently small $\epsilon > 0$ such that any $\epsilon$-perturbation of $f_*$ satisfies hypotheses \Cgroup.
\end{mylemma}
\begin{myproof}
Properties (2) and (3) from the definition of $\epsilon$-perturbations ensure that 
the estimates made in the proof of lemma \ref{lemma:Sigma_hypotheses} can be made uniformly in $t$.
Hence hypothesis \ref{C:pointwise_H} is satisfied.
By definition (according to property (1)) $\epsilon$-perturbations also satisfy \ref{C:asymptotic_constant}.

Left to check is that hypothesis \ref{C:growth_bound_f} holds when $\epsilon$ is chosen sufficiently small, i.e.\ that
\[
| \partial_t \alpha_t(x) F(x,u) + \partial_t H_t(x,u) | \leq C + \Theta | \alpha_t(x) F(x,u) + H_t(x,u) |,
\]
for some $C \geq 0$ and sufficiently small $\Theta$.
Here $H_t(x,u) = \int_0^u h_t(x,s) \d s$.
As noted in the proof of lemma \ref{lemma:Sigma_hypotheses}, for each $t$ the function $H_t$ is $o(F)$ as $|u| \to \infty$,
uniformly in $x$, and in light of property (3) from the definition of $\epsilon$-perturbations,
these estimates are also uniform in $t$.
Thus it suffices to see that
\[
| \partial_t \alpha_t(x) F(x,u) + \partial_t H_t(x,u) | \leq C + \Theta | \alpha_t(x) F(x,u) |.
\]
Since by assumption $| \partial_t \alpha_t(x) | < \epsilon$ and $| 1 - \alpha_t(x) | < \epsilon$,
given any $\Theta > 0$ we can find $\epsilon > 0$ such that
\[
| \partial_t \alpha_t(x) F(x,u) | \leq \frac{\Theta}{2} | \alpha_t(x) F(x,u) |.
\]
Furthermore, since for each $t$ we have $H_t = o(F)$ as $|u| \to \infty$, uniformly in $x$,
we also have $\partial_t H_t = o(F)$ as $|u| \to \infty$, uniformly in $x$.
Hence in particular, for any given $\Theta > 0$ there exists a $C \geq 0$ such that
\[
| \partial_t H_t(x,u) | \leq C + \frac{\Theta}{2} | \alpha_t(x) F(x,u) |.
\]
This proves that hypothesis \ref{C:growth_bound_f} holds, and the constant $\Theta$ can be made arbitrary small
by choosing $\epsilon$ sufficiently small.
\end{myproof}

Combining lemmata \ref{lemma:Sigma_hypotheses} and \ref{lemma:Sigma_hypotheses_nonautonomous}
with theorem \ref{thm:f_isomorphism} shows that any two nonlinearities (regular and for which $X^0$ is finitely generating)
which can be connected via an $\epsilon$-perturbation have isomorphic homologies.

Denote by $\Sigma(f_*)$ the set of all nonlinearities $f$ of the form \eqref{eq:perturbed_f}.
Endow $\Sigma(f_*)$ with the topology of $C^m_{\text{loc}}(\Omega\times\R)$.
Let $\Sigma_{\text{reg}}(f_*)$ consist of all those $f \in \Sigma(f_*)$ which are regular;
in light of theorem \ref{thm:genericity_hyperbolicity} $\Sigma_{\text{reg}}(f_*)$ is dense in $\Sigma(f_*)$.
Finally, denote by $\Sigma_{\text{fin}}(f_*)$ the set of those $f \in \Sigma_{\text{reg}}(f_*)$ for which $X^0$ is finitely generating.
Now note that any two nonlinearities $f_0 , f_1 \in \Sigma(f_*)$ 
can be connected via a concatenation of finitely many $\epsilon$-perturbations.
However, these $\epsilon$-perturbations can only be inducing isomorphisms of homologies
if the endpoints of each of the individual $\epsilon$-perturbations can be chosen to be elements of $\Sigma_{\text{fin}}(f_*)$.
Hence we arrive at the following theorem.
\begin{mythm}\label{thm:f_isomorphism_class}
  Fix arbitrary $f_0 , f_1 \in \Sigma_{\text{fin}}(f_*)$ and suppose $f_0$ and $f_1$ belong to the 
same path-component of $\clin{\Sigma_{\text{fin}}(f_*)}{\Sigma(f_*)}$.
Then
\[
\HTW_*(f_0) \iso \HTW_*(f_1).
\]
\end{mythm}

\begin{myremark}
  In the examples in section \ref{sec:application} we consider nonlinearities $f_*$ which are homogeneous in $u$.
For these nonlinearities it follows that $\Sigma_{\text{fin}}(f_*) = \Sigma_{\text{reg}}(f_*)$,
and since $\clin{\Sigma_{\text{fin}}(f_*)}{\Sigma(f_*)} = \Sigma(f_*)$ is path-connected,
it follows that any two $f_0 , f_1 \in \Sigma_{\text{reg}}(f_*)$ have isomorphic homologies.
It remains an open question whether any regular $f$ which satisfies hypotheses \Fgroup 
automatically has $X^0$ as a finitely generating isolating neighbourhood.
If this turns out to be true, it follows that any two regular nonlinearities $f_0 , f_1 \in \Sigma_{\text{reg}}(f_*)$
induce isomorphic homologies.
One can then proceed to define the homology for any nongeneric nonlinearity $f_*$ as the inverse limit over $\Sigma_{\text{reg}}(f_*)$.
\end{myremark}

\subsection{Direct sum property}
We conclude this section with an algebraic property of the travelling wave homology 
which will come in useful when applying the theory to concrete problems.

\begin{mylemma}
\label{lemma:direct_sum}
Let $f$ be a regular nonlinearity and $c > 0$.
Let $N \subset X^0$ be a finitely generating isolating neighbourhood for the dynamics associated with $(f,0,c)$.
Suppose $V_1, V_2 \subset X^0$ are isolating neighbourhoods for the dynamics associated with $(f,0,c)$, such that
$\BInv(N;f,0,c) = \BInv(V_1;f,0,c) \cup \BInv(V_2;f,0,c)$ and $\BInv(V_1;f,0,c) \cap \BInv(V_2;f,0,c) = \emptyset$.
Then
\[
  \HTW_*(N,f,c) = \HTW_*(V_1,f,c) \oplus \HTW_*(V_2,f,c).
\]
\end{mylemma}
\begin{myproof}
Note that without loss of generality we may assume $V_1 \cap V_2 = \emptyset$.
Fix any $k \in \Z$ and $E > E_{\text{crit}}(k)$.
We claim that, for any $g \in \calB_{\text{reg}}(\delta)$ with $0 < \delta \leq \delta_E$ sufficiently small,
any bounded solution $U$ of the dynamics of \eqref{eq:TWE} associated with $(f,g,c)$,
and for which $U(t) \in N^E$ for all $t \in \R$,
has either $U(t) \in \inter{ V_1^E }$ for all $t \in \R$, or $U(t) \in \inter{ V_2^E }$ for all $t \in \R$.
Suppose this is not the case.
Then one can find a sequence of $(g_n)_n$ with $g_n \to 0$ as $n\to\infty$,
bounded solutions $U_n$ of \eqref{eq:TWE} associated with $(f,g_n,c)$ and such that $U_n(t) \in N^E$ for all $t \in \R$,
and a sequence $(t_n)_n \subset \R$ such that $U_n(t_n) \not\in \inter{ V_1^E } \cup \inter{ V_2^E }$ for all $n \in \N$.
As discussed in remark \ref{remark:Palais-Smale_type_compactness},
 $U_n(\cdot + t_n)$ converges over a subsequence to a solution $U_\infty$ of \eqref{eq:TWE} corresponding to $(f,0,c)$.
But then $U_\infty(0) \not\in \inter{ V_1^E } \cup \inter{ V_2^E }$,
and since $\calE_f(U_\infty(0)) \leq E$, also $U_\infty(0) \not\in \inter{ V_1 } \cup \inter{ V_2 }$.
However, $U_\infty(t) \in N$ for all $t \in \R$.
Hence we have constructed a solution $U_\infty$ of \eqref{eq:TWE} corresponding to $(f,0,c)$ which is isolated by $N$, but not isolated by either $V_1$ or $V_2$.
This is in contradiction with the hypotheses of the lemma.

We conclude that whenever $g \in \calB_{\text{reg}}(\delta)$ and $0 < \delta \leq \delta_E$ is sufficiently small,
the sets $N^E$, $V_1^E$, and $V_2^E$ are isolating neighbourhoods for the dynamics of \eqref{eq:TWE} associated with $(f,g,c)$,
and
\[
\BInv(N^E;f,g,c) = \BInv(N^E;f,g,c) \cup \BInv(N^E;f,g,c).
\]
Hence the critical groups satisfy the direct sum property $C_*(N^E) = C_*(V_1^E) \oplus C_*(V_2^E)$
(here $C_*(N)$ denotes the chain group corresponding to a given isolating neighbourhood $N$)
and the boundary operator $\partial_k(N^E,f,g,c)$ factorizes through this direct sum.
Hence
\[
\HTW_k(N^E , f,g,c) = \HTW_k(V_1^E, f,g,c) \oplus \HTW_k(V_2^E, f,g,c).
\]
This is true for any $E$ sufficiently large and $g$ sufficiently small, hence the conclusion of the lemma follows.
\end{myproof}

 \section{Applications}
\label{sec:application}

In this section we will first compute the travelling wave homology for various classes of nonlinearities,
and finally give some examples of how this information can be used to derive conclusions about
the existence of travelling waves.

In this section we consider nonlinearities $f_{\text{odd},\pm}$ and $f_{\text{even},\pm}$ 
as introduced in \eqref{eq:f_odd} and \eqref{eq:f_even}.
It was already pointed out in remark \ref{remark:practical_restrictions_f} that these nonlinearities satisfy hypotheses \Fgroup.
We will compute the travelling wave homologies for these nonlinearities of this form,
after which we will show how this information can be used to prove existence of travelling waves in reaction-diffusion equations.

\begin{mythm}
\label{thm:computation_homology}
  For any regular nonlinearity $f = f_{\text{odd},\pm}$ or $f = f_{\text{even},\pm}$, the set $N = X^0$ is a finitely generating isolating neighbourhood.
  There exists a $k_0 \in \Z$ (depending on the chosen normalized Morse index $\mu$) such that
  \[
  \HTW_k( f_{\text{odd},-} ) \iso
  \left\{
    \begin{array}{l l}
      \Z_2 & \text{if } k = k_0, \\
      0 & \text{otherwise}.
    \end{array}
  \right.
  \]
  Furthermore,
  \begin{align*}
    \HTW_*( f_{\text{odd},+} ) &= 0, \\
    \HTW_*( f_{\text{even},-} ) &= 0, \\
    \HTW_*( f_{\text{even},+} ) &= 0.
  \end{align*}
\end{mythm}
\begin{myproof}
We begin with verifying that $X^0$ is finitely generating for each of the nonlinearities, so that the homologies are indeed well-defined.
We will in fact show that $\Sigma_{\text{fin}}(f) = \Sigma_{\text{reg}}(f)$.
This also shows that, in light of theorem \ref{thm:f_isomorphism_class}, 
all the nonlinearities in $\Sigma_{\text{reg}}(f)$ have isomorphic homologies.

When either $f = f_{\text{odd},-}$, or $f = f_{\text{even},-}$, or $f = f_{\text{even},+}$,
this is true because the set of solutions $z$ of
\[
\left\{
\begin{array}{l l}
\Delta z + f(\cdot,z) = 0 & \text{on } \Omega, \\
B(z) = 0 & \text{on } \bdy\Omega
\end{array}
\right.
\]
is compact in $H_B^2(\Omega)$ (see e.g.\  \cite{fiedler1998large}),
hence $\calS(f)$ is finite for regular $f$.
Hence $\Sigma_{\text{fin}}(f) = \Sigma_{\text{reg}}(f)$.

When $f = f_{\text{odd},+}$ the set $\calS(f)$ will typically not be finite.
Assume that $f$ is regular,
and recall from the definition of the normalized Morse index 
that there exists a constant $m_0$ such that for the index of a rest point $Z=(z,0) \in \calS(f)$ we have the following identity:
\begin{equation}
  \label{eq:59}
  \mu_{f}(Z) = m_0 - m_f(z),
\end{equation}
where $m_f$ is the classical Morse index
\[
m_f(z) := \# \big( \sigma(\Delta + f_u(x,z) ) \cap (0,\infty) \big).
\]
For any given $k \in \Z$, a classical result from Bahri and Lions (see \cite{bahri1992solutions}, 
but also \cite{yang1998nodal, ramos2009bahri, yu2014solutions, harrabi2012priori, hajlaoui2015morse}) 
then gives us a priori bounds on the $L^\infty$ norm of rest points $Z$ with a given morse index $m_f(z) \leq k$.
In light of \eqref{eq:59} this gives $L^\infty$ bounds on $Z$ with a given index $\mu_{f}(Z) \geq k$.
Thus $\calS_k(X^0,f)$ is finite for each $k$, i.e.\  $X^0$ is finitely generating for $f$.
Thus we again find that $\Sigma_{\text{fin}}(f) = \Sigma_{\text{reg}}(f)$.

We now proceed to the actual computation of the various homologies.
Computation of the homology for $f_{\text{odd},-}$ requires a different technique from the computation of the homology for $f_{\text{odd},+}$,
but these approaches work for any boundary condition.
The homology for the nonlinearities $f_{\text{even},-}$ and $f_{\text{even},+}$ can be computed using the same technique,
but in this case the chosen approach depends on the boundary data.

\paragraph{Computation for $f_{\text{odd},-}$.}

By theorem \ref{thm:f_isomorphism_class},
\[
\HTW_*( f_{\text{odd},-} ) \iso \HTW_*( f_\epsilon ),
\]
where
\[
f_\epsilon(x,u) = -|u|^{p-1} u - \epsilon u.
\]
For $\epsilon \geq 0$, suppose $z$ is a solution of
\begin{equation}
  \label{eq:60}
  \left\{
  \begin{array}{l l}
    \Delta z + f_\epsilon(x,z) = 0 & \text{on } \Omega, \\
    B(z) = 0 & \text{on } \bdy \Omega.
  \end{array}
\right.
\end{equation}
Then
\[
\int_\Omega | \nabla z |^2 + |z|^{p+1} + \epsilon |z|^2 \d x = - \int_\Omega \big( \Delta z + f_\epsilon(x,z) \big) z \d x = 0.
\]
Hence the only solution $z$ of \eqref{eq:60} is $z \equiv 0$.
Moreover, for $\epsilon > 0$ sufficiently small this rest point is hyperbolic, so that
\[
\HTW_k( f_\epsilon ) \iso
  \left\{
    \begin{array}{l l}
      \Z_2 & \text{if } k = \mu_{f_\epsilon}(0), \\
      0 & \text{otherwise}.
    \end{array}
  \right.
\]

\paragraph{Computation for $f_{\text{odd},+}$.}

By theorem \ref{thm:f_isomorphism_class},
\[
\HTW_*(f_{\text{odd},+}) \iso \HTW_*(f_\beta),
\]
where
\[
f_\beta(x,u) = |u|^{p-1} u + \beta u + \phi_\beta(x,u).
\]
Here $\phi_\beta$ is chosen such that any solution $z \in H_B^2(\Omega)$ of $\Delta z + f_\beta(\cdot,z) = 0$ is hyperbolic,
and $\| \partial_u \phi_\beta \|_{L^\infty(\Omega\times\R)} \leq 1$.
Such $\phi_\beta$ exist in light of theorem \ref{thm:genericity_hyperbolicity}.
Note that if $z$ is a solution of $\Delta z + f_{\beta}(\cdot,z) = 0$, then for the linearization one has
\[
\Delta + \ddi{f_{\beta}}{u}(\cdot,z) = \Delta + p |z|^{p-1} + \beta + \partial_u \phi_\beta(x,z) \succeq \Delta + \beta - 1,
\]
where ``$\succeq$'' denotes the partial ordering on $L^2(\Omega)$ induced by the cone of positive operators on $L^2(\Omega)$.
It then follows from the min-max characterization of eigenvalues of self-adjoint operators (see e.g.\  \cite{evans1998partial}) that
\[
m_{f_\beta}(z) \geq \# \big( \sigma( \Delta + \beta - 1 ) \cap (0,\infty) \big).
\]
So for any $Z=(z,0) \in \calS(f_{\beta})$ it follows that
\[
\mu_{f_{\beta}}(Z) = m_0 - m_{f_\beta}(z) \leq m_0 - \# \big( \sigma( \Delta + \beta - 1 ) \cap (0,\infty) \big).
\]
The right hand side tends to $-\infty$ as $\beta \to \infty$.
Hence, for each given $k$ we can choose $\beta$ sufficiently large so that $\calS_k(X^0,f_{\beta}) = \emptyset$,
 and therefore $\HTW_k(f_{\beta}) = 0$.

\paragraph{Computation for $f_{\text{even},\pm}$, Dirichlet boundary data.}

By theorem \ref{thm:f_isomorphism_class},
\[
\HTW_*(f_{\text{even},\pm}) \iso \HTW_*(f_{\mu,\pm}),
\]
where
\[
f_{\mu,\pm}(x,u) = \pm \mu |u \pm 1|^p.
\]
We will argue that for $\mu > 0$ sufficiently large, there is no solution to the equation $\Delta z + f_{\mu,\pm}(\cdot,z) = 0$.
Let us now first discuss the case where $f_{\mu,\pm} = f_{\mu,-}$.
Suppose $z \in H_B^2(\Omega)$ is a solution of
\begin{equation}
\label{eq:test_eq_even}
\Delta z + f_{\mu,-}(x,z) = 0.
\end{equation}
Then $z$ is subharmonic, and zero on $\bdy \Omega$, hence by the maximum principle $z \leq 0$.
So
\[
f_{\mu,-}(x,z) = \big( - |z-1|^{p-2}(z-1) \big) \big( \mu (z-1) \big) \leq \mu (z-1) \leq \mu z.
\]
Furthermore, $z \equiv 0$ is clearly not a solution of \eqref{eq:test_eq_even}.
Thus $z$ also satisfies
\begin{equation}
\label{eq:test_eq_even_simplified}
\left\{
\begin{array}{l l}
\Delta z + \mu z \geq 0 & \text{on } \Omega, \\
z \neq 0 & \text{on } \Omega, \\
z \leq 0 & \text{on } \Omega, \\
z = 0 & \text{on } \bdy \Omega.
\end{array}
\right.
\end{equation}
Let $\lambda_1$ be the fundamental eigenvalue of $\Delta$ with Dirichlet boundary data, and let $\phi_1$ be a corresponding eigenfunction.
Recall (see e.g.\  \cite{evans1998partial}) that $\lambda_1 < 0$ and we may assume that $\phi_1(x) > 0$ for all $x \in \inter \Omega$.
Now multiply \eqref{eq:test_eq_even_simplified} by $\phi_1$ and integrate to obtain
\begin{equation}
\label{eq:eigv_estimate}
0 \leq \int_\Omega \phi_1 \Delta z + \mu \phi_1 z \d x = \int_\Omega z \Delta \phi_1 + \mu \phi_1 z \d x = (\lambda_1 + \mu) \int_\Omega \phi_1 z \d x.
\end{equation}
Since $z \leq 0$ and $z \neq 0$ on $\Omega$, and $\phi_1 > 0$ on $\inter \Omega$, the last integral must be strictly negative.
But then $(\lambda_1 + \mu) \int_\Omega \phi_1 z \d x < 0$ for $\mu > -\lambda_1$, contradicting inequality \eqref{eq:eigv_estimate}.
Hence there can be no solution of \eqref{eq:test_eq_even} whenever $\mu > -\lambda_1$.
Similarly, for $\mu > -\lambda_1$ there are no solutions of $\Delta z + f_{\mu,+}(\cdot,z) = 0$ with Dirichlet boundary conditions.
Hence for $\mu > -\lambda_1$ we have $\HTW_*(f_{\mu,\pm}) = 0$.

\paragraph{Computation for $f_{\text{even},\pm}$, Neumann or periodic boundary data.}

By theorem \ref{thm:f_isomorphism_class},
\[
\HTW_*(f_{\text{even},\pm}) \iso \HTW_*(f_\pm),
\]
where
\[
f_\pm(x,u) = \pm |u|^p \pm 1.
\]
Now if $z \in H^2_B(\Omega)$ were a solution of $\Delta z + f_\pm(z) = 0$, one would find that
\[
\Vol(\Omega) \leq \pm \int_\Omega f_\pm(z) \d x = \mp \int_\Omega \Delta z \d x.
\]
But by Stokes' theorem and the chosen boundary data, the last integral equals $0$.
Hence $\Delta z + f_\pm(z) = 0$ does not have any solutions with Neumann or periodic boundary data.
Therefore, $\HTW_*(f_\pm) = 0$.
\end{myproof}

With the homologies computed, we can now apply this information to prove existence of travelling wave solutions of \eqref{eq:reaction_diffusion}.
\begin{mythm}
\label{thm:forced_existence_of_waves}

Consider any wave speed $c \neq 0$, and let $k \geq 1$.
Then the following holds:
\begin{itemize}

\item
If $f = f_{\text{odd},-}$ and \eqref{eq:TWE} has at least $2k$ distinct hyperbolic stationary solutions,
then \eqref{eq:reaction_diffusion} has at least $k$ distinct travelling wave solutions of wave speed $c$.
More precisely, to each given hyperbolic stationary solution $Z$ (but with the possible exception of at most one of them),
there corresponds at least one travelling wave $U$
such that $\alpha(U) = \{Z\}$ or $\omega(U) = \{Z\}$ (but it is possible that $\omega(U)$ resp.\ $\alpha(U)$ consist of non-hyperbolic stationary solutions).

\item
If either $f = f_{\text{odd},+}$, or $f = f_{\text{even},-}$, or $f = f_{\text{even},+}$,
and \eqref{eq:TWE} has at least $2k-1$ distinct hyperbolic stationary solutions,
then \eqref{eq:reaction_diffusion} has at least $k$ distinct travelling wave solutions of wave speed $c$.
More precisely, to each given hyperbolic stationary solution $Z$,
there corresponds at least one travelling wave $U$
such that $\alpha(U) = \{Z\}$ or $\omega(U) = \{Z\}$ (but it is possible that $\omega(U)$ resp.\ $\alpha(U)$ consist of non-hyperbolic stationary solutions).

\end{itemize}
Furthermore, in each of these cases there exists at least one more stationary solution (which might be non-hyperbolic).
\end{mythm}
\begin{myproof}
Let us first discuss the case where $f = f_{\text{odd},-}$.
Fix any $c > 0$.
Let $S_1$ consist of the $2k$ given hyperbolic stationary solutions of \eqref{eq:TWE}.

Suppose there exist two points $Z_1, Z_2 \in S_1$,
such that for both of them there exist no connecting orbit which converges to $Z_i$ in either forward or backward time.
This means that both $\{Z_1\}$ and $\{Z_2\}$ are connected components of $\BInv(X^0;f,0,c)$.
Hence we can find mutually disjoint isolating neighbourhoods $V_1$, $V_2$, and $N$ such that
$\{Z_i\} = \BInv( V_i ; f,0,c )$, and $\BInv(X^0;f,0,c) = \BInv(V_1;f,0,c) \cup \BInv(V_2;f,0,c) \cup \BInv(N ; f,0,c)$.
We claim we can choose a small perturbation $f_\epsilon$ of $f$ such that
\begin{enumerate}
\item $f_\epsilon$ is regular,
\item $V_1$, $V_2$, and $N$ are isolating neighbourhoods for \eqref{eq:TWE} associated with $(f_\epsilon,0,c)$, and
\item $\BInv(X^0;f_\epsilon,0,c) = \BInv(V_1;f_\epsilon,0,c) \cup \BInv(V_2;f_\epsilon,0,c) \cup \BInv(N ; f,0,c)$.
\end{enumerate}
Indeed, set $f_\epsilon := f + \phi_\epsilon$, where $\phi_\epsilon$ is an arbitrarily chosen $\phi_\epsilon \in \calF^m_{\text{reg}}$ 
with $\| \phi_\epsilon \|_{\calF^m} \leq \epsilon$
(recall that $\calF^m_{\text{reg}}$ and $\|\cdot\|_{\calF^m}$ are defined in section \ref{subsec:generic_hyperbolicity}).
The first property then follows from theorem \ref{thm:genericity_hyperbolicity}.
The fact that the other two properties hold for sufficiently small choices of $\epsilon$ follows from an argument identical to the one used in the proof of
lemma \ref{lemma:direct_sum}.

In light of theorem \ref{thm:computation_homology}, the isolating neighbourhood $X^0$ is finitely generating for $f_\epsilon$
(hence so are the isolating neighbourhoods $N$, $V_1$, and $V_2$),
and $\HTW_*(X^0,f_\epsilon,c)$ is of rank 1.
Since $V_1$ and $V_2$ each contain exactly one hyperbolic stationary solution of the unperturbed equation \eqref{eq:TWE} associated with $(f,0,c)$, 
it follows from the implicit function theorem that
(after choosing a sufficiently small perturbation and shrinking the neighbourhoods $V_1$ and $V_2$)
they each contain exactly one hyperbolic stationary solution of the perturbed equation \eqref{eq:TWE} associated with $(f_\epsilon,0,c)$.
Hence both $\HTW_*(V_1,f_\epsilon,c)$ and $\HTW_*(V_2,f_\epsilon,c)$ are of rank 1.

\begin{figure}[t]
  \centering
  \def\svgwidth{0.7\textwidth}
  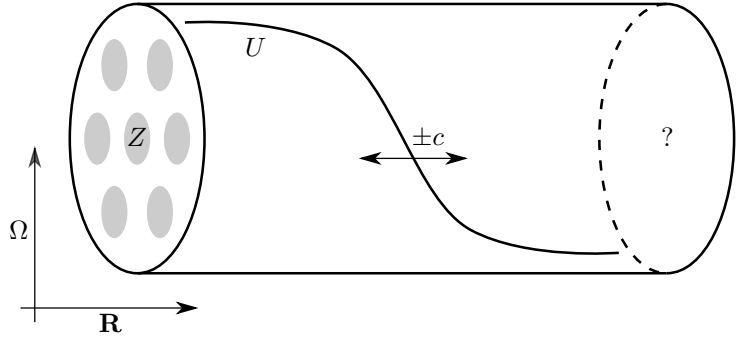
  \caption{Detection of a travelling wave with given wave speed $c$.
  When the index of the stationary solution $Z$ does not coincide with the full homology,
  a heteroclinic orbit $U$ connecting $Z$ to another, unknown, stationary solution must exist.
  Direction of propagation of the travelling wave depends on whether $Z$ is the $\alpha$- or $\omega$-limit set of $U$.
}
  \label{fig:detected_travelling_wave}
\end{figure}

By invariance and the direct sum property (lemma \ref{lemma:direct_sum}) of the homology,
\[
  \HTW_*(X^0,f_\epsilon,c)
  = \HTW_*(V_1,f_\epsilon,c) \oplus \HTW_*(V_2,f_\epsilon,c) \oplus \HTW_*(N,f_\epsilon,c).  
\]
We have arrived at a contradiction, since the homology on the left hand side has rank $1$, while the homology on the right hand side has rank at least $2$.
From this we conclude that one of the $Z_i$ must have a corresponding bounded solution $U$ of \eqref{eq:TWE} associated to the triplet $(f,0,c)$,
such that either $\alpha(U) = \{Z_i\}$ or $\omega(U) = \{Z_i\}$.
See also figure \ref{fig:detected_travelling_wave}.

For simplicity of the argument, say $Z_1$ is the point which is not isolated, 
and $U$ converges towards it in backward time, i.e., $\alpha(U) = \{Z_1\}$.
From lemma \ref{lemma:gradient_like} it follows that $\omega(U)$ consists of stationary solutions of \eqref{eq:TWE},
which can either be another one of the hyperbolic stationary solutions, or (a family of) non-hyperbolic solution(s).
In the first case (that is, $\{ Z_3 \} = \omega(U)$ is hyperbolic), set $S_2 := S_1 \setmin \{ Z_1 , Z_3 \}$.
In the latter case (that is, $\omega(U)$ consists of non-hyperbolic solutions), set $S_2 := S_1 \setmin \{ Z_1 \}$.
By repeating the preceding argument with $Z_1$ and $Z_2$ replaced by points $\widetilde Z_1 , \widetilde Z_2 \in S_2$,
we prove the existence of another connecting orbit which is distinct from the one previously found.
We can iterate this procedure $k$ times, at which point the iteration terminates since we can no longer guarantee that $\# S_k \geq 2$.

In the other cases, when either $f = f_{\text{odd},+}$, or $f = f_{\text{even},-}$, or $f = f_{\text{even},+}$,
a similar argument shows that the existence of a single isolated rest point $Z$ is already excluded.
This again relies on the direct sum property, combined with the fact that for these nonlinearities the full homology $\HTW_*(X^0,f ,c)$
is of rank $0$.
\end{myproof}


\appendix
\section{Fredholm theory}
\label{sec:appendix_fredholm}

Let $\frakL$ and $\frakL_{\text{hyp}}$ be as defined in section \ref{sec:fredholm}.
In this appendix we will fill in some details about the Fredholm theory for the operator
\begin{align*}
    &\calD_L : W^{1,2}(\R;X^0,X^1) \to L^2(\R;X^0), \\
    &\calD_L W = \partial_t W + L(t) W,
\end{align*}
where $L \in C^0(\R;\frakL)$ is such that the limits $L_\pm = \lim_{t\to\pm\infty} L(t)$ exist in the uniform operator topology on $\calL(X^1,X^0)$,
and $L_\pm \in \frakL_{\text{hyp}}$.

\subsection{Fredholm alternative for \texorpdfstring{$\calD_L$}{DL}}

We will use the results from \cite{rabier2004robbin}.
Let us first cite a simplified version of the main result from said article.
\begin{mythm}[\cite{rabier2004robbin}]
\label{thm:rabier}
  Let $H$ be a Hilbert space and $W \subset H$ a normed space.
Let $( L(t) )_{t\in\R}$ be a family of unbounded operators on $H$ with common domain $W$.
Assume that the following holds:
\begin{enumerate}
\item
$W$ is a Banach space and the embedding $W \hookrightarrow H$ is continuous, compact, and dense.

\item
$L \in C^0(\R, \calL(W,H))$.

\item
There are operators $L_-, L_+ \in GL(W,H)$ such that
\[
\lim_{t \to \infty} \| L(t) - L_+ \|_{\calL(W,H)} = \lim_{t\to -\infty} \| L(t) - L_- \|_{\calL(W,H)} = 0.
\]

\item
For every $t \in \R \cup \{\pm\infty\}$ there exist constants $C_0(t) > 0$ and $R_0(t) > 0$ such that
\[
 \| \lambda ( L(t) - i \lambda )^{-1} \|_{\calL(H)} \leq C_0(t) \qquad \text{for all} \quad \lambda \in \R,\; |\lambda| \geq R_0(t).
\]

\item
$\sigma(L_\pm) \cap i \R = \emptyset$.
\end{enumerate}
Then the operator $\partial_t + L(t)$ is Fredholm from $W^{1,p}(\R;H,W)$ to $L^p(\R;H)$ for every $p \in (1,\infty)$. 
\end{mythm}

Note that, in contrast to the classical Robbin-Salamon theorem \cite{robbin1995spectral}, the operators $L(t)$ do not have to be symmetric.
In fact, the spectrum may cross the imaginary axis, as long as we remain in control of the resolvent as per property (4).

In our case, $W = X^1$ and $H = X^0$.
It is then obvious that properties (1), (2), (3), and (5) hold.
The fact that also property (4) is satisfied is the content of lemma \ref{lemma:resolvent_decay}.

Theorem \ref{thm:rabier} combined with lemma \ref{lemma:resolvent_decay} shows that $\calD_L$ is a Fredholm operator.
The index is independent of the particular choice of the path $L$, but only depends on the hyperbolic limits $L_-$, $L_+$.
To see why, let $L' \in C^0(\R;\frakL)$ be another curve with $L(t) \to L_\pm$ as $t\to\pm\infty$,
convergence in the uniform operator topology on $\calL(X^1,X^0)$.
Then for each $t \in \R$
\[
L(t) - L'(t) =
\begin{pmatrix}
  0 & 0 \\
  L_1(t) - L'_1(t) & L_2(t) - L'_2(t)
\end{pmatrix},
\]
note that this is a bounded operator from $X^0$ to $X^0$ and therefore compact operator from $X^1$ to $X^0$.
Moreover, $L - L' \in C^0(\R, \calL(X^1,X^0))$, and $L(t) - L'(t) \to 0$ as $t \to \pm\infty$.
Hence the induced multiplication operator $L - L'$ is compact from $W^{1,p}(\R;X^0,X^1)$ to $L^2(\R;X^0)$, see \cite{rabier2004robbin}.
Consequently, $\ind( \calD_L ) = \ind( \calD_{L'} )$.

Summarising these observations, we have the following theorem.
\begin{mythm}
Let $L \in C^0(\R;\frakL)$ be such that $L(t) \to L_\pm$ as $t \to \pm\infty$ in the uniform operator topology on $\calL(X^1,X^0)$, where $L_\pm \in \frakL_{\text{hyp}}$.
  Then the map $\calD_L$ is Fredholm from $W^{1,2}(\R;X^0,X^1)$ to $L^2(\R;X^0)$, and its index only depends on the endpoints $L_-$, $L_+$.
\end{mythm}
This allows us to define a relative index:
\begin{align*}
  & \nu : \frakL_{\text{hyp}} \times \frakL_{\text{hyp}} \to \Z, \\
  & \nu(L_-,L_+) = \ind( \calD_L ).
\end{align*}

\subsection{Transitivity}

We now prove lemma \ref{lemma:cyclicity_index}.
\begin{mylemma}
  Let $L_\alpha, L_\beta, L_\gamma \in \frakL_{\text{hyp}}$.
  Then
  \[
     \begin{array}{l l}
      \nu(L_\alpha,L_\beta) = - \nu(L_\beta,L_\alpha) & \quad \text{(antisymmetry)}, \\
      \nu(L_\alpha , L_\gamma) = \nu(L_\alpha, L_\beta) + \nu(L_\beta, L_\gamma) & \quad \text{(cyclicity)}.
    \end{array}
  \]
\end{mylemma}
\begin{myproof}
This proof is an adaption of the argument given in \cite{robbin1995spectral}.
First we note that the antisymmetry follows from the cyclicity.
Indeed, since $\sigma(L_\alpha) \cap i\R = \emptyset$ and lemma \ref{lemma:resolvent_decay} is applicable, the operator
\[
\calD_{L_\alpha} = \partial_t + L_\alpha : W^{1,2}(\R;X^0,X^1) \to L^2(\R;X^0)
\]
is invertible, see \cite{rabier2003isomorphism} for details.
Consequently,
\[
\nu(L_\alpha,L_\beta) + \nu(L_\beta,L_\alpha) = \nu(L_\alpha,L_\alpha) = \ind(\calD_{L_\alpha}) = 0.
\]

To prove the cyclicity we first choose paths $L_{\alpha\beta}, L_{\beta\gamma} \in C^0(\R;\frakL)$ such that 
$L_{\alpha\beta}(t) = L_\alpha$ for $t \leq -1$, $L_{\alpha\beta}(t) = L_\beta$ for $t \geq 1$,
$L_{\beta\gamma}(t) = L_\beta$ for $t \leq -1$, $L_{\beta\gamma}(t) = L_\gamma$ for $t \geq 1$.
Moreover, given $T \geq 0$ let $L_{\alpha\gamma} \in C^0(\R;\frakL)$ be defined by
\[
L_{\alpha\gamma}(t) =
\begin{cases}
  L_{\alpha\beta}(t+T+1) & \text{for } t \leq 0, \\
  L_{\beta\gamma}(t-T-1) & \text{for } t \geq 0.
\end{cases}
\]

Consider the operators
\begin{align*}
  & M = \partial_t +
  \begin{pmatrix}
    L_{\alpha\beta} & 0 \\
    0 & L_{\beta\gamma}
  \end{pmatrix}, \\
  & N = \partial_t +
  \begin{pmatrix}
    L_{\alpha\gamma} & 0 \\
    0 & L_\beta
  \end{pmatrix}.
\end{align*}
These are bounded Fredholm operators from 
\[
\calX := W^{1,2}(\R;X^0,X^1) \times W^{1,2}(\R;X^0,X^1)
\]
to
\[
\calY := L^2(\R;X^0) \times L^2(\R;X^0).
\]
We have $\ind(M) = \nu(L_\alpha,L_\beta) + \nu(L_\beta,L_\gamma)$ and $\ind(N) = \nu(L_\alpha,L_\gamma) + \nu(L_\beta,L_\beta) = \nu(L_\alpha,L_\gamma)$.
Hence we need to prove that $\ind(M) = \ind(N)$.

Let $\eta \in C^\infty(\R)$ be such that $\eta(t) = 0$ for $t \leq -1$, and $\eta(t) = \pi/2$ for $t \geq 1$, and define
\[
R(t) =
\begin{pmatrix}
  \cos( \eta(t/T) ) & \sin( \eta(t/T) ) \\
  -\sin( \eta(t/T) ) & \cos( \eta(t/T) )
\end{pmatrix}.
\]
Then $R$ induces automorphisms of both $\calX$ and $\calY$.
Hence we can conjugate $N$ with $R$ without changing the Fredholm index.
Computing the conjugate yields
\[
 ( R^{-1} N R )(t) = \partial_t + I(t) + J(t) + K(t)
\]
where
\begin{gather*}
  I(t) = \frac{\eta'(t)}{T}
\begin{pmatrix}
  0 & -1 \\
  1 & 0
\end{pmatrix}, \\
  J(t) = 
\begin{pmatrix}
  \cos^2(\eta(t/T)) L_{\alpha\gamma}(t) + \sin^2(\eta(t/T)) L_\beta & 0 \\
  0 & \cos^2(\eta(t/T)) L_\beta + \sin^2(\eta(t/T)) L_{\alpha\gamma}(t)
\end{pmatrix}, \\
  K(t) =
\cos(\eta(t/T)) \sin(\eta(t/T))
\begin{pmatrix}
  0 & L_{\alpha\gamma}(t) - L_\beta \\
  L_{\alpha\gamma}(t) - L_\beta & 0
\end{pmatrix}.
\end{gather*}
Since $I \to 0$ in the uniform operator topology on $\calL(\calX,\calY)$ as $T \to \infty$, 
continuity of the Fredholm index implies that after choosing $T$ sufficiently large, $\ind(N) = \ind( R^{-1} N R ) = \ind(\partial_t + J + K)$.
Observe that $L_{\alpha\gamma}(t) = L_{\alpha\beta}(t+T+1)$ for $-\infty < t \leq T$, so that
\[
\cos^2( \eta(t/T) ) L_{\alpha\gamma}(t) + \sin^2( \eta(t/T) ) L_\beta(t) = L_{\alpha\beta}(t+T+1) + \sin^2( \eta(t/T) ) \big( L_\beta - L_{\alpha\beta}(t+T+1) \big)
\]
for all $t \in \R$.
Since $L_{\alpha\beta}(t+T+1) = L_\beta$ for $t \geq -T$, in fact
\[
\cos^2( \eta(t/T) ) L_{\alpha\gamma}(t) + \sin^2( \eta(t/T) ) L_\beta(t) = L_{\alpha\beta}(t+T+1).
\]
Using a similar computation for the other nonzero entry in $J(t)$ one sees that
\[
J(t) =
\begin{pmatrix}
  L_{\alpha\beta}(t+T+1) & 0 \\
  0 & L_{\beta\gamma}(t-T-1)
\end{pmatrix}.
\]
Similarly one can verify that $K(t) = 0$ for all $t \in [-T,T]$ (and hence $K(t) = 0$ for all $t \in \R$).

Let $S \in \calL(\calY)$ be the shift operator $S(U,V)(t) = ( U(t+T+1) , V(t-T-1) )$.
Note that $S$ is an automorphism of $\calY$ and restricts to an automorphism of $\calX$.
So we can let $S$ act on $\calL(\calX,\calY)$ via conjugation.
Note that $S$ commutes with $\partial_t$, and
\[
S^{-1} J S = 
\begin{pmatrix}
  L_{\alpha\beta} & 0 \\
  0 & L_{\beta\gamma}
\end{pmatrix},
\]
so that $S^{-1} ( \partial_t + J ) S = \partial_t + S^{-1} J S = M$.
Hence
\[
\ind(N) = \ind( S^{-1} R^{-1} N R S ) = \ind( S^{-1} ( \partial_t + J ) S ) = \ind(M),
\]
thus concluding the proof.
\end{myproof}

\section{Exponential dichotomy along heteroclinic orbits}
\label{sec:appendix_exp_decay}

Here we give some details as to why the linearization of \eqref{eq:TWE} along heteroclinic orbits possesses an exponential dichotomy.
We begin with citing a simplified version of the main theorem from \cite{peterhof1997exponential}.
\begin{mythm}[\cite{peterhof1997exponential}]
\label{thm:exp_dich}
  Let $X^0$ be a reflexive Banach space, and $L : \calD(L) \to X^0$ a closed, 
possibly unbounded operator such that $X^1 := \calD(L)$ is dense in $X^0$.
Let $X^1$ be equipped with the graph norm of $L$, i.e., $\| u \|_{X^1} = ( \| u \|_{X^0}^2 + \| L u \|_{X^0}^2 )^{1/2}$.
Let $J = [\tau_0,\infty)$ and suppose that $B \in C^0(J,\calL(X^0))$ is Lipschitz continuous.
Consider the abstract differential equation
\begin{equation}
  \label{eq:19}
  \partial_t W(t) + ( L + B(t) ) W(t) = 0, \qquad W \in C^0(  J ; X^1 ) \cap C^1( \inter J ; X^0 , X^1 ).
\end{equation}
Assume that the following four conditions are satisfied.
\begin{enumerate}
\item 
There exists a constant $C$ such that
\begin{equation}
  \label{eq:resolvent_estimate}
  \| ( L - i\mu )^{-1} \|_{\calL(X^0)} \leq \frac{ C }{ 1 + |\mu| }
\end{equation}
for all $\mu \in \R$.
Assume that there is a projection $P_- \in \calL(X^0)$ such that $L^{-1}$ and $P_-$ commute.
Furthermore, assume there exists a constant $\delta > 0$ such that $\real \lambda < -\delta$ for $\lambda \in \sigma( L P_- )$
and $\real \lambda > \delta$ for $\lambda \in \sigma( L ( 1 - P_- ) )$.

\item 
It holds that $\| B(t) \|_{\calL(X^0)} \to 0$ as $t \to \infty$.

\item
The operator $L$ has compact resolvent.

\item
The only solution $W$ of \eqref{eq:19} such that $\sup_{t\in J} \| W(t) \|_{X^0} < \infty$ and $W(0) = 0$ is the trivial solution $W \equiv 0$.
\end{enumerate}

Then \eqref{eq:19} has an exponential dichotomy in $X^0$ on the interval $J$ with rate $\gamma$, for any $0 \leq \gamma < \delta$.
In particular, there exists $K>0$ such that if $W$ is a solution of \eqref{eq:19}
with $\sup_{t\in J} \| W(t) \|_{X^0} < \infty$, it holds that
\[
\| W(t) \|_{X^0} \leq K e^{-\gamma |t - \tau|} \| W(\tau) \|_{X^0} \qquad \text{for} \quad t \geq \tau \geq \tau_0.
\]
\end{mythm}
We need to check that this theorem applies to our linearised equation.
Suppose hypotheses \Hgroup are satisfied.
Let $U$ be a solution of \eqref{eq:TWE} such that $U(t) \to Z$ in $X^0$ as $t \to \infty$, where $Z \in \calS$ is a hyperbolic rest point.
Define $J = [0,\infty)$.
We decompose
\[
\d A(U(t)) = L + B(t),
\]
where $L = \d A(Z_+) \in \calL(X^1,X^0)$, and
\[
B(t) =
\begin{pmatrix}
  0 & 0 \\
  f_u(x,u(t,x)) - f_u(x,z(x)) & 0
\end{pmatrix}
+
\begin{pmatrix}
  0 & 0 \\
  \partial_1 g(x,U(t)) & \partial_2 g(x,U(t))
\end{pmatrix}.
\]
Here $\partial_1 g(x,(u,v)) := \ddi{ g(x,(u,v)) }{ u }$ and $\partial_2 g(x,(u,v)) := \ddi{ g(x,(u,v)) }{ v }$.

Let us now construct the projections needed in condition (1) of theorem \ref{thm:exp_dich}.
First, let $\{ \mu_n \}_n$ be the eigenvalues (counting multiplicity) of $\Delta + f_u(x,z(x))$, arranged in decreasing order.
Let $k_0$ be such that $\mu_{k_0+1} < 0 < \mu_{k_0}$, let $k_c$ be such that $\mu_{k_c + 1} < c^2 / 4 \leq \mu_{k_c}$,
and let $k_{\text{def}}$ denote the number of eigenvalues which are equal to $c^2 / 4$.

Let $( \phi_n )_n$ be an orthonormal basis for $H_B^1(\Omega)$ consisting of eigenfunctions for $\Delta + f_u(\cdot,z)$, 
arranged so that $\phi_n$ is an eigenfunction corresponding to the eigenvalue $\mu_n$.
Then $L$ has eigenvalues
\[
\lambda^\pm_n = - \frac c 2 \pm \frac 1 2 \sqrt{ c^2 - 4 \mu_n },
\]
with corresponding eigenfunctions $\Psi^\pm_n$ given by
\[
\Psi^\pm_n =
\begin{pmatrix}
  \phi_n \\
  - \lambda^\pm_n \phi_n
\end{pmatrix}.
\]
A direct computation shows that these eigenfunctions are orthogonal in $X^0$.
If $\mu_n = c^2 / 4$ (i.e.\  when $k_{\text{def}} \geq 1$ and $k_c - k_{\text{def}} \leq n \leq k_c$),
then $\lambda_n^\pm = - c /2$ is a defective eigenvalue of $L$ and a corresponding generalized eigenfunction is given by
\[
\widetilde \Psi_n =
\begin{pmatrix}
  0 \\
  \phi_n
\end{pmatrix}.
\]
Also, note that for $1 \leq n \leq k_c$, we have
\[
\real \Psi^\pm_n =
\begin{pmatrix}
  \phi_n \\
  \frac c 2 \phi_n
\end{pmatrix}
\qquad \text{and} \qquad
\imag \Psi^\pm_n =
\begin{pmatrix}
  0 \\
  \mp \frac 1 2 \sqrt{ 4 \mu_n - c^2 } \phi_n
\end{pmatrix}.
\]
Define
\begin{equation*}
  \begin{split}
    E_- := \vspan \bigg( \set{ \real \Psi^-_n , \imag \Psi^-_n }{ n \geq 1 } &\cup \set{ \real \Psi^+_n , \imag \Psi^+_n }{ 1 \leq n \leq k_0 } \\
&\cup \set{ \widetilde \Psi_n }{ \mu_n = c^2 / 4 } \bigg)
  \end{split}
\end{equation*}
and
\[
E_+ := \vspan\set{ \Psi^+_n }{ n \geq k_0 + 1 }.
\]
Denote by $\cl{E_-}$ and $\cl{E_+}$ the closure of $E_-$ and $E_+$ in $X^0$.

We claim that $X^0 = \cl{E_-} \oplus \cl{E_+}$.
First note that if $(a_n)_n \in \ell^2(\N)$, then $\sum_{n=k_c+1}^\infty a_n \Psi^+_n$ converges in $X^0$.
Indeed,
\begin{equation}
  \label{eq:57}
  \left\| \sum_{n=i}^j a_n \Psi^+_n \right\|_{X^0}^2 
= \left\| \sum_{n=i}^j a_n \phi_n \right\|_{H_B^1(\Omega)}^2 + \left\| \sum_{n=i}^j a_n \lambda^+_n \phi_n \right\|_{L^2(\Omega)}^2
= \sum_{n=i}^j \big( 1 +  |\lambda^+_n|^2 \|\phi_n\|_{L^2(\Omega)}^2 \big)|a_n|^2,
\end{equation}
and we need to check that the right hand side tends to zero as $i,j\to\infty$.
Note that, since $\| \phi_n \|_{H^1_B(\Omega)} = 1$,
\[
-\mu_n \| \phi_n \|_{L^2(\Omega)}^2 = \| \nabla \phi_n \|_{L^2(\Omega)}^2 - \langle f_u(\cdot,z) \phi_n , \phi_n \rangle_{L^2(\Omega)} 
\leq 1 + \|f_u(\cdot,z)\|_{L^\infty(\Omega)},
\]
and since $-\mu_n \to \infty$ as $n\to\infty$ it follows that $\|\phi_n\|_{L^2(\Omega)} \to 0$ as $n\to\infty$.
Then, since
\begin{align*}
 - \mu_n \| \phi_n \|_{L^2(\Omega)}^2 &= \|\phi_n\|_{H^1_B(\Omega)}^2 - \|\phi_n\|_{L^2(\Omega)}^2 - \langle f_u(\cdot,z) \phi_n , \phi_n \rangle_{L^2(\Omega)} \\
&= 1 - \langle ( 1 + f_u(\cdot,z) ) \phi_n , \phi_n \rangle_{L^2(\Omega)},
\end{align*}
we see that $-\mu_n \| \phi_n \|_{L^2(\Omega)}^2 \sim 1$ as $n\to\infty$ (here ``$\sim$'' denotes asymptotic equivalence of sequences).
Since $|\lambda_n^+|^2 \sim -\mu_n$ as $n\to\infty$, it now follows that the right hand side in \eqref{eq:57} tends to $0$ as $i , j \to \infty$.
An identical computation shows that  $(a_n)_n \in \ell^2(\N)$, then $\sum_{n=k_c+1}^\infty a_n \Psi^-_n$ converges in $X^0$.
Thus $\cl{E_-} + \cl{E_+}$ contains elements of the form
\begin{multline*}
  \begin{pmatrix}
  x \\ y
\end{pmatrix}
=
\sum_{n=1}^{k_c} a_n
\begin{pmatrix}
  \phi_n \\
  \frac c 2 \phi_n
\end{pmatrix}
+
\sum_{n=1}^{k_c-k_{\text{def}}} b_n
\begin{pmatrix}
  0 \\
  \frac 1 2 \sqrt{ 4 \mu_n - c^2 } \phi_n
\end{pmatrix}
+
\sum_{n=k_c - k_{\text{def}}+1}^{k_c}
b_n
\begin{pmatrix}
  0 \\
  \phi_n
\end{pmatrix}
\\
+
\sum_{n=k_c+1}^\infty a_n
\begin{pmatrix}
  \phi_n \\ \lambda^+_n \phi_n
\end{pmatrix}
+
\sum_{n=k_c+1}^\infty b_n
\begin{pmatrix}
  \phi_n \\ \lambda^-_n \phi_n
\end{pmatrix},
\end{multline*}
where $(a_n)_n , (b_n)_n \in \ell^2(\N)$.
To see why any $(x,y) \in X^0$ is of this form, write $x = \sum_{n=1}^\infty c_n \phi_n$ and $y = \sum_{n=1}^\infty d_n \phi_n$,
where $(c_n)_n \in \ell^2(\N)$ and $( \|\phi_n\|_{L^2(\Omega)} d_n )_n \in \ell^2(\N)$.
Set
\begin{equation*}
  \left\{
    \begin{array}{l l}
      \left\{
      \begin{array}{l}
        a_n = c_n \\
        b_n = - \dfrac{c}{\sqrt{4 \mu_n - c^2}} c_n + \dfrac{2}{\sqrt{4 \mu_n - c^2}} d_n
      \end{array}
      \right. & n \leq k_c - k_{\text{def}}, \\ \\
      \left\{
      \begin{array}{l}
        a_n = c_n \\
        b_n = - \dfrac c 2 c_n + d_n
      \end{array}
      \right. & k_c - k_{\text{def}} + 1 \leq n \leq k_c, \\ \\
      \left\{
      \begin{array}{l}
        a_n = \dfrac{ \lambda_n^- }{ \lambda_n^- - \lambda_n^+ } c_n - \dfrac{ 1 }{ \lambda_n^- - \lambda_n^+ } d_n \\
        b_n = - \dfrac{ \lambda_n^+ }{ \lambda_n^- - \lambda_n^+ } c_n + \dfrac{ 1 }{ \lambda_n^- - \lambda_n^+ } d_n
      \end{array}
      \right. & n \geq k_c + 1.
    \end{array}
  \right.
\end{equation*}
Note that, as $n\to\infty$,
\[
\left| \frac{ \lambda_n^\pm }{ \lambda_n^- - \lambda_n^+ } \right|^2 \sim \frac 1 4 \qquad \text{and} \qquad
\left| \frac{ 1 }{ \lambda_n^- - \lambda_n^+ } \right|^2 \sim - \frac{1}{4 \mu_n} \sim \frac 1 4 \|\phi_n\|_{L^2(\Omega)}^2.
\]
Hence $(a_n)_n , (b_n)_n \in \ell^2(\N)$, 
and since $a_n \phi_n + b_n \phi_n = c_n \phi_n$ and $a_n \lambda_n^+ \phi_n + b_n \lambda_n^- \phi_n = d_n \phi_n$,
this proves that $(x,y) \in \cl{E_-} + \cl{E_+}$.
Since $E_-$ and $E_+$ are orthogonal in $X^0$, it follows that $\cl{E_-} \cap \cl{E_+} = \{0\}$.
Thus $X^0 = \cl{E_-} \oplus \cl{E_+}$.

Let $P_-$ be the projection onto $\cl{E_-}$ along $\cl{E_+}$.
Then $P_-$ commutes with $L^{-1}$.
The construction of the sets $\cl{E_\pm}$ ensures that
$\real \lambda < -\delta$ for $\lambda \in \sigma(AP_-)$
and $\real \lambda > \delta$ for $\lambda \in \sigma( A(1-P_-) )$.
Finally, estimate \eqref{eq:resolvent_estimate} is a special case of the result from lemma \ref{lemma:resolvent_decay}.
Hence condition (1) of theorem \ref{thm:exp_dich} is satisfied.

Note that $B \in C^1(J,\calL(X^0))$ and $\| B(t) \|_{\calL(X^0)} \to 0$ as $t \to \infty$, hence $B$ is Lipschitz continuous.
Thus condition (2) of theorem \ref{thm:exp_dich} is also satisfied.
Since the embedding $X^1 \hookrightarrow X^0$ is compact, $L$ has compact resolvent.
Hence condition (3) of theorem \ref{thm:exp_dich} is satisfies.

To ensure condition (4) of theorem \ref{thm:exp_dich} is satisfied we need to assume that the nonlinearity $f$ is of class $C^4$.
Recalling from section \ref{sec:compactness} that $U \in C_b^4(\R;X^0,\dots,X^3)$, we see that 
\[
B \in C^3(\inter J,\calL(X^0)) \cap C^2(\inter J,\calL(X^1)) \cap C^1(\inter J,\calL(X^2)) \cap C^0(J,\calL(X^3)).
\]
In turn elliptic regularity theory implies that $w \in C^4(\inter J,L^2(\Omega))$ (where $W = (w,\partial_t w)$).
Together with the mean value theorem this ensures that if $W(0) = 0$, then $w$ satisfies the decay estimates \eqref{eq:decay_conditions} 
from lemma \ref{lemma:ineq_continuation} around $t = 0$.
Hence by lemma \ref{lemma:ineq_continuation}, $w(t) = 0$ and hence $W(t) = 0$ for $t$ in a neighbourhood of $0$.
By an argument similar to the proof of theorem \ref{thm:unique_continuation} one then finds that $W \equiv 0$.

The preceding discussion shows that the linearized equation posesses an exponential dichotomy on $X^0$ 
with some rate $\gamma > 0$ on the time interval $J = [\tau_0,\infty)$, provided that the nonlinearity $f$ is of class $C^m$ with $m \geq 4$.


\bibliography{main}
\bibliographystyle{abbrv}

\end{document}